\documentclass[twoside,11pt]{article}
\usepackage{jmlr2e}
\usepackage{wan}
\usepackage{graphicx} % Required for inserting images
\usepackage{tabularx}
\usepackage{array}
\usepackage{algorithm}
\usepackage{algorithmic}
\usepackage{subcaption}
\usepackage{xpatch,xcolor}
\usepackage{rotating}
\usepackage{doi}

\usepackage{lastpage}

\firstpageno{1}

\begin{document}

\title{High-Dimensional Model Averaging via Cross-Validation}

\author{\name Zhengyan Wan \email zhengyanwan066@gmail.com \\
       \addr School of Statistics\\
       East China Normal University\\
       Shanghai, 200062, China
       \AND
       \name Fang Fang\thanks{Corresponding author} \email ffang@sfs.ecnu.edu.cn \\
       \addr KLATASDS-MOE, School of Statistics\\
       East China Normal University\\
       Shanghai, 200062, China
       \AND
       \name Binyan Jiang\email by.jiang@polyu.edu.hk\\
       \addr Department of Data Science and Artificial Intelligence\\
       The Hong Kong Polytechnic University\\
       Hung Hom, Kowloon, Hong Kong}

\maketitle

\begin{abstract}%   
 Model averaging is an important alternative to model selection with attractive prediction accuracy. However, its application to high-dimensional data remains under-explored. We propose a high-dimensional model averaging method via cross-validation under a general framework and systematically establish its theoretical properties. Each candidate model is fitted using a flexible loss function paired with a general regularizer, and the optimal weights are determined by minimizing a cross-validation criterion.    When all candidate models are misspecified, we establish a non-asymptotic upper bound and a minimax lower bound for our weight estimator. The asymptotic optimality is also derived, showing that the proposed weight estimator achieves the lowest possible prediction risk asymptotically. When the correct models are included in the candidate model set, the proposed method asymptotically assigns all weights to the correct models, and the model averaging estimator achieves a nearly-oracle convergence rate. Further, we introduce a post-averaging debiased estimator and establish Gaussian and bootstrap approximation to construct simultaneous confidence intervals. A fast greedy model averaging (FGMA) algorithm is proposed to solve the simplex-constrained optimization problem, which has a descent property empirically and a faster convergence rate, compared to the original greedy model averaging algorithm. Empirical results demonstrate the strong competitiveness of the proposed method in prediction and inference, compared to other existing model averaging and selection methods.
\end{abstract}


\vspace{2cm}
\section{Introduction}

Model averaging is a widely used strategy for addressing model uncertainty, which integrates all candidate models into a weighted estimator, reducing uncertainty in model selection \citep{draper1995assessment,yuan2005combining} and achieving lower prediction error due to its asymptotic optimality, as established by the seminal work of
\cite{hansen2007least} and many others.
Model averaging comprises two main approaches: Bayesian model averaging (BMA) and frequentist model averaging (FMA). BMA computes a weighted average of posterior distributions across candidate models, with weights determined by their posterior model probabilities \citep{hoeting1999bayesian}. Conversely, FMA employs frequentist methods to optimize weight selection for specific criteria \citep[e.g.,][]{yang2001adaptive,leung2006information,hansen2007least,hansen2012jackknife,ando2014model,zhang2020parsimonious,peng2024optimality}.

Over the past decade, some studies have developed model averaging methods incorporating screening steps in high-dimensional settings. \cite{ando2014model} grouped covariates to prepare candidate linear models following a screening step, and selected model weights by minimizing a cross-validation criterion. This approach has been extended to various models. Examples include \cite{ando2017weight}, \cite{xie2021model}, \cite{yan2021optimal} and \cite{wang2023jackknife}. \cite{zhang2020parsimonious} constructed nested candidate linear models using the covariates from the adaptive Lasso solution path, selecting weights by minimizing a BIC-type criterion. These studies establish asymptotic optimality in the sense that the selected weights can yield the lowest in-sample prediction risk among all weight choices asymptotically. However, these asymptotic optimality results require the number of candidate models to be smaller than the sample size. In addition, the screening step in these studies may remove some truly important covariates especially when the signals are non-sparse. Alternatively, for high-dimensional survival data, \cite{he2020functional} and \cite{he2023rank} ranked the importance of each covariate and partitioned all the covariates into disjoint groups; each group is used for constructing a candidate model. \cite{he2023rank} established an upper bound for the estimation error for the weight estimator, allowing the number of candidate models to grow exponentially with the sample size. However, this upper bound depends on the infinity norm of the gradient function, which might not converge to zero when the number of candidate models is fixed. Moreover, the theoretical results they developed require the optimal weight to be unique, which is a relatively strong assumption when the candidate model set includes more than one correct model. For ultrahigh-dimensional longitudinal data, \cite{jiang2024robust} used forward selection to build nested candidate models and minimized a penalized weighted minimum absolute deviation criterion to select weights. The results of weight consistency and risk consistency are developed in their work under the boundedness assumption of the covariates. Overall speaking, the model averaging with high-dimensional data is far from fully understood. 

In this paper, we develop a unified high-dimensional model averaging framework. Different from the above developments, in our work, we divide the covariates into two parts, a nested part and a nonnested part, after ranking the importance of covariates. This procedure reduces more approximation error than using only a low-dimensional nested structure or a candidate model set with non-overlapping groups, which is more robust against model underfitting and coefficient non-sparsity. In addition, for each candidate model, our approach allows us to estimate candidate models by penalized empirical risk minimization with a general regularizer, whereas most of the existing papers require candidate models to be fitted using unpenalized methods. We select model weights by minimizing a $J$-fold cross-validation criterion. Unlike \cite{zhang2023model}, in our framework, we allow the use of general loss functions rather than just squared loss. For computation, we develop a fast greedy model averaging (FGMA) algorithm to solve this optimization problem under a simplex constraint, which empirically has a descent property. The convergence analysis for the proposed FGMA algorithm is derived, which matches the optimal rate of first-order algorithms \citep{nesterov2018lectures}. In contrast, the original greedy model averaging algorithms introduced in \cite{dai2012deviation} and \cite{he2020functional,he2023rank} do not exhibit monotonic decrease with iterations and have a slower convergence rate than our proposed algorithm. 

Theoretically, we develop non-asymptotic bounds for our proposed cross-validation weight estimator under weaker conditions on the covariates. When all models are wrong, we provide a minimax lower bound for the weight estimator, which nearly matches the upper bound we developed up to the square root of a compatibility factor. To the best of our knowledge, this is the first result regarding the lower bound of weight estimation in the optimal model averaging literature. When there are multiple correct models, we demonstrate that the proposed weight estimator assigns all weights to the correctly specified models asymptotically. In this scenario, the optimal weight is not unique, which requires additional efforts to make theoretical analysis. In both scenarios, the optimality results or the risk consistency for the proposed model averaging estimator are derived in terms of out-of-sample prediction risk. These theories are established under a general framework, allowing the number of candidate models to diverge exponentially with the sample size and to be greater than the sample size.

Beyond estimation and prediction, statistical inference is crucial in many applications. However, model averaging has been long term criticized for its limited development of robust inference methods.
 Some efforts have been made in a local misspecification framework, for example, \cite{liu2015distribution}. Under the more realistic framework of fixed parameters, \cite{zhang2019inference} and \cite{yu2024post} investigated the asymptotic distribution of model averaging estimators and post-averaging inference for linear and generalized linear models under the scenario where correct models are included in the candidate model set. The authors proposed to construct confidence intervals using a simulation-based method. Similar theoretical results have also been established in \cite{Racine2023optimal}, \cite{tu2025quantile}, and \cite{wong2025valid}, among other works, in a case-by-case manner. For high-dimensional models, however, there is still no theory or methodology for post-averaging inference. In our paper, we consider the post-averaging simultaneous inference problem and propose a one-step debiased estimator, inspired by the debiasing technique developed by many authors \citep[see, for example,][]{zhang2014confidence,ning2017general}. Following its Bahadur representation, we establish the high-dimensional Gaussian and bootstrap approximation results for this estimator based on the high-dimensional central limit theorem and high-dimensional multiplier bootstrap \citep{chernozhuokov2022improved}. These results allow us to construct simultaneous confidence intervals using the proposed model averaging estimator using the multiplier bootstrap technique.

To summarize, we list our main contributions as follows.
\begin{enumerate}
  %  \item We develop a high-dimensional model averaging framework with a general loss function. To handle high-dimensional covariates, our approach allows us to estimate candidate models by penalized empirical risk minimization with a general regularizer.
    \item We propose a general model averaging framework for high-dimensional data, employing a flexible loss function and an efficient strategy for selecting candidate models based on covariate importance ranking. This approach facilitates the estimation of candidate models via penalized empirical risk minimization with a general regularizer.

    \item We establish novel theoretical results under two scenarios: (i) all candidate models are misspecified, and (ii) at least one correct model is included in the candidate set.   These results hold under general assumptions on the form of penalized loss functions and accommodate settings where the number of candidate models exceeds the sample size. 
    \begin{itemize}
    \item
    When all candidate models are misspecified, we establish an oracle inequality for the $\ell_1$ estimation error of the proposed weight estimator and develop the asymptotic optimality in terms of out-of-sample prediction risk. Additionally, we establish a minimax lower bound for the weight estimator, making the first such result in the model averaging literature.

  %  \item We investigate the post-averaging simultaneous inference problem and develop a multiplier bootstrap procedure based on the Gaussian and bootstrap approximation results.
  \item   When there exists at least one correctly specified model, we show that the weight estimator asymptotically assigns all weights to the correct models. The convergence rate of the proposed model averaging estimator achieves a nearly-optimal rate. Moreover, risk consistency of the proposed model averaging method is established. To enable simultaneous inference, we introduce a one-step debiased estimator, deriving its Bahadur representation and establishing Gaussian and bootstrap approximation results for the Kolmogorov-Smirnov distance. These results facilitate the construction of simultaneous confidence intervals using the multiplier bootstrap technique.

\end{itemize}

 %   \item We propose a fast greedy model averaging algorithm to minimize the cross-validation criterion under a simplex constraint and derive its convergence analysis, which matches the optimal rate of first-order algorithms.

  \item   We develop a fast greedy model averaging algorithm to minimize the cross-validation criterion under a simplex constraint, with a convergence analysis demonstrating an optimal rate matching that of first-order algorithms. Empirically, our algorithm exhibits a descent property and achieves faster convergence than other widely used greedy model averaging algorithms \citep{dai2012deviation,he2020functional,he2023rank}.

\end{enumerate}

% The remainder of the paper is organized as follows. \cref{sec: set-up} introduces the high-dimensional model averaging framework using a cross-validation criterion and the preparation of candidate models. \cref{sec: theoretical result} presents some theoretical guarantees on our proposed estimator. For the scenario where all candidate models are misspecified, we state non-asymptotic bounds on the estimation error of our proposed weight estimator, and establish the asymptotic optimality of the weight estimator. For the scenario where there is at least one correctly specified  candidate model, the consistency of the weight estimator is provided, and Gaussian and bootstrap approximation results for the post-averaging debiased estimator are derived. In \cref{sec:implementation}, we introduce the fast greedy model averaging algorithm and demonstrate the implementation details of our methods. \cref{sec: simulation} investigates the finite sample performance of the proposed methods via simulations. \cref{sec: real data} applies our approach to analyze a genomic dataset. \cref{sec: concluding remarks} provides the concluding remark. Technical proofs and some further discussion are presented in the Appendix. All the implementation codes for our proposed method including estimation, prediction and inference are available at: \url{https://github.com/WanZhengyan/HDMA}.

The remainder of the paper is organized as follows. \cref{sec: set-up} introduces the high-dimensional model averaging framework using a $J$-fold cross-validation criterion and the preparation of candidate models. \cref{sec: theoretical result} provides the theoretical guarantees of our proposed methods. In \cref{sec:implementation}, we introduce the fast greedy model averaging algorithm and demonstrate the implementation details of our methods. \cref{sec: simulation} investigates the finite sample performance of the proposed methods via simulations. \cref{sec: real data} applies our approach to analyze a genomic dataset. \cref{sec: concluding remarks} provides the concluding remark. Technical proofs and some further discussion are presented in the Appendix. All the implementation codes for our proposed method including estimation, prediction and inference are available at: \url{https://github.com/WanZhengyan/HDMA}.

\vspace{0.2cm}
{\sc Notation.} We write $a_n\lesssim b_n$ if $a_n=O(b_n)$, and $a_n\gtrsim b_n$ if $b_n=O(a_n)$. We say $a_n\asymp b_n$ if $a_n\lesssim b_n$ and $a_n\gtrsim b_n$. For an integer $N$, we denote $\brac{N}=\set{1,2,\dots,N}$. We write $a\wedge b=\min(a,b)$ and $a\vee b=\max(a,b)$. Let $\lfloor a\rfloor$ be the floor of a real number $a$. Write $\norm{\beta}_q=(\sum_{j=1}^p\abs{\beta_j}^q)^{1/q}$ for some vector $\beta\in\R^p$, where $q\in(0,\infty).$ The sub-gaussian norm of a sub-gaussian random variable $u\in \R$ is $\norm{u}_{\psi_2}=\inf \{t>0|\E \exp(u^2/t^2)\leq 2\}$ and the sub-gaussian norm of a sub-gaussian random vector $V\in\R^p$ is $\norm{V}_{\psi_2}=\sup_{v\in\R^p:\ltwo{v}=1}\norm{\<V,v\>}_{\psi_2}$. Similarly, the sub-exponential norm of a sub-exponential random variable $u\in \R$ is $\norm{u}_{\psi_1}=\inf \{t>0|\E \exp(\abs{u}/t)\leq 2\}$ and the sub-exponential norm of a sub-exponential random vector $V\in\R^p$ is $\norm{V}_{\psi_1}=\sup_{v\in\R^p:\ltwo{v}=1}\norm{\<V,v\>}_{\psi_1}$. For a matrix $H\in\R^{m\times m}$, $\lambda_{\min}(H)$ and $\lambda_{\max}(H)$ stand for the minimum and maximum eigenvalues of the matrix $H$, respectively. Moreover, for any $H\in\R^{m\times n}$, $\matrixnorm{H}_p=\sup_{x\in\R^n:\norm{x}_p=1}\norm{Hx}_p$ and $\matrixnorm{H}_{\max}=\max_{1\leq i\leq m;1\leq j\leq n}\abs{H_{i,j}}$. Throughout this paper, some positive constants $C, C_1, C_2, C_3, C_4, c, c_1, c_2$ which are not related with $n$, are allowed to change from line to line.

\section{High-Dimensional Model Averaging via Cross-Validation}\label{sec: set-up}
Consider $n$ independent and identically distributed random variables $\set{(Y_i,\mathbf{x}_i)}_{i=1}^n$ where $\mathbf{x}_i=(\mathrm{x}_{i1},\mathrm{x}_{i2},\dots)^\top$ is an infinite-dimensional vector and $Y_i\in\mc{Y}$ is a response variable. We assume that the data are generated by
\begin{align}
    Y_i=g\parenBig{\sum_{j=1}^\infty \beta_j\mathrm{x}_{ij}, \eps_i},~1\leq i\leq n,\label{model}
\end{align}
where $\bm{\beta}=(\beta_1,\beta_2,\dots)^\top$ is a vector of unknown coefficients, $\eps_i$ is a random error independent of $\mathbf{x}_i$, and $g(\cdot)$ is a known function. {Many common models, such as the linear model and the logistic model, are included in \eqref{model}.} As mentioned in \cite{hansen2007least}, the vector $\mathbf{x}_i$ can include transformations of the original covariates and model \eqref{model} also includes nonparametric models.  Let $L: \mc{Y}\times \R\to \R$ be a loss function which is convex in the second argument. We assume that the prediction risk $\E L(Y_1,\sum_{j=1}^\infty \theta_j\mathrm{x}_{1j})$ is minimized by $\bm{\beta}$ over the set $\set{\theta\in\R^\N|\sum_{j=1}^\infty\theta_j^2<\infty}$.

 In practice, we only have access to a $p$-dimensional covariate vector $X_i=(\mathrm{x}_{i1},\dots,\mathrm{x}_{ip})^\top\in\R^{p}$ and $p$ could be larger than the sample size $n$. Our goal is to use a linear function $\sum_{j=1}^p\hat{\beta}_j\mathrm{x}_{ij}$ to approximate $\sum_{j=1}^\infty\beta_j\mathrm{x}_{ij}$ with an estimator $\hat{\beta}=(\hat{\beta}_1,\dots,\hat{\beta}_p)^\top\in\R^p$ given the data $\mc{D}=\set{(Y_i,X_i)}_{i=1}^n$. In high-dimensional settings, we usually consider the estimator $\hat{\beta}$ based on penalized empirical risk minimization such as MCP, SCAD and Lasso, having the form
\begin{align*}
    \hat{\beta}=\argmin_{\beta\in\R^p}\dn\sumn L(Y_i,\langle \beta,X_i\rangle)+r_{\lambda_n}\paren{\beta},
\end{align*}
where $r_\lambda:\R^p\to\R$ is a regularizer and $\lambda_n$ is a tuning parameter. We assume that the regularizer $r_\lambda$ is decomposable across coordinates, and with a little abuse of notation, we write $r_\lambda(\beta)=\sum_{j=1}^p r_\lambda(\beta_j)$ for $\beta\in\R^p$. 

%%%%%%%%%%%%%%%%%%%%%%%%%%%%%

\subsection{Model Averaging Estimation}
We consider a sequence of approximating models $k=1,\dots,K$, where the $k$-th candidate model uses the covariates in $\mc{A}_k\subseteq\brac{p}$ with the index set size $p_k=\abs{\mc{A}_k}$ and $K$ can go to infinity as the sample size increases. We write the $k$-th approximating candidate model as
\begin{align*}
    Y_i=g\parenBig{\sum_{j\in\mc{A}_k}\beta_j\mathrm{x}_{ij},\eps_i},~ 1\leq i\leq n.
\end{align*}The estimator under the $k$-th candidate models with regularizer $r_\lambda$ is given by
\begin{align*}
    \hat{\beta}_{(k)}=\argmin_{\beta\in\R^p, \beta_{\mc{A}^c_k}=0}\dn\sumn L(Y_i,\langle \beta,X_i\rangle)+r_{\lambda_{n,k}}(\beta),~ k=1,
    \dots,K,
\end{align*}
where $\mc{A}^c_k=[p]\backslash\mc{A}_k$ and $\lambda_{n,k}$ is a tuning parameter. To aggregate $K$ candidate estimators $\set{\hat{\beta}_{(k)}}_{k=1}^K$, we would like to find a model averaging estimator which has the form $\hat{\beta}_{MA}(w)=\sum_{k=1}^Kw_k\hat{\beta}_{(k)}=\hat{B}w$, where $\hat{B}=(\hat{\beta}_{(1)},\dots,\hat{\beta}_{(K)})$ and the weight vector
\begin{align*}
    w\in\mc{W}^K\overset{\triangle}{=}\setBig{w\in[0,1]^K\Big|\sum_{k=1}^Kw_k=1}.
\end{align*}
The following procedure describes how to choose the weights by minimizing a cross-validation criterion.
\begin{enumerate}
    \item Split the data $\mc{D}=\set{(Y_i,X_i)}_{i=1}^n$ into $J$ ($J\ge2$) folds of equal size with the index sets $\set{\mc{I}_m}_{m=1}^J$, so that there are $n_1=n/J$ observations in each group.
    \item For $m=1,\dots,J$ and $k=1,\dots,K$, using the observations in $\mc{D}\backslash\set{(Y_i,X_i)}_{i\in\mc{I}_m}$ to calculate the estimator $\hat{\beta}_{(k)}^{\brac{-m}}$, which is defined by
    \begin{align*}
        \hat{\beta}_{(k)}^{\brac{-m}}=\argmin_{\beta\in\R^p, \beta_{\mc{A}^c_k}=0}\frac{1}{n_{-1}}\sum_{i\in\brac{n}\backslash\mc{I}_m} L(Y_i,\langle \beta,X_i\rangle)+r_{\lambda_{n,k}^{\brac{-m}}}(\beta),
    \end{align*}
    where $n_{-1}=n-n_1$ and $\lambda_{n,k}^{\brac{-m}}$ is a tuning parameter. 
    \item Define a $J$-fold cross-validation criterion
    \begin{align*}
        CV(w)=\sum_{m=1}^J\sum_{i\in\mc{I}_m}L\paren{Y_i,\<\hat{B}^{\brac{-m}}w,X_i\>},
    \end{align*}
    where $\hat{B}^{\brac{-m}}=\paren{\hat{\beta}_{(1)}^{\brac{-m}},\hat{\beta}_{(2)}^{\brac{-m}},\dots,\hat{\beta}_{(K)}^{\brac{-m}}}$.
    \item Choose the weights by minimizing the cross-validation criterion
    \begin{align}\label{eq:optimization problem}
        \hat{w}=\argmin_{w\in\mc{W}^K}CV(w),
    \end{align}
    and the model averaging estimator is given by $\hat{\beta}_{MA}(\hat{w})=\hat{B}\hat{w}$.
\end{enumerate}

Since the cross-validation criterion $CV(w)$ is a convex function of the weight vector, the cross-validation weights can be computed efficiently. In our work, we propose a fast greedy model averaging algorithm to solve this simplex-constrained optimization problem. We present this algorithm and the implementation details in \cref{sec:implementation}.

\subsection{Preparation of Candidate Models}

Without any constraints, there are $2^p-1$ nonempty possible candidate models, which is huge in high-dimensional settings. Thus, we have to discuss how to prepare candidate models. In \cite{ando2014model,ando2017weight} and \cite{wang2023jackknife}, the variables that survive in a screening procedure are split into some disjoint groups, and each candidate model only uses one group of covariates. However, such a strategy implies that there is only one candidate model that can contain the most important variables, potentially introducing uncertainty and more bias compared to a nested structure, although the constraint that the weights sum to one is relaxed. In our work, after ranking the covariates by importance in descending order (c.f. Section \ref{sec:pred}), we separate them into two parts: the nested part and the non-nested part. For the nested part, given a pre-specified positive integer $K_{ne}$, we divide the top $p_0$ ranked variables into $K_{ne}$ equal-sized groups $\set{\mc{M}_k}_{k\in\brac{K_{ne}}}$, each containing $d_1=p_0/K_{ne}$ variables (assuming $d_1$ is an integer for simplicity), and then we construct nested models by group-wise adding, that is, $\mc{A}_k=\mc{M}_1\cup\dots\cup \mc{M}_k, k=1,\dots,K_{ne}$. For the non-nested part which contains the remaining $p-p_0$ variables, we divide it into $K-K_{ne}$ groups $\set{\mc{M}_k}_{k=K_{ne}}^{K}$, each with $d_2=(p-p_0)/(K-K_{ne})$ variables (assuming $d_2$ is an integer), and then construct candidate models separately using each group: $\mc{A}_k=\mc{M}_k, k=K_{ne}+1,\dots,K$. Further details are provided in   \cref{sec:implementation}.

Similar to \cite{he2020functional,he2023rank}, we incorporate all covariates in our construction of candidate models. By considering both nested and non-nested parts, the size of the largest model in the nested part can be smaller than $p$, which facilitates the theoretical developments. Also, if the true coefficients are non-sparse, our proposed model averaging approach can reduce the approximation error compared to the methods using only the nested part, and it is more robust against model misspecification and parameter non-sparsity.

\section{Theoretical Results}\label{sec: theoretical result}
Following \cite{white1982maximum} and \cite{lv2014model}, define the pseudo-true values of the candidate models
\begin{align*}
    \beta^*_{(k)}=\argmin_{\beta\in\R^p,\beta_{\mc{A}_k^c}=0}\E L(Y_1,\langle \beta,X_1\rangle),~ k=1,\dots,K,
\end{align*}
and denote $B^*=(\beta_{(1)}^*,\beta_{(2)}^*,\dots,\beta_{(K)}^*)$. 
We also define the optimal weight $$w^*\in\argmin_{w\in\mc{W}^K}\E L\paren{Y_1,\<B^*w,X_1\>}.$$ We say that the $k$-th candidate model is correct if $\beta_j=0$ for any $j>p$ and $\beta^*_{(k)}=\beta^{*}\overset{\triangle}{=}(\beta_1,\dots,\beta_p)^\top$; otherwise, the model is called a wrong model. In this section, we investigate the theoretical properties for two scenarios, one of which is that all models are wrong, and the other is that correct models are included. 

%%%%%%%%%%%%%%%%%%%%%%%%%%%%%%%%%%%%%%%%%%%%%%%%%%%%%%%%
\subsection{All Candidate Models Are Wrong}
We first consider the scenario that all candidate models are wrong. In this subsection, we establish non-asymptotic bounds and asymptotic optimality for our weight estimator. We require the following assumptions.
\begin{assumption}\label{lasso:lossfunc}
$\inf_{y\in\mc{Y}}\frac{\partial^2L(y,\mu)}{\partial \mu^2}$ is positive and continuous in $\mu$. $\frac{\partial L(y,\mu)}{\partial\mu}$ is Lipschitz continuous in $\mu$: $\abs{\frac{\partial L(y,\mu_1)}{\partial\mu}-\frac{\partial L(y,\mu_2)}{\partial\mu}}\leq L_u\abs{\mu_1-\mu_2}$ for any $y\in\mc{Y}$ and $\mu_1,\mu_2\in\R$, where $L_u$ is a positive constant. 
\end{assumption}

\begin{assumption}\label{lasso:covariate}
     $X_1$ is sub-gaussian with $\norm{X_1}_{\psi_2}= K_x$ and $\frac{\partial L(Y_1,0)}{\partial \mu}$ is sub-gaussian with $\norm{\frac{\partial L(Y_1,0)}{\partial \mu}}_{\psi_2}=K_y$, where $K_x$ and $K_y$ are positive constants. There exists a positive constant $\kappa_l$ such that $\lambdamin{\mathbf{\Sigma}}\ge\kappa_l> 0$, where $\mathbf{\Sigma}$ is the covariance matrix of $X_1$.
\end{assumption}

\begin{assumption}\label{HD-MA:beta-deviation} Define
\begin{align*}
    &~T_{n,q}^{(2)}(\delta)=\begin{cases}
        C\sqrt{s_u}(\delta+\sqrt{\frac{\log p_u}{n}}),&~~\text{if }q=0\\
        Ch_{q,u}^{\frac{q}{2}}(\delta+\sqrt{\frac{\log p_u}{n}})^{1-\frac{q}{2}},&~~\text{if }q\in(0,1]
    \end{cases},\\
    &~T_{n,q}^{(1)}(\delta)=\begin{cases}
        Cs_u(\delta+\sqrt{\frac{\log p_u}{n}}),&~~\text{if }q=0\\
        Ch_{q,u}^{q}(\delta+\sqrt{\frac{\log p_u}{n}})^{1-q},&~~\text{if }q\in(0,1]
    \end{cases},
\end{align*}where $\delta\in(0,1)$, $p_u=\max_{k\in\brac{K}}\abs{\mc{A}_k}, s_u=\max_{k\in\brac{K}}\lzero{\beta^*_{(k)}}$ and $h_{q,u}=\max_{k\in\brac{K}}\norm{\beta^*_{(k)}}_q$ with $q\in(0,1]$. There exists a constant $q\in[0,1]$ such that for each $k\in\brac{K}$, 
\begin{align*}
    \ltwo{\hat{\beta}_{(k)}-\beta^*_{(k)}}\leq T_{n,q}^{(2)}(\delta) ~\text{ and }~ 
    \lone{\hat{\beta}_{(k)}-\beta^*_{(k)}}\leq T_{n,q}^{(1)}(\delta)
\end{align*}
    hold simultaneously with probability at least $1-4\exp\paren{-n\delta^2}$. In addition, $\ltwo{\beta^*_{(k)}}$ are uniformly bounded by a positive constant $K_\beta$ for $k\in\brac{K}$.
\end{assumption}
\cref{lasso:lossfunc} is mild as it is satisfied by a large class of loss functions. \cref{lasso:covariate} is similar to Condition GLM1 in \cite{negahban2012unified} and Condition 3.1 in \cite{li2024estimation}. These conditions are commonly used in the field of high-dimensional data analysis.  \cref{HD-MA:beta-deviation} are some non-asymptotic upper bounds for estimators based on candidate models, which is similar to Condition 2.1 in \cite{battey2018distributed}, Condition C2 in \cite{he2023rank} and Condition C5 in \cite{he2020functional}. \cref{HD-MA:beta-deviation} can be proved if the pseudo-true values of candidate models are sufficiently sparse for some $q\in[0,1]$. We allow the pseudo-true values to have either hard sparsity $(q=0)$ or soft sparsity $(q\in(0,1])$. To verify the upper bounds given in \cref{HD-MA:beta-deviation}, we need to impose some conditions on the regularizer in terms of the univariate function $r_\lambda:\R\to\R$.

\begin{assumption}\label{con:regularizer}(i) $r_\lambda(0)=0$ and $r_\lambda(t)=r_\lambda(-t)$ for $t\in\R$; (ii) On the nonnegative real line, $r_\lambda$ is nondecreasing and $r_\lambda(t)/t$ is nonincreasing in $t$; (iii) $r_\lambda$ is differentiable for all $t\neq0$ and subdifferentiable at $t=0$ with $r_\lambda(0+)=\lambda$. In addition, there exists a positive constant $\kappa_r>0$ such that $r_\lambda(t)+\frac{\kappa_r}{2}t^2$ is convex, where $\kappa_r<c\kappa_l$ for some $c>0$ only depending on $K_\beta$ in \cref{HD-MA:beta-deviation}
and the loss function. Here, $\kappa_l$ is defined in \cref{lasso:covariate}.
\end{assumption}

\cref{con:regularizer} is similar to Assumption 1 in \cite{loh2015regularized}. Many common regularizers satisfy \cref{con:regularizer} such as Lasso and the folded concave penalties MCP and SCAD. The following proposition shows that the upper bounds in \cref{HD-MA:beta-deviation} holds under \cref{lasso:lossfunc,lasso:covariate,con:regularizer} and some conditions on the pseudo-true values. 

\begin{proposition}\label{lasso}
    Under \cref{lasso:lossfunc,lasso:covariate,con:regularizer}, suppose that $\lambda_{n,k}=c\paren{\delta+\sqrt{\frac{\log p_k}{n}}}$ for each $k\in\brac{K}$ and $\ltwo{\beta^*_{(k)}}$ are uniformly bounded by some positive constant $K_\beta$, where $p_k=\abs{\mc{A}_k}$. In addition, we assume that $\beta_{(k)}^*$ and $\hat{\beta}_{(k)}$ are feasible in the $\ell_1$-ball $\set{\beta\in\R^p: \lone{\beta}\leq R/2}$ with $\lambda_{n,k}R\leq c$ uniformly for each $k\in\brac{K}$ (For the Lasso penalty, this condition is unnecessary).  For each $k\in\brac{K}$, we have the following results.
    \begin{enumerate}
        \item If $\sqrt{s_k}\lambda_{n,k}\leq c$, then for large enough $n$,
    \begin{align*}
        \ltwo{\hat{\beta}_{(k)}-\beta^*_{(k)}}\leq C\sqrt{s_k}\paren{\delta+\sqrt{\frac{\log p_k}{n}}} ~\text{ and }~
        \lone{\hat{\beta}_{(k)}-\beta^*_{(k)}}\leq Cs_k\paren{\delta+\sqrt{\frac{\log p_k}{n}}}
    \end{align*} hold simultaneously with probability at least $1-2\exp\paren{-n\delta^2}-2\exp\paren{-cn}$, where $s_k=\lzero{\beta^*_{(k)}}$.
    \item If $h_{q,k}^{\frac{q}{2}}\lambda_{n,k}^{1-\frac{q}{2}}\leq c$, then for large enough $n$,
    \begin{align*}
        \ltwo{\hat{\beta}_{(k)}-\beta^*_{(k)}}\leq Ch_q^{\frac{q}{2}}\paren{\delta+\sqrt{\frac{\log p_k}{n}}}^{1-\frac{q}{2}} ~\text{ and }~
        \lone{\hat{\beta}_{(k)}-\beta^*_{(k)}}\leq Ch_{q,k}^{q}\paren{\delta+\sqrt{\frac{\log p_k}{n}}}^{1-q}
    \end{align*} hold simultaneously with probability at least $1-2\exp\paren{-n\delta^2}-2\exp\paren{-cn}$, where $h_{q,k}=\norm{\beta^*_{(k)}}_q$ with $q\in(0,1]$.
    \end{enumerate}
\end{proposition}

\medskip
In practice, for some low-dimensional candidate models, the unpenalized estimator may be preferred. Denote the candidate models with model dimension not greater than $d$ by $\mc{A}^d=\set{k\in\brac{K}|~\abs{\mc{A}_k}\leq d}$. When $d$ is fixed, the upper bounds in \cref{HD-MA:beta-deviation} still holds if we take $\lambda_{n,k}=0$ for $k\in\mc{A}^d$.
\begin{proposition}\label{unpenalized}
    Under \cref{lasso:lossfunc,lasso:covariate}, suppose that $\lambda_{n,k}=0$ for each $k\in\mc{A}^d$ and $\ltwo{\beta^*_{(k)}}$ are uniformly bounded. For any $k\in\mc{A}^d$, if $\sqrt{d}(\delta+\sqrt{\frac{\log d}{n}})\leq c$, then for large enough $n$,
    \begin{align*}
        \ltwo{\hat{\beta}_{(k)}-\beta^*_{(k)}}\leq C\sqrt{d}\paren{\delta+\sqrt{\frac{\log d}{n}}} ~\text{ and }~
        \lone{\hat{\beta}_{(k)}-\beta^*_{(k)}}\leq Cd\paren{\delta+\sqrt{\frac{\log d}{n}}}
    \end{align*}
    hold simultaneously with probability at least $1-2\exp\paren{-n\delta^2}-2\exp\paren{-cn}$.
\end{proposition}

\subsubsection{Upper Bound}
We first introduce some notations. Denote 
\begin{align*}
    &~\Psi_{n,q}=\begin{cases}
        \sqrt{s_u}\sqrt{\frac{\log (pK)}{n}},&~~\text{if }q=0\\
        h_{q,u}^{\frac{2}{2-q}}\sqrt{\frac{\log (pK)}{n}}\vee h_{q,u}^{\frac{q}{2}}(\frac{\log (pK)}{n})^{\frac{1}{2}-\frac{q}{4}},&~~\text{if }q\in(0,1)
    \end{cases},\\
    &~T_{n,q}=\begin{cases}
        \sqrt{s_u}\sqrt{\frac{\log (pK)}{n}},&~~\text{if }q=0\\
        h_{q,u}^{\frac{q}{2}}(\frac{\log (pK)}{n})^{\frac{1}{2}-\frac{q}{4}},&~~\text{if }q\in(0,1)
    \end{cases}.
\end{align*}
Given a set $S_*\subseteq\brac{K}$, we take the set of restrictions
\begin{align*}
    \mathbb{C}(S_*,w^*)=\setBig{\Delta\in\mc{W}^K-w^*\Big| \lone{\Delta_{S_*^c}}-2\lone{w^*_{S_*^c}}\leq \lone{\Delta_{S_*}}}.
\end{align*}
In our proof of \cref{MA}, we can show that with high probability, $\hat{w}-w^*$ lies in the set $\mathbb{C}(S_*,w^*)$ for some set $S_*\subseteq\brac{K}$. Define the compatibility factor (restricted $\ell_1$-eigenvalue) for a set $S_*\subseteq\brac{K}$ as
\begin{align*}
    \phi^2(S_*, B^*, w^*)=\min\setBig{\frac{\abs{S_*}~\ltwo{B^*\Delta}^2}{\lone{\Delta_{S_*}}^2}\Big|~\Delta\in\mathbb{C}(S_*,w^*)}.
\end{align*}
This quantity is similar to the compatibility factor introduced in \cite{van2009conditions,zhang2012general} and \cite{javanmard2014confidence}. In the following, we drop $B^*$ and $w^*$ and rewrite $\phi^2(S_*, B^*, w^*)$ as $\phi^2(S_*)$ for brevity. Note that $w^*$ is not necessarily unique. In fact, $\hat{w}$ could be close to any point in the set $\argmin_{w\in\mc{W}^K}\E L(Y_1,\<B^*w,X_1\>)$ if the diameter of this set goes to zero as $n\to \infty$. See \cref{remark:uniqueness} for more discussion.  The following theorem presents the upper bound of $\lone{\hat{w}-w^*}$ for any $w^*\in\argmin_{w\in\mc{W}^K}\E L(Y_1,\<B^*w,X_1\>)$.

\begin{theorem}\label{MA} Under \cref{lasso:covariate,lasso:lossfunc,HD-MA:beta-deviation} with some constant $q\in[0,1]$, if $T_{n,q}\leq c_1$, then for any $w^*\in\argmin_{w\in\mc{W}^K}\E L(Y_1,\<B^*w,X_1\>)$, with probability at least $1-C_2(pK)^{-1}-C_2e^{-cn}$,
\begin{align}
\lone{\hat{w}-w^*}\leq&~C_1\frac{\abs{S_*}}{\phi^2(S_*)}T_{n,q}+C_1\lone{w^*_{S_*^c}}+C_1\sqrt{\frac{\abs{S_*}}{\phi^2(S_*)}\Psi_{n,q}\lone{w^*_{S_*^c}}}\label{eq:oracle-inequality}
\end{align}
uniformly holds for any $S_*\subseteq\brac{K}$ with cardinality $\abs{S_*}\leq c_2\phi^2(S_*)\Psi_{n,q}^{-1}$.
\end{theorem}

\begin{remark}{\it
    Since $\lone{\hat{w}}=1$, $\hat{w}$ can be viewed as a Lasso estimator under a simplex constraint. The standard theory for $M$-estimators \citep{negahban2012unified,loh2015regularized} requires that the gradient of the population risk at its unique minimizer is zero. However, the gradient of the risk $\E L(Y_1,\<B^*w,X_1\>)$ at the point $w=w^*$ may not be zero since the simplex $\mc{W}^K$ does not contain any neighborhood of $w^*$ in $\R^K$. Thus, the classic proof technique would give the upper bound \eqref{eq:oracle-inequality} a quantity $\linf{\E\parenBig{\frac{\partial L(Y_1,\<B^*w^*,X_1\>)}{\partial \mu}B^{*\top}X_1}}$, which is a deterministic positive constant if we fix the candidate model and the dimension of covariates. A similar quantity appears in Theorem 1 of \cite{he2023rank}. In contrast, in our theory, we remove the gradient term by noting the optimality condition:
 \begin{align}\label{eq:nonnegative grad}
     \Big\<\E\parenBig{\frac{\partial L\paren{Y_1,\<B^*w^*,X_1\>}}{\partial\mu}B^{*\top}X_1},w-w^*\Big\>\ge 0
 \end{align}
 for any feasible $w\in\mc{W}^K$. By the inequality \eqref{eq:nonnegative grad}, the upper bound \eqref{eq:oracle-inequality} can be derived by bounding the following quantity:
 \begin{align*}
     \linfBig{\dn\sum_{m=1}^J\sum_{i\in \mc{I}_m}\frac{\partial L\paren{Y_i,\<\hat{B}^{\brac{-m}}w^*,X_i\>}}{\partial \mu}\hat{B}^{\brac{-m}\top}X_i-\E\parenBig{\frac{\partial L\paren{Y_1,\<B^*w^*,X_1\>}}{\partial\mu}B^{*\top}X_1}}.
 \end{align*}
 The quantity $T_{n,q}$ in \eqref{eq:oracle-inequality} is the uniform $\ell_2$ estimation error of $\hat{\beta}_{(k)}$'s by \cref{HD-MA:beta-deviation}. When $q=0$ or $q\in(0,1]$ with $\frac{\log p}{n}\leq h_{q,u}^{2-8/(2q-q^2)}$, we have $T_{n,q}=\Psi_{n,q}$ and the upper bound \eqref{eq:oracle-inequality} can be simplified to $C_1\abs{S_*}T_{n,q}/\phi^2(S_*)+C_1\lone{w^*_{S^c_*}}$.}
\end{remark}

  It is worth noting that the inequality \eqref{eq:oracle-inequality} provides a family of upper bounds without any unnatural conditions and does not depend on the sparsity of $w^*$ compared to Theorem 1 in \cite{he2023rank}. By choosing $S_*$ to balance the terms $\abs{S_*}T_{n,q}/\phi^2(S_*)$ and $\lone{w^*_{S^c_*}}$, we can show that $\hat{w}$ is consistent. For example, when we fix the dimension of covariates and the sets of candidate models, we can immediately obtain that $\hat{w}$ converges to $w^*$ at a $n^{-1/2}$ rate and $\hat{\beta}_{MA}(\hat{w})$ is consistent to $\beta^*(w^*)$ by taking $S_*=\text{supp}(w^*)$, where $\beta^*(w^*)=B^*w^*$ is the minimizer of the prediction risk over the convex hull of $\beta^*_{(k)}$'s. \cite{fang2023asymptotic} shows that the least squares model averaging with nested models asymptotically assigns weight one to the largest candidate model when all models are wrong. This inspires us to choose $S_*=\set{K_{ne}}$. By the definition of $\phi^2(S_*)$, we have $\phi^2(S_*)\ge\sum_{j\in\mc{A}_{K_{ne}}\backslash\mc{A}_{K_{ne}-1}}\beta_j^2$. Assume that $\sum_{j\in\mc{A}_{K_{ne}}\backslash\mc{A}_{K_{ne}-1}}\beta_j^2\ge n^{-\tau}$ for some constant $\tau\in[0,1/2)$, which is similar to Condition (C3) in \cite{fang2023asymptotic}. If $\lone{w^*_{S_*^c}}n^{-\tau}\lesssim T_{n,q}$ and $\Psi_{n,q}=T_{n,q}$, then we have
    \begin{align*}
        1-\hat{w}_{K_{ne}}\leq\lone{\hat{w}-w^*}=O_p(T_{n,q}n^{\tau})
    \end{align*}
    and
    \begin{align*}
        \ltwo{\hat{\beta}_{MA}(\hat{w})-\beta^*_{(K_{ne})}}\leq K_\beta\lone{\hat{w}-w^*}+\max_{k\in\brac{K}}\ltwo{\hat{\beta}_{(k)}-\beta^*_{(k)}}+2K_\beta\lone{w^*_{S_*^c}}=O_p(T_{n,q}n^{\tau}).
    \end{align*}
    In this case, if $T_{n,q}n^{\tau}\to0$, then the weight estimator $\hat{w}$ is consistent. We discuss the support of $w^*$ in \cref{remark:support of w^*} and give some specific examples for the convergence rate of $\hat{w}$ in \cref{remark:rate of w} to better understand \cref{MA}.

\subsubsection{Lower Bound}
\cref{MA} provides the upper bound for the $\ell_1$ error of our proposed weight estimator. A natural question is whether this upper bound is optimal. Here, we provide the minimax lower bound for weight estimation under the quadratic loss function. Let $S^*=\text{supp}(w^*)$. Define the parameter space
\begin{align*}
    \Theta=\setBig{(w^*, B^*)\Big|\lzero{w^*}\leq s_w, w^*\in\mc{W}^K, \max_{k\in\brac{K}}\ltwo{\beta^*_{(k)}}\leq C, \phi^2(S^*)\ge \phi^2_n}.
\end{align*}
 We focus on the following subset of the distribution space that satisfies \cref{lasso:covariate,lasso:lossfunc}:
\begin{align*}
    \mc{P}_{Y,X}^{1:n}=\setBig{\set{(Y_i,X_i)}_{i=1}^n\overset{i.i.d.}{\sim} Q\Big|Y_1|X_1\sim N(\sum_{j=1}^p\beta_jX_{1j},1), X_1\sim N(0,I_p), (w^*,B^*)\in\Theta}.
\end{align*}
The following theorem gives a minimax lower bound among all potential weight estimators regarding the class $\mc{P}_{Y,X}^{1:n}$.

\begin{theorem}\label{HD-MA:minimax}
There exists $K$ pairwise disjoint index sets $\mc{A}_k\subseteq\brac{p}, k=1,\dots,K$ such that if $s_w\ge 2$, $s_w\leq\frac{2}{5}K$ and $\phi_n^2\leq c(s_w^2\log (K/s_w)/n)^{1/3}$, then
    \begin{align}\label{eq:minimax lower bound}
\inf_{\hat{w}}\sup_{Q\in\mc{P}_{Y,X}^{1:n}}\P_{(Y_i,X_i)\sim Q}\parenBig{\lone{\hat{w}-w^*}\ge Cs_w\sqrt{\frac{\log(K/s_w)}{\phi_n^2n}}}\ge\frac{1}{2}.
    \end{align}
\end{theorem}
To prove this lower bound, by Fano's lemma, we need to find a separated set on the parameter space $\Theta$, which is a challenging task as the oracle weights $w^*$ and the pseudo-true values $\beta^*_{(k)}$'s are related to the underlying parameter $\bm{\beta}$. Due to the simplex constraint, the condition $\phi_n^2\leq c(s_w^2\log (K/s_w)/n)^{1/3}$ allows us to transform this task into the problem of finding a separated set on the simplex $\mc{W}^K$. For any $K$ pairwise disjoint index sets $\mc{A}_k\subseteq\brac{p}, k=1,\dots,K$, under the distribution class $\mc{P}_{Y,X}^{1:n}$, if the maximum model size $s_u$ is fixed and $K\asymp p$, then by \cref{MA}, the following upper bound
\begin{align}\label{eq: a upper bound}
    \lone{\hat{w}-w^*}\leq C\frac{s_w}{\phi_n^2}\sqrt{\frac{\log K}{n}}
\end{align}
holds with high probability. When $s_w\lesssim K^\alpha$ for some $\alpha\in[0,1)$, the lower bound in \eqref{eq:minimax lower bound} has the order $O(s_w\sqrt{\frac{\log K}{n\phi_n^2}})$, which only differs from the upper bound by a multiple of the square root of the compatibility factor.

\subsubsection{Asymptotic Optimality}
Next, we will establish the result of asymptotic optimality. Let $(Y_{n+1},\mathbf{x}_{n+1})$ be a new sample drawn from the same distribution as $(Y_1,\mathbf{x}_1)$ and $X_{n+1}=(\mathrm{x}_{n+1,1},\dots,\mathrm{x}_{n+1,p})^\top$. The conditional risk based on the estimators $\set{\hat{\beta}_{(k)}}_{k=1}^K$ and the pseudo-true values $\set{\beta^*_{(k)}}_{k=1}^K$ are
\begin{align*}
    R\paren{w}=\E \bracBig{L\paren{Y_{n+1},\<\hat{B}w,X_{n+1}\>}\Big|\set{(Y_i,X_i)}_{i=1}^n}
\end{align*}
and 
\begin{align*}
    R^*\paren{w}=\E L\paren{Y_{n+1},\<B^*w,X_{n+1}\>},
\end{align*}
respectively. Define $\xi_n\overset{\triangle}{=}\inf_{w\in\mc{W}}R^*(w)-R_{Bayes}$, where $R_{Bayes}=\E L(Y_{n+1},\sum_{j=1}^\infty\beta_j\mathrm{x}_{n+1,j})$ is the Bayes risk.
To show asymptotic optimality in the sense that the weight estimator $\hat{w}$ minimizes excess prediction risk $R(w)-R_{Bayes}$, we require the following assumptions.

\begin{assumption}\label{HD-MA:lossfunc}
    One of the following conditions holds:
    \begin{enumerate}
        \item $L\paren{y,\mu}=\frac{1}{2}\paren{y-\mu}^2$.
        \item There exists a positive constant $L_u^\prime$ such that $\abs{L\paren{y,\mu_1}-L\paren{y,\mu_2}}\leq L_u^\prime\abs{\mu_1-\mu_2}$ for any $y\in\mc{Y}$, $\mu_1,\mu_2\in\R$.
    \end{enumerate}
\end{assumption}

\begin{assumption}\label{HD-MA:xi}
       $T_{n,q}\to 0, \xi_n^{-1}T_{n,q}\to 0$ as $n\to\infty$.
\end{assumption}

 \cref{HD-MA:lossfunc} is an alternative to \cref{lasso:lossfunc}. It requires that the loss function $L$ is either the $\ell_2$ loss or Lipschitz continuous in the second argument. There are many loss functions satisfying \cref{HD-MA:lossfunc}, such as quantile regression loss, cross entropy loss, Huber loss and hinge loss. This assumption is used to show that 
 \begin{align*}
     \sup_{w\in\mc{W}^K}\absBig{\dn\sumn L(Y_i,\<B^*w,X_i\>)-\E L(Y_i,\<B^*w,X_i\>)}
 \end{align*}tends to $0$ in expectation. The expectation of such a quantity can be bounded by the Rademacher complexity using a symmetrization argument. If this quantity has the convergence rate of order $O_p(T_{n,q})$, then the result of asymptotic optimality  in \cref{HD-MA:asymptotical-optimality} still holds under \cref{lasso:lossfunc} but not \cref{HD-MA:lossfunc}. \cref{HD-MA:xi} requires that the divergence rate of $n\xi_n$ is fast enough. This corresponds to the case where the approximation error dominates the uniform estimation error of $\hat{\beta}_{(k)}$'s. It is similar to Assumption 5 of \cite{zhang2023model}, Condition 7 of \cite{ando2014model}, Condition C.6 of \cite{zhang2016optimal}, and Condition A3 of \cite{ando2017weight}. In the case where $\xi_n$ is bounded away from zero, the sufficient condition for this assumption is $T_{n,q}\to0$. A necessary condition of \cref{HD-MA:xi} is that all candidate models are misspecified since $\xi_n=0$ if there is a correct model. To save space here, we give an example for the convergence rate of $\xi_n$ in \cref{remark:xi}. The following theorem provides the asymptotic optimality of the weights estimator with respect to the excess prediction risk.

\begin{theorem}\label{HD-MA:asymptotical-optimality}
    Under \cref{HD-MA:lossfunc,lasso:covariate,HD-MA:beta-deviation,HD-MA:xi}, we have
    \begin{align}
        \frac{R\paren{\hat{w}}-R_{Bayes}}{\inf_{w\in\mc{W}^K}R\paren{w}-R_{Bayes}}=1+O_p(\xi_n^{-1}T_{n,q}).\label{eq:AOPrate}
    \end{align}
     Furthermore, if
    \begin{align*}
        \xi_n^{-1}\sup_{w\in\mc{W}^K}\absBig{R(w)-R^*(w)}
    \end{align*}
    is uniformly integrable, then we have
    \begin{align*}
        \frac{\Bar{R}\paren{\hat{w}}-R_{Bayes}}{\inf_{w\in\mc{W}^K}\Bar{R}\paren{w}-R_{Bayes}}= 1+o_p(1),
    \end{align*}
    where $\Bar{R}(w)=\E R(w)$.
\end{theorem}
The condition of uniform integrability in \cref{HD-MA:asymptotical-optimality} is similar to Assumption 3 in \cite{zhang2023model} and is only imposed to ensure the convergence of expectation. According to \cref{HD-MA:asymptotical-optimality}, the model averaging estimator with weights $\hat{w}$ is asymptotically equivalent to the estimator with the best weights. To have a better performance, we require the size of the nested part $p_0$ to diverge such that the denominator on the left-hand side of \eqref{eq:AOPrate} is small enough. To the best of our knowledge, \cref{HD-MA:asymptotical-optimality} is the first result of asymptotic optimality that allows the dimension $p$ and the number of candidate models $K$ to exceed the sample size $n$. 

To understand \cref{HD-MA:asymptotical-optimality} more clearly, assume that the conditional density function of $Y_1$ given $\mathbf{x}_1$ is $f_{Y|\mathbf{x}}(y|\sum_{j=1}^\infty \beta_j\mathrm{x}_j)$ and $L(y,\mu)=-\log f_{Y|\mathrm{x}}(y|\mu)$. Then the excess prediction risk is
\begin{align*}
    \bar{R}(w)-R_{Bayes}=\E\setBig{\E\bracBig{-\log \frac{f_{Y|\mathbf{x}}(Y_{n+1}|\<\hat{B}w,X_{n+1}\>)}{f_{Y|\mathbf{x}}(Y_{n+1}|\sum_{j=1}^\infty\beta_j\mathrm{x}_{n+1,j})}\Big|\hat{B},\mathbf{x}_{n+1}}},
\end{align*}
and the proposed estimator is asymptotically optimal in the sense of minimizing the KL divergence. Consider a linear model $Y_1=\sum_{j=1}^\infty\beta_j\mathrm{x}_{1j}+\eps_1$. By a simple calculation, we have $2\set{\bar{R}(w)-R_{Bayes}}=\E\brac{\<\hat{B}w,X_{n+1}\>-\E(Y_{n+1}|\mathbf{x}_{n+1})}^2$, which is the same as the risk function in \cite{zhang2023model}. Thus, we obtain the same result as Theorem 1 of \cite{zhang2023model} in the high-dimensional setting. 

%%%%%%%%%%%%%%%%%%%%%%%%%%%%%%%%%%%%%%%%%%%%%%%%%%%%%%%
\subsection{There Exists at Least One Correct Model}
 Now we consider the scenario where at least one of the candidate models is correctly specified. In this subsection, we show that the proposed estimator assigns all weights on the correct models with probability tending to one. We propose a post-averaging debiased estimator. Following its Bahadur representation, the Gaussian and bootstrap approximation results are obtained, which allows us to construct simultaneous confidence intervals using multiplier bootstrap.
 
 \subsubsection{Weight Consistency} 
 Note that the optimal weights $w^*$ are not unique when there exist two or more correct models. In this subsection, we assume that the $k$-th candidate model is correct for $k= K^*,\dots,K_{ne}$ while the remaining candidate models are wrong, that is, the correct models are nested within the nested part of the candidate model set. This means that the first $K^*-1$ candidate models in the nested part and all candidate models in the non-nested part are misspecified and $\beta^*_{(k)}=\beta^*$ for $k\in K^*,\dots,K_{ne}$. Assume that $\tilde{w}^*\overset{\triangle}{=}(w_1^*,\dots,w^*_{K^*-1},\sum_{k=K^*}^{K_{ne}}w^*_k,w^*_{K_{ne}+1},\dots,w^*_K)$ is the unique minimizer of $\E L\paren{Y_1,\<\tilde{B}^*\tilde{w},X_1\>}$ over $\tilde{w}\in\mc{W}^{K_0}$, where $K_0=K-(K_{ne}-K^*)$ and $\tilde{B}^*=(\beta_{(1)}^*,\dots,\beta_{(K^*)}^*,\beta^*_{(K_{ne}+1)},\dots,\beta^*_{(K)})$. Consequently, we have $\tilde{w}^*_{K^*}=1$, which indicates that the sum of the optimal weights on the correct models is equal to one. Let $\breve{\beta}^{\brac{-m}}_{(K^*)}=\frac{\sum_{k=K^*}^{K_{ne}}\hat{w}_k\hat{\beta}^{\brac{-m}}_{(k)}}{\sum_{k=K^*}^{K_{ne}}\hat{w}_k}$. Then,
\begin{align*}
    \breve{w}\overset{\triangle}{=}&~(\hat{w}_1,\dots,\hat{w}_{K^*-1},\sum_{k=K^*}^{K_{ne}}\hat{w}_k,\hat{w}_{K_{ne}+1},\dots,\hat{w}_K)\\
    =&~\argmin_{w\in\mc{W}^{K_0}}\sum_{m=1}^J\sum_{i\in\mc{I}_m}L\paren{Y_i,w_{K^*}\breve{\beta}_{(K^*)}^{\brac{-m}\top}X_i+\sum_{k\neq K^*,\dots,K}w_k\hat{\beta}_{(k)}^{\brac{-m}\top}X_i},
\end{align*}
and the consistency of $\breve{w}$ can be seen as a special case of \cref{MA} to some extent.
However, since $\breve{\beta}^{\brac{-m}}_{(K^*)}$ is related to $\set{(Y_i,X_i)}_{i\in\mc{I}_m}$ for each $m\in\brac{J}$, we cannot apply \cref{MA} directly to derive the consistency of $\breve{w}$, i.e., $\lone{\breve{w}-\tilde{w}^*}=2\abs{1-\sum_{k=K^*}^{K_{ne}}\hat{w}_k}=o_p(1)$. To deal with this issue, we need to generalize the restricted strong convexity result in \cite{negahban2012unified}. We also need the following assumption. 

\begin{assumption}\label{HD-MA:NOunique} One of the following conditions holds.
\begin{enumerate}
    \item $\inf_{\mu\in\R}\inf_{y\in\mc{Y}}\frac{\partial^2L(y,\mu)}{\partial\mu^2}\ge\eta_l>0$.
    \item \cref{HD-MA:beta-deviation} holds with $q\in[0,1)$. In addition, $s_u\leq c (\log (pK))^{-1}\sqrt{n}$ if $q=0$; $h_{q,u}\leq c (\log (pK))^{\frac{q-3}{2q}}n^{\frac{1-q}{2q}}$ if $q\in(0,1)$.
\end{enumerate}
\end{assumption}
\cref{HD-MA:NOunique} requires that either the loss function is strongly convex or there are some sparsity conditions on candidate models. In the case where $\abs{\text{supp}(\beta^*)}=s$ and the size of models in the non-nested part $d_2$ is fixed, the second condition can be easily satisfied if $s\lesssim\sqrt{n}/\log p$.

Similar to the compatibility factor defined in the previous subsection, denote \begin{align*}
    \tilde{\phi}^2=\min\setBig{\frac{\ltwo{\tilde{B}^*\Delta}^2}{\abs{\Delta_{K^*}}^2}\Big|~\Delta\in\mc{W}^{K_0}-\tilde{w}^*}.
\end{align*}
\begin{theorem}\label{weights2one}
    Under \cref{lasso:covariate,lasso:lossfunc,HD-MA:beta-deviation,HD-MA:NOunique}, if $\Psi_{n,q}\leq c_1\tilde{\phi}^2$ with some fixed $q\in[0,1)$, then with probability greater than $1-C_2(pK)^{-1}-C_2e^{-cn}$, we have
\begin{align*}
1-\sum_{k=K^*}^{K_{ne}}\hat{w}_k\leq C_1T_{n,q}/\tilde{\phi}^2.
\end{align*}
\end{theorem}

Note that $\tilde{\phi}$ can be bounded by the minimum singular value of the matrix $\tilde{B}^*_{ne}=(\beta^*_{(1)},\beta^*_{(2)},\dots,\beta^*_{(K^*)})$ from below. In addition, if $\lambda_{\min}(\tilde{B}^{*\top}_{ne}\tilde{B}^*_{ne})\ge c$, $d_2$ is fixed and the sparsity level of the underlying coefficients $\abs{\text{supp}(\beta^*)}=s$, then we have $1-\sum_{k=K^*}^{K_{ne}}\hat{w}_k=O_p(\sqrt{\frac{s\log p}{n}})$. The following result indicates that the model averaging estimator is consistent.
\begin{corollary}\label{weights2one:cor}
      Suppose that the minimum singular value of $\tilde{B}^*_{ne}$ is bounded by some constant from below, $\abs{\text{supp}(\beta^*)}=s$ and $d_2$ is fixed.  Under the assumptions in \cref{weights2one}, if $\max_{k\in\brac{K}}\lone{\beta^*_{(k)}}\leq C$, then we have 
      \begin{align}\label{eq: estimation error}
          \ltwo{\hat{\beta}_{MA}(\hat{w})-\beta^*}=O_p(\sqrt{\frac{s\log p}{n}})~~\text{ and }~~~ \lone{\hat{\beta}_{MA}(\hat{w})-\beta^*}=O_p(s\sqrt{\frac{\log p}{n}}).
      \end{align}
\end{corollary}
This corollary can be obtained immediately by \cref{weights2one} and the triangle inequality. Thus we omit the proof. When there is at least one correct candidate model, the model averaging estimator $\hat{\beta}_{MA}(\hat{w})=\hat{B}\hat{w}$ has the same convergence rate as the Lasso estimator, which is nearly minimax rate-optimal. Since \cref{HD-MA:xi} fails in this case, the result of asymptotic optimality in \cref{HD-MA:asymptotical-optimality} does not hold; however, we have the following risk consistency result as a parallel to \cref{HD-MA:asymptotical-optimality}.
\begin{theorem}\label{riskconsistency}
    Under the assumptions in \cref{weights2one:cor}, we have
    \begin{align*}
        R(\hat{w})-R_{Bayes}=O_p\parenBig{\frac{s\log p}{n}}.
    \end{align*}

\end{theorem}

\subsubsection{Gaussian and Bootstrap Approximation}
In \cref{weights2one:cor}, we have established the consistency of the proposed model averaging estimator. However, it cannot be used directly for statistical inference. In fact, the penalized estimators are generally biased in the absence of the beta-min condition, which may not be easy to verify. Inspired by the one-step debiased estimator introduced in \cite{zhang2014confidence} and \cite{ning2017general}, we consider the following one-step debiased estimator:
\begin{align}\label{debiased lasso}
    \tilde{\beta}(\hat{w})=\hat{\beta}_{MA}(\hat{w})-\hat{W}\dn\sumn\frac{\partial L(Y_i,
    \<\hat{\beta}_{MA}(\hat{w}),X_i\>)}{\partial \mu}X_i,
\end{align}
where $\hat{W}$ is an estimator of the inverse of Hessian $J=\E\brac{\frac{\partial^2 L(Y_1,\<\beta^*,X_1\>)}{\partial \mu^2}X_1X_1^\top}$. The second term in \eqref{debiased lasso} is the shrinkage bias, which need to be removed according to the KKT conditions \citep[see, for example, ][]{van2014asymptotically}.

Motivated by \cite{yan2023confidence} and \cite{cai2025statistical}, for technical convenience, we use a modified version of CLIME estimator \citep{cai2011constrained} to estimate $J^{-1}$:
\begin{align}\label{CLIME estimator}
    \hat{W}=\argmin_{W\in\R^{p\times p}}\matrixnorm{W}_1\text{   ~~~ s.t. } \matrixnorm{W\dn\sumn\frac{\partial^2 L(Y_i,\hat{\beta}_{MA}(\hat{w}),X_i\>)}{\partial \mu^2}X_iX_i^\top-I_p}_{\max}\leq \gamma_n,
\end{align}
where $I_p$ is the $p$-dimensional identical matrix and $\gamma_n$ is a tuning parameter. In \cref{lemma: CLIME} in \cref{appendix:lemma}, we show that $J^{-1}$ is feasible for this optimization problem with probability tending to one. Note that the estimator $\hat{W}$ obtained in \eqref{CLIME estimator} might not be symmetric. To address this issue, similar to \cite{cai2011constrained}, in practice, we can define a symmetrized estimator $\hat{W}^*=(\hat{w}^*_{i,j})_{1\leq i,j\leq p}$ with $\hat{w}^*_{i,j}=\hat{w}_{i,j}\indic{\abs{\hat{w}_{i,j}}\leq\abs{\hat{w}_{j,i}}}+\hat{w}_{j,i}\indic{\abs{\hat{w}_{i,j}}>\abs{\hat{w}_{j,i}}}.$ Now, we impose some regularity conditions on the inverse of Hessian $J^{-1}$.
\begin{assumption}\label{inference: hessian}
    $J$ is a positive definite matrix. There exists some positive constant $M$ such that $\matrixnorm{J^{-1}}_1\leq M$. Moreover, $J^{-1}\overset{\triangle}{=}(\tilde{\mathbf{b}}_1,\dots,\tilde{\mathbf{b}}_p)^\top=(\tilde{b}_{i,j})_{1\leq i,j\leq p}$ is row-wisely sparse, that is, $\max_{1\leq i\leq p}\sum_{j=1}^p\abs{\tilde{b}_{i,j}}^{q^\prime}\leq s_0$ for some $q^\prime\in[0,1)$. The tuning parameter $\gamma_n$ satisfies $\gamma_n\sqrt{n}/(s\log^2(p))=o(1)$.
\end{assumption}
This assumption requires the inverse Hessian matrix $J^{-1}$ to be 
row-wisely $\ell_{q^\prime}$-sparse \citep{cai2011constrained}, which is similar to Assumption 5 in \cite{yan2023confidence} and Condition (B1) in \cite{cai2025statistical}. We allow $s_0$ to diverge with the sample size $n$. We also need to impose a smoothness condition to control moments of higher-order derivatives of the loss function.
\begin{assumption}\label{inference:smoothness}
There exists a positive constant $C>0$ such that $\abs{\frac{\partial^2 L(y,\mu_1)}{\partial \mu^2}-\frac{\partial^2 L(y,\mu_2)}{\partial \mu^2}}\leq C\abs{\mu_1-\mu_2}$ for any $y\in\mc{Y}$ and $\mu_1,\mu_2\in\R$.
\end{assumption}
According to the context of optimization, the Hessian Lipschitz condition is common in the convergence analysis in Newton-type methods. \cref{inference:smoothness} is always satisfied for linear and logistic models where the third-order derivative is bounded. We now establish the uniform Bahadur representation for the debiased estimator $\tilde{\beta}(\hat{w})$, which is crucial to obtain the Gaussian and bootstrap approximation results.

\begin{theorem}\label{inference: bahadur representation}
     Let $\eta_n=s_0\setBig{s\log^{(5-4q^\prime)/(2-2q^\prime)}(p)/\sqrt{n}}^{1-q^\prime}+s^2\log^{5/2}(p)/\sqrt{n}.$ Suppose that $\hat{\beta}(\hat{w})$ satisfies the estimation error bound in \eqref{eq: estimation error}. Under \cref{lasso:lossfunc,lasso:covariate,inference: hessian,inference:smoothness}, if $\eta_n=o(1)$, then we have
    \begin{align*}
        \linfBig{\sqrt{n}(\tilde{\beta}(\hat{w})-\beta^*)+\sqrt{n}J^{-1}\dn\sumn\frac{\partial L(Y_i,\<\beta^*,X_i\>)}{\partial \mu}X_i}=O_p\parenBig{\eta_n}.
    \end{align*}
\end{theorem}

Next, we will derive the Gaussian approximation result for the debiased estimator $\tilde{\beta}(\hat{w})$. Denote $S=\text{diag}(J)$ and $\Lambda=\E\brac{(\frac{\partial L(Y_1,\<\beta^*,X_1\>)}{\partial \mu})^2X_1X_1^\top}$. 
\begin{assumption}\label{inference:Lambda}
    $\lambda_{\min}(\Lambda)\ge c$ and $\lambda_{\min}(J)/\lambda_{\max}(J)\ge c$ for some positive constant $c>0.$
\end{assumption}
This assumption requires that the condition number of the Hessian matrix $J$ is bounded from above and the minimum eigenvalue of the covariance of the gradient is bounded away from zero. For example, in the linear model, if $\E X_1=0$ and $\frac{\partial L(Y_1,\<\beta^*,X_1\>)}{\partial \mu}$ is a noise independent of $X_1$ with mean zero and variance $\sigma^2$, then $\lambda_{\min}(\Lambda)\ge\sigma^2\lambda_{\min}(\mathbf{\Sigma})\ge \sigma^2\kappa_l$ and $\lambda_{\min}(J)/\lambda_{\max}(J)=\lambda_{\min}(\mathbf{\Sigma})/\lambda_{\max}(\mathbf{\Sigma})\ge c$ by \cref{lasso:covariate}. Let $SJ^{-1}\overset{\triangle}{=}(\mathbf{a}_1,\dots,\mathbf{a}_p)^\top$. Since $\min_{j\in\brac{p}}\ltwo{\mathbf{a}_j}^2\ge\lambda_{\min}(SJ^{-2}S)\ge c^2$ and $\max_{j\in\brac{p}}\E\brac{\frac{\partial^2 L(Y_1,\<\beta^*,X_1\>)}{\partial \mu^2}X_{1,j}^2}\lesssim\norm{X_1}_{\psi_2}^2$, then $\min_{j\in\brac{p}}\mathbf{a}_j^\top\Lambda\mathbf{a}_j\ge c^3$ and $\max_{j\in\brac{p}}\norm{\mathbf{a}_j^\top\frac{\partial L(Y_1,\<\beta^*,X_1\>)}{\partial \mu}X_1}_{\psi_1}\lesssim \max_{j\in\brac{p}}\ltwo{\mathbf{a}_j}\lesssim \matrixnorm{J^{-1}}_1$ by \cref{lasso:lossfunc,lasso:covariate,inference: hessian,inference:Lambda}. By Gaussian approximation theorem \citep[see Theorem 2.1 in][]{chernozhuokov2022improved}, we have
\begin{align}\label{eq: bahadur approximation}
\sup_{t\ge0}\absBig{\P\parenBig{\normBig{SJ^{-1}\frac{1}{\sqrt{n}}\sumn\frac{\partial L(Y_i,\<\beta^*,X_i\>)}{\partial \mu}X_i}_{\mc{G}}\leq t}-\P\parenBig{\norm{\mathbf{Z}}_{\mc{G}}\leq t}}\lesssim \parenBig{\frac{\log^{5}(np)}{n}}^{1/4},
\end{align}
where $\mathbf{Z}$ is a Gaussian random variable with mean zero and covariance $SJ^{-1}\Lambda J^{-1}S$, $\mc{G}$ is a subset of $\brac{p}$ and $\norm{\cdot}_{\mc{G}}$ is a semi-norm defined by $\norm{a}_{\mc{G}}=\max_{j\in\mc{G}}\abs{a_j}$ for $a\in\R^{p}$. Similar to \cite{cai2025statistical}, here, the purpose of introducing the diagonal matrix $S$ in \eqref{eq: bahadur approximation} is to obtain simultaneous confidence intervals with varying length. Denote the estimate of $S$ as
\begin{align*}
    \hat{S}=\text{diag}\setBig{\dn\sumn\frac{\partial^2 L(Y_i,\<\hat{\beta}_{MA}(\hat{w}),X_i\>)}{\partial\mu^2}X_iX_i^\top}.
\end{align*}
By anti-concentration inequality \citep[see, for example, Theorem 3 in ][]{chernozhukov2015comparison} and \cref{inference: bahadur representation}, we have the following Gaussian approximation result.

\begin{theorem}\label{Gaussian approximation}
Under \cref{inference:Lambda} and the conditions in \cref{inference: bahadur representation}, for any $\mc{G}\subseteq\brac{p}$, if $\eta_n=o(\log^{-1/2}(p))$ and $\log^5(np)=o(n)$, we have
    \begin{align*}
\sup_{t\ge0}\absBig{\P\parenBig{\norm{\sqrt{n}\hat{S}(\tilde{\beta}(\hat{w})-\beta^*)}_{\mc{G}}\leq t}-\P\parenBig{\norm{\mathbf{Z}}_{\mc{G}}\leq t}}=o(1).
\end{align*}
\end{theorem}

This theorem shows that the distribution of $\sqrt{n}\hat{S}(\tilde{\beta}(\hat{w})-\beta^*)$ can be approximated by the Gaussian random variable $\mathbf{Z}$ uniformly. However, the covariance of $\mathbf{Z}$ is unknown. We consider approximating the distribution of $\mathbf{Z}$ by using multiplier bootstrap. Let 
\begin{align*}
      \hat{\mathbf{Z}}^{(boot)}=\hat{S}\hat{W}\frac{1}{\sqrt{n}}\sumn\frac{\partial L(Y_i,\<\hat{\beta}_{MA}(\hat{w}),X_i\>)}{\partial \mu}X_ie_i,
\end{align*}
where $e_1,\dots,e_n$ are $N(0,1)$ standard random variables independent of the data $\set{(Y_i,X_i)}_{i=1}^n$. By Gaussian comparison result \citep[see Proposition 2.1 in ][]{chernozhuokov2022improved}, we have the following bootstrap approximation result.

\begin{theorem}\label{bootstrap approximation}
Under \cref{inference:Lambda} and the conditions in \cref{inference: bahadur representation}, for any $\mc{G}\subseteq\brac{p}$, if $\eta_n=o(\log^{-3/2}(p))$ and $\log^5(np)\log^4(p)=o(n)$, we have
\begin{align*}
    \sup_{t\ge0}\absBig{\P\parenBig{\norm{\hat{\mathbf{Z}}^{(boot)}}_{\mc{G}}\leq t\Big|\set{(Y_i,X_i)}_{i=1}^n}-\P\parenBig{\norm{\mathbf{Z}}_{\mc{G}}\leq t}}=o_p(1).
\end{align*}
\end{theorem}

\section{Practical Implementation}\label{sec:implementation}
In this section, we propose a fast greedy model averaging algorithm (FGMA) to solve the simplex-constrained problem \eqref{eq:optimization problem}, and discuss the implementation details for prediction and post-averaging inference.
\subsection{Fast Greedy Model Averaging Algorithm}
In order to obtain the weight estimator, we need to minimize the cross-validation criterion \eqref{eq:optimization problem}, which is a simplex-constrained optimization problem. If the loss function that we considered is a $\ell_2$ loss, then the objective function is a quadratic function, which can be efficiently solved using ADMM \citep{boyd2011distributed} algorithms. However, when $CV$ is not a quadratic function of $w$, there is no closed-form solution for the $w$-minimization step in the ADMM algorithm. An alternative popular method in the model averaging literature is the greedy model averaging algorithm introduced in \cite{dai2012deviation}, \cite{he2020functional} and \cite{he2023rank}. The greedy model averaging algorithm in our framework is described in \cref{greedy model averaging}, which can be viewed as a variant of Frank-Wolfe algorithms. In the classical Frank-Wolfe algorithm \citep[a.k.a. conditional gradient descent; see, for example, GMA-1 algorithm in][]{dai2012deviation}, in each step, the algorithm computes the steepest descent direction:
\begin{align*}
    G^{(N)}=&~\argmin_{G\in\mc{W}^K}\<\nabla CV(\hat{w}^{(N-1)}),G^{(N)}\>,\\
    \hat{w}^{(N)}=&~(1-\alpha_N)\hat{w}^{(N-1)}+\alpha_NG^{(N)},
\end{align*}
where $\hat{w}^{(N)}$ is the weight after the $N$-th iteration and $\alpha_N$ is the step size (learning rate) in the $N$-th step. The greedy model averaging algorithm described in \cref{greedy model averaging} and the GMA-0 algorithm proposed in \cite{dai2012deviation} minimize the natural convex upper bound $CV(w)$ of the first-order approximation $CV(\hat{w}^{(N-1)})+\<\nabla CV(\hat{w}^{(N-1)}),w-\hat{w}^{(N-1)}\>$ over the finite set $(1-\alpha_N)\hat{w}^{(N-1)}+\alpha_N\set{e^{(1)},\dots,e^{(K)}}$ in each step, where $\set{e^{(k)}; k=1,\dots,K}$ is the canonical basis in $\R^K$, which has only linear time complexity. Unfortunately, these algorithms not only fail to guarantee a monotonically decreasing optimization trajectory but also typically exhibit slow convergence rates.

To address these issues, inspired by majorization-minimization (MM) \citep{hunter2004tutorial} technique and proximal gradient method \citep{parikh2014proximal}, we locally majorize $CV$ by an isotropic quadratic function:
\begin{align}\label{second-order-approx}
CV(w)\leq CV(w_0)+\<\nabla CV(w_0),w-w_0\>+\frac{L_{CV}}{2}\ltwo{w-w_0}^2\overset{\triangle}{=}G_L(w|w_0),
\end{align}
where $L_{CV}$ is some upper bound of the operator norm of the Hessian $\nabla^2 CV(w_0)$. Such a technique is also used in \cite{fan2018lamm,fan2023communication} and \cite{liu2025communication}. Define an operator $\Pi_{\mc{W}^K}(u)=\argmin_{w\in\mc{W}^K}\ltwo{w-u}$, which is a Euclidean projection onto the simplex $\mc{W}^K$. Let
\begin{align*}
    p_L(w_0)\overset{\triangle}{=}\Pi_{\mc{W}^K}\setBig{w_0-\frac{1}{L}\nabla CV(w_0)}=\argmin_{w\in\mc{W}^K}G_L(w|w_0).
\end{align*}
Thus, we can update $\hat{w}^{(N-1)}$ by $\hat{w}^{(N)}=(1-\alpha_N)\hat{w}^{(N-1)}+\alpha_N p_{L_{CV}}(\hat{w}^{(N-1)})$ with a step size $\alpha_N$ close to one in each step, which gives us a decreasing optimization trajectory.

Inspired by \cite{nesterov2013gradient}, \cite{beck2009fast} and \cite{fan2018lamm}, we can accelerate this algorithm, and the unknown parameter $L_{CV}$ can be approximated by backtracking line search. This fast greedy model averaging (FGMA) algorithm is summarized in \cref{alg:FGMA}. Given a vector $z\in\R^K$, the projection $\Pi_{\mc{W}^K}(z)$ can be obtained with a linear time complexity \citep[see Figure 2 in][]{duchi2008efficient}, which is the same as the complexity of each step in simplex-constrained Frank-Wolfe-type algorithms. We give the convergence analysis for our proposed fast greedy model averaging algorithm in \cref{thm:FGMA}. Compared with the greedy model averaging algorithm described in \cref{greedy model averaging}, \cref{alg:FGMA} matches the optimal convergence rate of first-order algorithm with smooth convex objective function (without the strong convexity condition), that is, $CV(\hat{w}^{(N)})/n=O(N^{-2})$ \citep[see, Theorem 2.1.7 in][]{nesterov2018lectures}.

We then introduce a stopping criterion for the fast greedy model averaging algorithm. Note the exact solution $\hat{w}$ satisfies the first-order optimality condition: for any $w\in\mc{W}^K$,
\begin{align*}
    \<\nabla CV(\hat{w}),w-\hat{w}\>\ge0.
\end{align*}
Thus, for a pre-specified optimization tolerance $\eps>0$, we can stop the algorithm when $\<\nabla CV(\hat{w}^{(N)}),w-\hat{w}^{(N)}\>\ge-\eps$ for any $w\in\mc{W}^K$, which is equivalent to the condition $\min_k\set{\nabla CV(\hat{w}^{(N)})}_k\ge\<\nabla CV(\hat{w}^{(N)}),\hat{w}^{(N)}\>-\eps$. For the initial point $\hat{w}^{(0)}$ of this algorithm, it is natural to choose $\hat{w}^{(0)}=e^{(\varsigma)}$, where $\varsigma=\argmin_{k\in\brac{K}}CV(e^{(k)})$. More details can be found in \cref{alg:FGMA}.

\begin{algorithm}
\caption{Fast greedy model averaging}
\label{alg:FGMA}
\begin{algorithmic}[1]
\INPUT{Initial values $L_0>0$, $A_1=1$, $z^{(1)}=\hat{w}^{(0)}=e^{(\varsigma)}$ (where $\varsigma=\argmin_{k\in\brac{K}}CV(e^{(k)})$), tolerance $\eps>0$, backtracking parameter $\gamma>1$.}
\FORALL{$N=1,2,\dots$ until $\min_k\set{\nabla CV(\hat{w}^{(N)})}_k\ge\<\nabla CV(\hat{w}^{(N)}),\hat{w}^{(N)}\>-\eps$}
\STATE Find the smallest nonnegative integers $j_N$ such that $CV(p_{\tilde{L}}(z^{(N)}))\leq G_{\tilde{L}}(p_{\tilde{L}}(z^{(N)})|z^{(N)})$ with $\tilde{L}=\gamma^{j_N}L_{N-1}$, where $p_{\tilde{L}}(z^{(N)})=\Pi_{\mc{W}^K}(z^{(N)}-\frac{1}{{\tilde{L}}}\nabla CV(z^{(N)}))$.
\STATE Set $L_{N}=\tilde{L}$, $\hat{w}^{(N)}=p_{L_N}(z^{(N)})$, $A_{N+1}=\frac{1+\sqrt{1+4A_{N}^2}}{2}$ and $z^{(N+1)}=\hat{w}^{(N)}+\frac{A_N-1}{A_{N+1}}(\hat{w}^{(N)}-\hat{w}^{(N-1)})$.
\ENDFOR
\OUTPUT{The approximate weight estimator $\hat{w}^{(N)}$.}
\end{algorithmic}
\end{algorithm}

\begin{theorem}\label{thm:FGMA}
    For any $N\in\N$, we have
    \begin{align*}
        \frac{CV(\hat{w}^{(N)})}{n}\leq\frac{CV(\hat{w})}{n}+\frac{2\gamma L_{CV}\ltwo{\hat{w}^{(0)}-\hat{w}}^2}{(N+1)^2},
    \end{align*}
    where $L_{CV}=\max_{\lone{z}\leq5}\matrixnorm{\nabla^2 CV(z)/n}_2.$
\end{theorem}

\vspace{0.2cm}
\subsection{Implementation Details}

\subsubsection{Prediction}\label{sec:pred}
In our framework, a general penalty function can be used to obtain an estimator for each candidate model. If $r_\lambda$ is a folded concave regularizer, we propose to use the local linear approximation \citep{zou2008one,fan2014strong} algorithm to solve penalized empirical risk minimization problems. This algorithm can be found in Algorithm 1 of \cite{fan2014strong}, and we use the three-step solution in practice. Thus, estimators based on concave penalties can be obtained by iteratively computing the solution of the weighted Lasso efficiently using some R packages such as \texttt{glmnet}.

Now, we introduce more details about the preparation of candidate models. Analogous to most high-dimensional model averaging literature, we rank the importance of each predictor in descending order. In \cite{ando2014model,ando2017weight}, \cite{xie2021model}, \cite{wang2023jackknife} and \cite{he2023rank}, they use marginal information to rank the importance of each predictor. \cite{zhang2020parsimonious} proposed ranking predictors based on their order of entry into the adaptive Lasso solution path in the LARS algorithm. Here, we sort the predictors using a two-step procedure. We first compute an initial estimator $\hat{\beta}=\argmin_{\beta\in\R^p}\dn\sumn L(Y_i,\<\beta,X_i\>)+r_{\lambda_n}(\beta)$ and rank the predictors in the support of the initial estimator $\hat{\beta}$ by the absolute values of the estimated coefficients according to the initial estimator, where $\lambda_n$ is obtained by cross-validation. Then we order the remaining predictors by their marginal utilities \citep[e.g. the $p$-values of marginal regression; see,][]{ando2014model,ando2017weight}. Given the pre-fixed $K_{ne}$ and $d_2$, we choose $d_1=2\lceil \lzero{\hat{\beta}}/K_{ne}\rceil, p_0=K_{ne}d_1$ and $K=K_{ne}+\lfloor(p-p_0)/d_2\rfloor.$ We can construct candidate models using the ranked predictors: $\mc{A}_k=[kd_1]$ for $k\in\brac{K_{ne}}$ and $\mc{A}_k=\set{p_0+(k-K_{ne}-1)d_2+1,\dots,p_0+(k-K_{ne})d_2}$ for $k\in\set{K_{ne}+1,\dots,K}$. Therefore, the model size of the largest candidate model in the nested part satisfies $p_0\in[2\lzero{\hat{\beta}},3\lzero{\hat{\beta}})$, which can include more important variables in the nested part.

\begin{algorithm}[htbp]
\caption{High-dimensional model averaging (HDMA) via cross-validation}
\label{alg:crossvalidation}
\begin{algorithmic}[1]
\INPUT{Data $\mc{D}=\set{(Y_i,X_i)}_{i=1}^n$, candidate model set $\set{\mc{A}_k}_{k\in\brac{K}}$, tuning parameters $\set{\lambda_{n,k}}_{k\in\brac{K}}$, number of iterations $N$ in \cref{alg:FGMA}, regularizer $r_\lambda(\cdot)$.}

\FORALL{$k\in\brac{K}$}
\STATE Compute $\hat{\beta}_{(k)}=\argmin_{\beta\in\R^p, \beta_{\mc{A}^c_k}=0}\dn\sumn L(Y_i,\langle \beta,X_i\rangle)+r_{\lambda_{n,k}}(\beta)$.
\ENDFOR
\STATE Split the data into $J$ folds of equal size with the index sets $\mc{I}_1,\dots,\mc{I}_J$.
\FORALL{$m\in\brac{J}, k\in\brac{K}$}
    \STATE Compute $\hat{\beta}_{(k)}^{\brac{-m}}=\argmin_{\beta\in\R^p, \beta_{\mc{A}^c_k}=0}\frac{1}{n_{-1}}\sum_{i\in\brac{n}\backslash \mc{I}_m} L(Y_i,\langle \beta,X_i\rangle)+r_{\lambda_{n,k}}(\beta)$ using the out-of-$m$-fold data $\mc{D}\backslash\set{(Y_i,X_i)}_{i\in\mc{I}_m}$.
\ENDFOR
\STATE Compute $\hat{w}^{(N)}$ using the FGMA algorithm described in \cref{alg:FGMA}.
\OUTPUT{The approximate MA estimator $\hat{\beta}_{MA}(\hat{w}^{(N)})=\hat{B}\hat{w}^{(N)}$.}
\end{algorithmic}
\end{algorithm}

Next, we discuss how to choose the tuning parameters for each candidate model. In our work, the number of tuning parameters is the same as that of the candidate models. It is time-consuming to choose the tuning parameter for each candidate model separately. To address this issue, we propose an adaptive method for selecting the tuning parameters. Recall that $\lambda_n$ is used to calculate $\hat{\beta}$. Note that if we calculate the largest candidate model in the nested part $\hat{\beta}_{(K_{ne})}$ using the tuning parameter $\lambda_n$, then $\hat{\beta}_{(K_{ne})}$ gives us the same solution as the initial estimator $\hat{\beta}$. For the $k$-th candidate model, we recommend using a uniform choice $\lambda_{n,k}=\lambda_{n,k}^{\brac{-m}}=\sqrt{\frac{\log \abs{\mc{A}_k}}{\log p_0}}\lambda_n$, $m\in\brac{J}$ in practice. 

Our model averaging procedure is presented in \cref{alg:crossvalidation}.

%\vspace{0.2cm}
\subsubsection{Post-Averaging Simultaneous Inference}

Following \cref{Gaussian approximation,bootstrap approximation}, if we assume that the true coefficients are sparse and that there is at least a correct model in the candidate model set,
%(for example, in practice, one may also consider a nested model structure: $\mc{A}_1\subseteq\mc{A}_2\subseteq\cdots\subseteq\mc{A}_K=\brac{p}$)
we can construct simultaneous confidence intervals by Gaussian multiplier bootstrap, which is described in \cref{alg:inference}. 

\begin{algorithm}[htbp]
\caption{Multiplier bootstrap procedures for post-averaging simultaneous confidence intervals}
\label{alg:inference}
\begin{algorithmic}[1]
\INPUT{Data $\mc{D}=\set{(Y_i,X_i)}_{i=1}^n$, the model averaging estimator $\hat{\beta}_{MA}(\hat{w})$, the estimated inverse-Hessian matrix $\hat{W}$, the index set $\mc{G}$, the bootstrap replications 
$B$, the significance level $\alpha$.}
\STATE Compute the one-step debiased estimator $\tilde{\beta}(\hat{w})=\hat{\beta}_{MA}(\hat{w})-\hat{W}\dn\sumn\frac{\partial L(Y_i,\<\hat{\beta}_{MA}(\hat{w}),X_i\>)}{\partial \mu}X_i$ and $\hat{S}=\text{diag}\parenBig{\dn\sumn\frac{\partial^2L(Y_i,\<\hat{\beta}_{MA}(\hat{w}),X_i\>)}{\partial \mu^2}X_iX_i^\top}$.
\FORALL{$b=1, \dots, B$}
\STATE Generate standard normal random variables $\set{e_{(b),i}}_{i=1}^n$.
\STATE Compute $\hat{Z}^{(boot)}_{b}=\hat{S}\hat{W}\frac{1}{\sqrt{n}}\sumn\frac{\partial L(Y_i,\<\hat{\beta}_{MA}(\hat{w}),X_i\>)}{\partial \mu}X_ie_{(b),i}$.
\ENDFOR
\STATE Compute the upper-$\alpha$ quantile of $\set{\norm{\hat{Z}_{(b)}^{(boot)}}_{\mc{G}}}_{b=1}^B$, denoted by $\hat{Q}_{1-\alpha}$.
\OUTPUT{The post-averaging simultaneous confidence intervals 
$$\Big(\tilde{\beta}_j(\hat{w})-(\hat{S}_{jj})^{-1}\hat{Q}_{1-\alpha}/\sqrt{n},\tilde{\beta}_j(\hat{w})+(\hat{S}_{jj})^{-1}\hat{Q}_{1-\alpha}/\sqrt{n}\Big), j\in\mc{G}.$$}
\end{algorithmic}
\end{algorithm}

Given the model averaging estimator $\hat{\beta}_{MA}(\hat{w})$, we first calculate the one-step debiased estimator $\tilde{\beta}(\hat{w})=\hat{\beta}_{MA}(\hat{w})-\hat{W}\dn\sumn\frac{\partial L(Y_i,\<\hat{\beta}_{MA}(\hat{w}),X_i\>)}{\partial \mu}X_i$, where the matrix $\hat{W}$ is obtained using the R package \texttt{flare} with the default choice of the tuning parameter, as suggested in \cite{cai2025statistical}. For the prefixed bootstrap replications $B$, index set $\mc{G}$ and significance level $\alpha\in(0,1)$, we compute the bootstrap samples $\hat{Z}^{(boot)}_{b}=\hat{S}\hat{W}\frac{1}{\sqrt{n}}\sumn\frac{\partial L(Y_i,\<\hat{\beta}_{MA}(\hat{w}),X_i\>)}{\partial \mu}X_ie_{(b),i}, b\in\brac{B}$ and compute the upper-$\alpha$ quantile $\hat{Q}_{1-\alpha}$ of $\set{\norm{\hat{Z}_{(b)}^{(boot)}}_{\mc{G}}}_{b=1}^B$. Then, the $1-\alpha$ post-averaging simultaneous confidence intervals of $\beta^*$ is given by
\begin{align*}
    \Big(\tilde{\beta}_j(\hat{w})-(\hat{S}_{jj})^{-1}\hat{Q}_{1-\alpha}/\sqrt{n},\tilde{\beta}_j(\hat{w})+(\hat{S}_{jj})^{-1}\hat{Q}_{1-\alpha}/\sqrt{n}\Big), ~j\in\mc{G}.
\end{align*}

\section{Simulation}\label{sec: simulation}
In this section, we perform numerical studies to assess the prediction and inference performance of the proposed high-dimensional model averaging method.
\subsection{Prediction}
In this subsection, we compare the performance of our approach with both model selection methods and model averaging methods in the settings of linear regression and logistic regression. In each setting, we generate $X_i$'s independently from a $p$-dimensional normal distribution $N(\mathbf{0},\mathbf{\Sigma})$. We consider the following two scenarios for the correlation structure of $X_i$'s.
\begin{enumerate}
    \item AR(1) design: $\mathbf{\Sigma}_{j,j^\prime}=0.5^{\abs{j-j^\prime}}$ for $j,j^\prime=1,2,\dots,p$;
    \item Band design: $\mathbf{\Sigma}_{j,j^\prime}=0.5^{\abs{j-j^\prime}}\indic{\abs{j-j^\prime}\leq 1}$ for $j,j^\prime=1,2,\dots,p$.
\end{enumerate} We also consider different settings with varying sample size $n\in\set{100,200}$ and dimension of covariates $p\in\set{1000,2000}$.

The proposed method described in \cref{alg:crossvalidation} is denoted by HDMA. To implement our method, we consider three common regularizers: Lasso, SCAD and MCP. For simplicity and time-saving, in our implementation of the proposed methods, we set the number of candidate models in the nested part $K_{ne}=4$, the size of candidate models in the nonnested part $d_2=10$ and the number of folds $J=5$. We use these parameters in both \cref{sec: simulation} and \cref{sec: real data}. To compare the prediction performance of our method with those of other methods, we generate out-of-sample observations that do not vary across the replications in each setting.

%%%%%%%%%%%%%%%%%%%%%%%%%%%%%%%%%%%%%%%%
\subsubsection{Linear Regression}
The data are generated by $Y_i=\sum_{j=1}^p\mathrm{x}_{i,j}\beta_{j}+\eps_i$, $i=1,\dots,n$. We take the distribution of $\varepsilon_i$ to be a normal distribution with mean zero and standard deviation $0.5$. The following three different coefficient settings are considered to assess the performance of the proposed methods.
\begin{enumerate}
    \item Sparse design: $\bm{\beta}=(\mathbbm{1}_5^\top,0.2\mathbbm{1}_{10}^\top,\mathbbm{1}_5^\top,0,0,\dots)^\top$, where $\mathbbm{1}_d$ is a $d$-dimensional vector of ones;
    \item Polynomially decaying design: $\beta_j=5j^{-2},j=1,\dots,p$;
    \item Exponentially decaying design: $\beta_j=5e^{-0.3j},j=1,\dots,p$.
\end{enumerate}

\begin{table}[tbp]
\centering
\resizebox{0.92\linewidth}{!}{\footnotesize
\begin{tabularx}{\textwidth}{p{16pt}p{10pt}p{10pt}XXXXXXXXXXXX}
\toprule
{\footnotesize$\mathbf{\Sigma}$} &{\footnotesize $p$} &{\footnotesize $n$} & &{\footnotesize Lasso }&{\footnotesize PLasso }&{\footnotesize SCAD }&{\footnotesize MCP }& {\footnotesize ENet }& {\footnotesize AnL }&{\footnotesize PMA }&{\footnotesize HDMA (Lasso) }&{\footnotesize HDMA (SCAD)} &{\footnotesize HDMA (MCP)} \\ \hline\\[-8pt]
 & & &  & \multicolumn{10}{c}{$\bm{\beta}=(\mathbbm{1}_5^\top,0.2\mathbbm{1}_{10}^\top,\mathbbm{1}_5^\top,0,0,\dots)^\top$}  \\
AR(1)   & 1000 & 100 & Mean    & 0.292 & 0.314  & 0.317 & 0.334 & 0.330 & 0.747 & 0.351 & \textbf{0.264} & 0.275          & 0.265          \\
      &      &     & SD      & 0.068 & 0.141  & 0.192 & 0.242 & 0.088 & 0.479 & 0.082 & 0.060          & 0.041          & 0.044          \\
      &      & 200 & Mean    & 0.170 & 0.177  & 0.207 & 0.207 & 0.179 & 0.290 & 0.239 & \textbf{0.162} & 0.172          & 0.165          \\
      &      &     & SD      & 0.013 & 0.029  & 0.020 & 0.016 & 0.015 & 0.114 & 0.020 & 0.013          & 0.017          & 0.016          \\
      & 2000 & 100 & Mean    & 0.368 & 0.353  & 0.451 & 0.425 & 0.442 & 1.171 & 0.416 & 0.329          & 0.295          & \textbf{0.279} \\
      &      &     & SD      & 0.086 & 0.147  & 0.499 & 0.441 & 0.115 & 0.728 & 0.099 & 0.107          & 0.048          & 0.046          \\
      &      & 200 & Mean    & 0.195 & 0.202  & 0.220 & 0.216 & 0.206 & 0.382 & 0.246 & \textbf{0.183} & 0.195          & 0.186          \\
      &      &     & SD      & 0.017 & 0.034  & 0.019 & 0.019 & 0.020 & 0.145 & 0.023 & 0.017          & 0.018          & 0.019          \\
Band    & 1000 & 100 & Mean    & 0.396 & 0.378  & 0.802 & 0.800 & 0.474 & 0.891 & 0.425 & 0.365          & 0.221          & \textbf{0.217} \\
      &      &     & SD      & 0.116 & 0.182  & 0.193 & 0.211 & 0.146 & 0.361 & 0.115 & 0.129          & 0.064          & 0.062          \\
      &      & 200 & Mean    & 0.197 & 0.203  & 0.240 & 0.180 & 0.209 & 0.460 & 0.259 & 0.181          & 0.154          & \textbf{0.153} \\
      &      &     & SD      & 0.017 & 0.030  & 0.125 & 0.065 & 0.020 & 0.087 & 0.027 & 0.015          & 0.013          & 0.011          \\
      & 2000 & 100 & Mean    & 0.543 & 0.492  & 0.854 & 0.864 & 0.704 & 1.466 & 0.586 & 0.505          & 0.256          & \textbf{0.253} \\
      &      &     & SD      & 0.199 & 0.312  & 0.176 & 0.217 & 0.277 & 0.961 & 0.209 & 0.222          & 0.095          & 0.100          \\
      &      & 200 & Mean    & 0.204 & 0.205  & 0.286 & 0.192 & 0.224 & 0.537 & 0.249 & 0.184          & 0.148          & \textbf{0.147} \\
      &      &     & SD      & 0.028 & 0.041  & 0.154 & 0.109 & 0.035 & 0.112 & 0.032 & 0.025          & 0.013          & 0.013        \\ \hline\\[-8pt]
 &  &&  & \multicolumn{10}{c}{$\beta_j=5j^{-2},j=1,\dots,p$}  \\
AR(1)    & 1000 & 100 & Mean    & 0.188 & 0.197  & 0.279 & 0.261 & 0.213 & 0.296 & 0.237 & \textbf{0.181} & 0.186          & 0.187          \\
      &      &     & SD      & 0.021 & 0.036  & 0.044 & 0.039 & 0.029 & 0.081 & 0.027 & 0.022          & 0.021          & 0.025          \\
      &      & 200 & Mean    & 0.160 & 0.170  & 0.253 & 0.233 & 0.169 & 0.207 & 0.230 & \textbf{0.157} & 0.161          & 0.161          \\
      &      &     & SD      & 0.009 & 0.017  & 0.024 & 0.025 & 0.011 & 0.030 & 0.019 & 0.009          & 0.011          & 0.010          \\
      & 2000 & 100 & Mean    & 0.197 & 0.210  & 0.256 & 0.241 & 0.227 & 0.374 & 0.239 & \textbf{0.187} & 0.193          & 0.193          \\
      &      &     & SD      & 0.023 & 0.031  & 0.031 & 0.038 & 0.032 & 0.104 & 0.032 & 0.022          & 0.022          & 0.024          \\
      &      & 200 & Mean    & 0.160 & 0.171  & 0.267 & 0.245 & 0.171 & 0.231 & 0.211 & \textbf{0.156} & 0.158          & 0.157          \\
      &      &     & SD      & 0.010 & 0.019  & 0.033 & 0.037 & 0.014 & 0.046 & 0.016 & 0.011          & 0.011          & 0.011          \\
Band     & 1000 & 100 & Mean    & 0.196 & 0.203  & 0.270 & 0.236 & 0.220 & 0.451 & 0.243 & 0.189          & 0.187          & \textbf{0.186} \\
      &      &     & SD      & 0.021 & 0.031  & 0.051 & 0.025 & 0.027 & 0.163 & 0.022 & 0.023          & 0.024          & 0.025          \\
      &      & 200 & Mean    & 0.162 & 0.174  & 0.261 & 0.233 & 0.168 & 0.291 & 0.231 & \textbf{0.158} & 0.159          & \textbf{0.158} \\
      &      &     & SD      & 0.008 & 0.016  & 0.036 & 0.017 & 0.010 & 0.091 & 0.018 & 0.008          & 0.008          & 0.010          \\
      & 2000 & 100 & Mean    & 0.199 & 0.212  & 0.253 & 0.221 & 0.242 & 0.624 & 0.238 & 0.188          & 0.182          & \textbf{0.181} \\
      &      &     & SD      & 0.028 & 0.038  & 0.060 & 0.041 & 0.042 & 0.241 & 0.030 & 0.029          & 0.022          & 0.025          \\
      &      & 200 & Mean    & 0.168 & 0.179  & 0.250 & 0.232 & 0.184 & 0.325 & 0.221 & 0.163          & \textbf{0.160} & \textbf{0.160} \\
      &      &     & SD      & 0.012 & 0.020  & 0.026 & 0.017 & 0.016 & 0.098 & 0.017 & 0.011          & 0.009          & 0.011          \\ \hline\\[-8pt]
 &  & &  & \multicolumn{10}{c}{$\beta_j=5e^{-0.3j},j=1,\dots,p$} \\
AR(1)     & 1000 & 100 & Mean    & 0.219 & 0.219  & 0.439 & 0.368 & 0.241 & 0.780 & 0.274 & 0.206          & 0.211          & \textbf{0.204} \\
      &      &     & SD      & 0.025 & 0.039  & 0.119 & 0.091 & 0.036 & 0.526 & 0.035 & 0.023          & 0.031          & 0.028          \\
      &      & 200 & Mean    & 0.168 & 0.178  & 0.371 & 0.316 & 0.174 & 0.341 & 0.242 & \textbf{0.163} & 0.165          & 0.164          \\
      &      &     & SD      & 0.011 & 0.018  & 0.086 & 0.062 & 0.013 & 0.154 & 0.018 & 0.010          & 0.012          & 0.010          \\
      & 2000 & 100 & Mean    & 0.233 & 0.225  & 0.481 & 0.390 & 0.267 & 1.136 & 0.276 & 0.216          & 0.209          & \textbf{0.203} \\
      &      &     & SD      & 0.039 & 0.042  & 0.126 & 0.080 & 0.051 & 0.728 & 0.050 & 0.043          & 0.044          & 0.038          \\
      &      & 200 & Mean    & 0.168 & 0.173  & 0.370 & 0.324 & 0.178 & 0.443 & 0.222 & 0.162          & 0.157          & \textbf{0.153} \\
      &      &     & SD      & 0.013 & 0.021  & 0.070 & 0.063 & 0.015 & 0.225 & 0.018 & 0.014          & 0.014          & 0.013          \\
Band     & 1000 & 100 & Mean    & 0.259 & 0.243  & 0.399 & 0.323 & 0.303 & 1.354 & 0.305 & 0.236          & \textbf{0.198} & 0.200          \\
      &      &     & SD      & 0.045 & 0.048  & 0.103 & 0.052 & 0.060 & 0.802 & 0.052 & 0.044          & 0.026          & 0.028          \\
      &      & 200 & Mean    & 0.175 & 0.177  & 0.342 & 0.303 & 0.186 & 0.578 & 0.239 & 0.167          & \textbf{0.153} & 0.154          \\
      &      &     & SD      & 0.013 & 0.020  & 0.043 & 0.035 & 0.015 & 0.312 & 0.020 & 0.014          & 0.009          & 0.011          \\
      & 2000 & 100 & Mean    & 0.275 & 0.250  & 0.379 & 0.316 & 0.336 & 1.754 & 0.312 & 0.243          & \textbf{0.199} & \textbf{0.199} \\
      &      &     & SD      & 0.045 & 0.051  & 0.089 & 0.053 & 0.070 & 1.196 & 0.052 & 0.047          & 0.030          & 0.036          \\
      &      & 200 & Mean    & 0.182 & 0.191  & 0.335 & 0.291 & 0.195 & 0.690 & 0.233 & 0.174          & 0.161          & \textbf{0.159} \\
      &      &     & SD      & 0.014 & 0.022  & 0.047 & 0.026 & 0.017 & 0.359 & 0.021 & 0.014          & 0.010          & 0.010\\  \bottomrule
\end{tabularx}}
\caption{Out-of-sample prediction errors (PE$_1$) of estimators for linear regression (The smallest mean is bolded for each setting)}
\label{sim: simulaton-lin}
\end{table}

For comparison, in addition to our approach, we also consider the following methods: (1) Lasso \citep{tibshirani1996regression}; (2) MCP \citep{zhang2010nearly}; (3) SCAD \citep{fan2001variable}; (4) AnL \citep[the leave-one-out cross-validation method in][]{ando2014model}; (5) PMA \citep{zhang2020parsimonious}; (6) Elastic net \citep{zou2005regularization} with ratio 0.5, denoted by ENet; (7) OLS Post-fit Lasso \citep{belloni2013least}, denoted by PLasso. The tuning parameters of the penalized methods are selected by using a 10-fold cross-validation. All penalized methods are implemented by the R package \texttt{glmnet} or \texttt{ncvreg} with default settings. We evaluate the performance of our method and the above methods with respect to the out-of-sample prediction error 
\begin{align}\label{eq:prediction error1}
    \text{PE}_1=\frac{1}{R}\sum_{r=1}^{R}\frac{1}{n_{test}}\sum_{i=1}^{n_{test}} (Y^o_i-\<\hat{\beta}^{(r)},X_i^o\>)^2,
\end{align}
where $R=100$ is the number of replications, $n_{test}=1000$ is the number of test samples, $\hat{\beta}^{(r)}$ denotes the estimator of the $r$-th replication and $\set{(Y^o_i,X_i^o)}_{i=1}^{n_{test}}$ is the set of out-of-sample observations. 

The simulation results with 100 replications are shown in \cref{sim: simulaton-lin}. Simulations show that our HDMA methods have very comparable performance and are uniformly superior to other model averaging and regularized methods, especially when the coefficient vector is not sparse. Moreover, the HDMA methods consistently achieve a low standard deviation similar to Lasso and significantly lower than the competing methods. The Lasso-based HDMA estimator has the best overall performance in the case where the coefficient vector is polynomially decaying. In the case where the covariance matrix has the band structure, the HDMA estimators with folded concave penalties have overall good performance.

% %%%%%%%%%%%%%%%%%%%%%%%%%%%%%%%%%%%%%%%%
\subsubsection{Logistic Regression}
We generate $Y_i$ from Bernoulli$(p_{i0})$ with $p_{i0}=1/\set{1+\exp(-\sum_{j=1}^p\mathrm{x}_{i,j}\beta_j)}, i=1,\dots,n.$ In this subsection, we consider the following coefficient settings:
\begin{enumerate}
    \item Sparse design: $\bm{\beta}=(3\mathbbm{1}_5^\top,\mathbbm{1}_{10}^\top,-0.2\mathbbm{1}_5^\top,0,0,\dots)^\top$;
    \item Polynomially decaying design: $\beta_j=5(\indic{j\leq5}+(j-5)^{-4}\indic{j>5}),j=1,\dots,p$;
    \item Exponentially decaying design: $\beta_j=5(\indic{j\leq5}+e^{-0.5(j-5)}\indic{j>5}),j=1,\dots,p$.
\end{enumerate}

For comparison, we also implement the following approaches: (1) Lasso; (2) MCP; (3) SCAD; (4) Elastic net (ENet) with ratio 0.5; (5) AnL \citep[the leave-one-out cross-validation method in][]{ando2017weight}; (6) Post Lasso, denoted by PLasso. For implementing post Lasso, we refit a Lasso using the variables selected by the first-stage Lasso, since the MLE of logistic regression might not exist \citep[see, for example,][]{candes2020phase}. The tuning parameters of the penalized methods are chosen via 10-fold cross-validation. We evaluate the performance of our method and the above methods with respect to the following out-of-sample prediction error
\begin{align}\label{eq:prediction error2}
    \text{PE}_2=\frac{1}{R}\sum_{r=1}^{R}\frac{1}{n_{test}}\sum_{i=1}^{n_{test}} \log\set{1+\exp(\<\hat{\beta}^{(r)},X_i^o\>)}-Y^o_i\<\hat{\beta}^{(r)},X_i^o\>,
\end{align}

\begin{table}[htbp]
\centering
\resizebox{0.93\linewidth}{!}{\footnotesize
\begin{tabularx}{\textwidth}{p{16pt}p{10pt}p{10pt}XXXXXXXXXXXX}
\toprule
{\footnotesize$\mathbf{\Sigma}$} &{\footnotesize $p$} &{\footnotesize $n$} & &{\footnotesize Lasso }&{\footnotesize PLasso }&{\footnotesize SCAD }&{\footnotesize MCP }& {\footnotesize ENet }& {\footnotesize AnL }&{\footnotesize HDMA (Lasso) }&{\footnotesize HDMA (SCAD)} &{\footnotesize HDMA (MCP)} \\ \hline\\[-8pt]
 & & &  & \multicolumn{9}{c}{$\bm{\beta}=(3\mathbbm{1}_5^\top,\mathbbm{1}_{10}^\top,-0.2\mathbbm{1}_5^\top,0,0,\dots)^\top$}  \\
AR(1)    & 1000 & 100 & Mean    & 0.407 & 0.528  & 0.421          & 0.440          & 0.419 & 0.592 & \textbf{0.402} & 0.448          & 0.450          \\
      &      &     & SD      & 0.040 & 0.119  & 0.043          & 0.069          & 0.036 & 0.104 & 0.047          & 0.056          & 0.054          \\
      &      & 200 & Mean    & 0.291 & 0.392  & 0.313          & 0.316          & 0.318 & 0.349 & \textbf{0.281} & 0.295          & 0.299          \\
      &      &     & SD      & 0.025 & 0.071  & 0.026          & 0.044          & 0.021 & 0.049 & 0.036          & 0.050          & 0.054          \\
      & 2000 & 100 & Mean    & 0.394 & 0.533  & 0.415          & 0.425          & 0.414 & 0.819 & \textbf{0.388} & 0.436          & 0.440          \\
      &      &     & SD      & 0.043 & 0.114  & 0.059          & 0.064          & 0.037 & 0.154 & 0.050          & 0.061          & 0.060          \\
      &      & 200 & Mean    & 0.299 & 0.406  & 0.319          & 0.321          & 0.329 & 0.500 & \textbf{0.283} & 0.306          & 0.309          \\
      &      &     & SD      & 0.025 & 0.073  & 0.035          & 0.034          & 0.020 & 0.083 & 0.030          & 0.053          & 0.060          \\
Band    & 1000 & 100 & Mean    & 0.449 & 0.603  & \textbf{0.426} & 0.439          & 0.473 & 0.709 & 0.444          & 0.441          & 0.452          \\
      &      &     & SD      & 0.042 & 0.139  & 0.055          & 0.064          & 0.037 & 0.136 & 0.056          & 0.055          & 0.058          \\
      &      & 200 & Mean    & 0.326 & 0.442  & 0.323          & 0.319          & 0.359 & 0.435 & \textbf{0.314} & 0.331          & 0.333          \\
      &      &     & SD      & 0.024 & 0.082  & 0.036          & 0.035          & 0.019 & 0.052 & 0.034          & 0.058          & 0.055          \\
      & 2000 & 100 & Mean    & 0.461 & 0.621  & \textbf{0.436} & 0.441          & 0.492 & 1.006 & 0.460          & \textbf{0.436} & 0.439          \\
      &      &     & SD      & 0.053 & 0.144  & 0.083          & 0.078          & 0.046 & 0.163 & 0.066          & 0.059          & 0.059          \\
      &      & 200 & Mean    & 0.320 & 0.431  & 0.304          & \textbf{0.299} & 0.360 & 0.575 & 0.311          & 0.324          & 0.330          \\
      &      &     & SD      & 0.025 & 0.077  & 0.039          & 0.043          & 0.022 & 0.091 & 0.037          & 0.045          & 0.044         \\ \hline\\[-8pt]
 & & &  & \multicolumn{9}{c}{$\beta_j=5(\indic{j\leq5}+(j-5)^{-4}\indic{j>5}),j=1,\dots,p$}  \\ 
AR(1)    & 1000 & 100 & Mean    & 0.295 & 0.366  & 0.310          & 0.310          & 0.342 & 0.418 & \textbf{0.258} & 0.294          & 0.301          \\
      &      &     & SD      & 0.048 & 0.087  & 0.044          & 0.056          & 0.038 & 0.109 & 0.061          & 0.068          & 0.075          \\
      &      & 200 & Mean    & 0.188 & 0.229  & 0.186          & 0.170          & 0.245 & 0.225 & 0.172          & 0.125          & \textbf{0.123} \\
      &      &     & SD      & 0.017 & 0.034  & 0.045          & 0.047          & 0.016 & 0.041 & 0.023          & 0.050          & 0.050          \\
      & 2000 & 100 & Mean    & 0.289 & 0.351  & 0.315          & 0.329          & 0.343 & 0.570 & \textbf{0.263} & 0.288          & 0.299          \\
      &      &     & SD      & 0.043 & 0.082  & 0.074          & 0.226          & 0.034 & 0.157 & 0.054          & 0.067          & 0.073          \\
      &      & 200 & Mean    & 0.208 & 0.255  & 0.206          & 0.185          & 0.268 & 0.350 & 0.196          & \textbf{0.150} & 0.155          \\
      &      &     & SD      & 0.021 & 0.040  & 0.043          & 0.042          & 0.017 & 0.069 & 0.027          & 0.069          & 0.064          \\
Band    & 1000 & 100 & Mean    & 0.324 & 0.422  & \textbf{0.270} & 0.275          & 0.382 & 0.508 & 0.290          & 0.336          & 0.275          \\
      &      &     & SD      & 0.040 & 0.101  & 0.052          & 0.060          & 0.036 & 0.141 & 0.061          & 0.665          & 0.061          \\
      &      & 200 & Mean    & 0.218 & 0.277  & 0.200          & 0.199          & 0.288 & 0.290 & 0.197          & \textbf{0.152} & 0.154          \\
      &      &     & SD      & 0.017 & 0.042  & 0.024          & 0.030          & 0.017 & 0.063 & 0.023          & 0.043          & 0.046          \\
      & 2000 & 100 & Mean    & 0.345 & 0.451  & \textbf{0.284} & 0.291          & 0.406 & 0.750 & 0.312          & 0.298          & 0.318          \\
      &      &     & SD      & 0.044 & 0.097  & 0.066          & 0.071          & 0.033 & 0.176 & 0.053          & 0.081          & 0.094          \\
      &      & 200 & Mean    & 0.232 & 0.289  & 0.214          & 0.210          & 0.302 & 0.403 & 0.215          & \textbf{0.156} & 0.162          \\
      &      &     & SD      & 0.022 & 0.045  & 0.033          & 0.036          & 0.019 & 0.101 & 0.027          & 0.044          & 0.047        \\ \hline\\[-8pt]
 & & &  & \multicolumn{9}{c}{$\beta_j=5(\indic{j\leq5}+e^{-0.5(j-5)}\indic{j>5}),j=1,\dots,p$}  \\
AR(1)    & 1000 & 100 & Mean    & 0.296 & 0.361  & 0.313          & 0.321          & 0.335 & 0.390 & \textbf{0.277} & 0.306          & 0.308          \\
      &      &     & SD      & 0.040 & 0.087  & 0.034          & 0.050          & 0.031 & 0.096 & 0.054          & 0.074          & 0.072          \\
      &      & 200 & Mean    & 0.211 & 0.263  & 0.228          & 0.212          & 0.259 & 0.244 & 0.193          & \textbf{0.163} & 0.166          \\
      &      &     & SD      & 0.022 & 0.048  & 0.050          & 0.041          & 0.019 & 0.044 & 0.027          & 0.041          & 0.048          \\
      & 2000 & 100 & Mean    & 0.332 & 0.411  & 0.348          & 0.360          & 0.368 & 0.659 & \textbf{0.313} & 0.349          & 0.364          \\
      &      &     & SD      & 0.041 & 0.089  & 0.039          & 0.074          & 0.032 & 0.143 & 0.054          & 0.065          & 0.088          \\
      &      & 200 & Mean    & 0.226 & 0.287  & 0.242          & 0.229          & 0.277 & 0.371 & 0.209          & \textbf{0.178} & 0.180          \\
      &      &     & SD      & 0.023 & 0.049  & 0.034          & 0.040          & 0.019 & 0.073 & 0.027          & 0.051          & 0.048          \\
Band     & 1000 & 100 & Mean    & 0.352 & 0.455  & \textbf{0.299} & 0.302          & 0.402 & 0.558 & 0.333          & 0.326          & 0.329          \\
      &      &     & SD      & 0.036 & 0.085  & 0.063          & 0.059          & 0.038 & 0.135 & 0.059          & 0.071          & 0.078          \\
      &      & 200 & Mean    & 0.271 & 0.352  & 0.250          & 0.243          & 0.328 & 0.365 & 0.250          & \textbf{0.216} & 0.219          \\
      &      &     & SD      & 0.023 & 0.063  & 0.029          & 0.031          & 0.019 & 0.053 & 0.031          & 0.042          & 0.045          \\
      & 2000 & 100 & Mean    & 0.364 & 0.465  & \textbf{0.291} & 0.309          & 0.417 & 0.769 & 0.347          & 0.316          & 0.329          \\
      &      &     & SD      & 0.040 & 0.093  & 0.073          & 0.149          & 0.029 & 0.174 & 0.060          & 0.083          & 0.080          \\
      &      & 200 & Mean    & 0.240 & 0.306  & 0.214          & 0.207          & 0.306 & 0.435 & 0.221          & \textbf{0.179} & 0.182          \\
      &      &     & SD      & 0.025 & 0.054  & 0.034          & 0.038          & 0.021 & 0.087 & 0.034          & 0.048          & 0.046  \\  \bottomrule
\end{tabularx}}
\caption{Out-of-sample prediction errors (PE$_2$) of estimators for logistic regression (The smallest mean is bolded for each setting)}
\label{sim: simulaton-log}
\end{table}

\noindent where $R=100$ is the number of replications, $n_{test}=1000$ is the number of test samples, $\hat{\beta}^{(r)}$ denotes the estimator of the $r$-th replication and $\set{(Y^o_i,X_i^o)}_{i=1}^{1000}$ is the set of out-of-sample observations.

\cref{sim: simulaton-log} presents the simulation results for logistic regression. It is evident that our proposed HDMA methods perform better than other approaches, except when $n=100$ and the covariance matrix is band. When the underlying coefficient vector is sparse and the covariance has the AR(1) structure, the HDMA estimator with the Lasso penalty has the best performance. In the case where the coefficient vector is not sparse and the sample size $n=200$, our HDMA estimators with folded concave penalties have the lowest prediction error among all competing methods.

%%%%%%%%%%%%%%%%%%%%%%%%%%%%%%%%%%%%%%%%%%%
\subsection{Post-Averaging Simultaneous Inference}
In this subsection, we assess the performance of our proposed post-averaging inference procedure in the settings of linear regression and logistic regression. The data generating processes and covariance structure we considered are the same as in the previous subsection. Due to theoretical constraints and computational burden, we consider the sample size $n\in\set{100, 200}$ and the dimension $p\in\set{100, 200}$. We are interested in making simultaneous inference with the three certain index sets: $\set{1,2,\dots,5}, \set{1,2,\dots,\frac{p}{5}}$ and $ \set{1,2,\dots,p}$. The underlying coefficient vector is given by $\bm{\beta}=(2,0.5,1,0,0\dots,)^\top$. In our simulation, we set the multiplier bootstrap replications $B=500$. 

For comparison, we also use Lasso, SCAD and MCP as initial estimators to obtain plug-in debiased estimators, and construct simultaneous confidence intervals via \cref{alg:inference}. \cref{sim: simulation-inference_lin} and \cref{sim: simulation-inference_log} present the empirical coverage rates (CR) and the average length (AL) of the simultaneous confidence intervals in linear and logistic regression, respectively. For each setting, these results are computed based on $500$ replications.

\begin{sidewaystable}[tbp]
\centering
\begin{minipage}{\textheight}
\resizebox{\linewidth}{!}{\fontsize{8pt}{10pt}\selectfont
\begin{tabularx}{\textwidth}{Xp{3mm}p{3mm}p{3mm}XXXXXXXXXXXXXXXXXX}
\toprule\\[-13pt]
{\scriptsize$\mathbf{\Sigma}$} &{\scriptsize $p$} &{\scriptsize $n$} & &{\fontsize{6pt}{7pt}\selectfont Lasso }&{\fontsize{6pt}{7pt}\selectfont SCAD }&{\fontsize{6pt}{7pt}\selectfont MCP }&{\fontsize{6pt}{7pt}\selectfont HDMA (Lasso) }& {\fontsize{6pt}{7pt}\selectfont HDMA (SCAD) }& {\fontsize{6pt}{7pt}\selectfont HDMA (MCP) }&{\fontsize{6pt}{7pt}\selectfont Lasso }&{\fontsize{6pt}{7pt}\selectfont SCAD }&{\fontsize{6pt}{7pt}\selectfont MCP }&{\fontsize{6pt}{7pt}\selectfont HDMA (Lasso) }& {\fontsize{6pt}{7pt}\selectfont HDMA (SCAD) }& {\fontsize{6pt}{7pt}\selectfont HDMA (MCP) }&{\fontsize{6pt}{7pt}\selectfont Lasso }&{\fontsize{6pt}{7pt}\selectfont SCAD }&{\fontsize{6pt}{7pt}\selectfont MCP }&{\fontsize{6pt}{7pt}\selectfont HDMA (Lasso) }& {\fontsize{6pt}{7pt}\selectfont HDMA (SCAD) }& {\fontsize{6pt}{7pt}\selectfont HDMA (MCP) }\\ \hline\\[-6pt]
 &    &  &    & \multicolumn{6}{c}{$\set{1,2,3,4,5}$}                          & \multicolumn{6}{|c}{$\set{1,2,\dots,\frac{p}{5}}$}                          & \multicolumn{6}{|c}{$\set{1,2,\dots,p}$}                          \\ 
AR(1)     & 100 & 100 & CR & 0.860 & 0.900 & 0.918 & 0.924 & 0.916 & 0.924 & 0.892 & 0.934 & 0.950 & 0.946 & 0.936 & 0.938 & 0.920 & 0.964 & 0.958 & 0.952 & 0.950 & 0.950 \\
      &     &     & AL & 0.310 & 0.324 & 0.322 & 0.312 & 0.312 & 0.315 & 0.381 & 0.402 & 0.399 & 0.388 & 0.388 & 0.391 & 0.456 & 0.482 & 0.478 & 0.463 & 0.464 & 0.468 \\
      &     & 200 & CR & 0.892 & 0.938 & 0.952 & 0.938 & 0.938 & 0.946 & 0.924 & 0.942 & 0.948 & 0.944 & 0.946 & 0.940 & 0.948 & 0.956 & 0.962 & 0.964 & 0.962 & 0.962 \\
      &     &     & AL & 0.225 & 0.228 & 0.228 & 0.228 & 0.228 & 0.228 & 0.273 & 0.276 & 0.277 & 0.276 & 0.276 & 0.277 & 0.321 & 0.325 & 0.326 & 0.326 & 0.326 & 0.326 \\
      & 200 & 100 & CR & 0.810 & 0.916 & 0.944 & 0.928 & 0.930 & 0.936 & 0.914 & 0.946 & 0.956 & 0.966 & 0.962 & 0.964 & 0.912 & 0.960 & 0.964 & 0.956 & 0.962 & 0.958 \\
      &     &     & AL & 0.294 & 0.311 & 0.309 & 0.299 & 0.299 & 0.303 & 0.393 & 0.419 & 0.417 & 0.404 & 0.404 & 0.409 & 0.461 & 0.492 & 0.490 & 0.474 & 0.474 & 0.480 \\
      &     & 200 & CR & 0.872 & 0.922 & 0.940 & 0.932 & 0.924 & 0.930 & 0.948 & 0.964 & 0.964 & 0.972 & 0.972 & 0.972 & 0.946 & 0.968 & 0.966 & 0.978 & 0.972 & 0.974 \\
      &     &     & AL & 0.222 & 0.229 & 0.227 & 0.224 & 0.224 & 0.224 & 0.290 & 0.299 & 0.299 & 0.294 & 0.294 & 0.295 & 0.336 & 0.346 & 0.346 & 0.341 & 0.341 & 0.341 \\
Band    & 100 & 100 & CR & 0.844 & 0.934 & 0.962 & 0.954 & 0.956 & 0.956 & 0.918 & 0.956 & 0.964 & 0.966 & 0.976 & 0.966 & 0.952 & 0.954 & 0.962 & 0.970 & 0.974 & 0.974 \\
      &     &     & AL & 0.340 & 0.354 & 0.352 & 0.344 & 0.343 & 0.345 & 0.426 & 0.447 & 0.443 & 0.433 & 0.434 & 0.436 & 0.515 & 0.542 & 0.539 & 0.526 & 0.526 & 0.531 \\
      &     & 200 & CR & 0.922 & 0.978 & 0.984 & 0.974 & 0.978 & 0.974 & 0.974 & 0.974 & 0.974 & 0.970 & 0.972 & 0.966 & 0.988 & 0.988 & 0.994 & 0.988 & 0.990 & 0.992 \\
      &     &     & AL & 0.270 & 0.271 & 0.273 & 0.271 & 0.272 & 0.272 & 0.339 & 0.342 & 0.342 & 0.341 & 0.342 & 0.341 & 0.407 & 0.411 & 0.412 & 0.411 & 0.411 & 0.411 \\
      & 200 & 100 & CR & 0.742 & 0.894 & 0.932 & 0.938 & 0.942 & 0.946 & 0.880 & 0.952 & 0.970 & 0.976 & 0.972 & 0.970 & 0.924 & 0.962 & 0.964 & 0.982 & 0.980 & 0.978 \\
      &     &     & AL & 0.318 & 0.331 & 0.330 & 0.319 & 0.320 & 0.325 & 0.428 & 0.449 & 0.449 & 0.435 & 0.435 & 0.442 & 0.508 & 0.533 & 0.533 & 0.515 & 0.516 & 0.524 \\
      &     & 200 & CR & 0.810 & 0.946 & 0.954 & 0.950 & 0.954 & 0.946 & 0.946 & 0.946 & 0.950 & 0.958 & 0.962 & 0.954 & 0.968 & 0.958 & 0.964 & 0.970 & 0.974 & 0.974 \\
      &     &     & AL & 0.250 & 0.256 & 0.255 & 0.252 & 0.253 & 0.253 & 0.341 & 0.351 & 0.351 & 0.345 & 0.345 & 0.347 & 0.405 & 0.417 & 0.417 & 0.410 & 0.411 & 0.411 \\ \hline
\end{tabularx}}
\caption{Coverage rates (CR) and average lengths (AL) of simultaneous confidence intervals in linear regression}
\label{sim: simulation-inference_lin}
\end{minipage}
\end{sidewaystable}

\begin{sidewaystable}[tbp]
\centering
\rotatebox{180}{\begin{minipage}{\textheight}
\resizebox{\linewidth}{!}{\fontsize{8pt}{10pt}\selectfont
\begin{tabularx}{\textwidth}{Xp{3mm}p{3mm}p{3mm}XXXXXXXXXXXXXXXXXX}
\toprule\\[-13pt]
{\scriptsize$\mathbf{\Sigma}$} &{\scriptsize $p$} &{\scriptsize $n$} & &{\fontsize{6pt}{7pt}\selectfont Lasso }&{\fontsize{6pt}{7pt}\selectfont SCAD }&{\fontsize{6pt}{7pt}\selectfont MCP }&{\fontsize{6pt}{7pt}\selectfont HDMA (Lasso) }& {\fontsize{6pt}{7pt}\selectfont HDMA (SCAD) }& {\fontsize{6pt}{7pt}\selectfont HDMA (MCP) }&{\fontsize{6pt}{7pt}\selectfont Lasso }&{\fontsize{6pt}{7pt}\selectfont SCAD }&{\fontsize{6pt}{7pt}\selectfont MCP }&{\fontsize{6pt}{7pt}\selectfont HDMA (Lasso) }& {\fontsize{6pt}{7pt}\selectfont HDMA (SCAD) }& {\fontsize{6pt}{7pt}\selectfont HDMA (MCP) }&{\fontsize{6pt}{7pt}\selectfont Lasso }&{\fontsize{6pt}{7pt}\selectfont SCAD }&{\fontsize{6pt}{7pt}\selectfont MCP }&{\fontsize{6pt}{7pt}\selectfont HDMA (Lasso) }& {\fontsize{6pt}{7pt}\selectfont HDMA (SCAD) }& {\fontsize{6pt}{7pt}\selectfont HDMA (MCP) }\\ \hline\\[-6pt]
 &    &  &    & \multicolumn{6}{c}{$\set{1,2,3,4,5}$}                          & \multicolumn{6}{|c}{$\set{1,2,\dots,\frac{p}{5}}$}                          & \multicolumn{6}{|c}{$\set{1,2,\dots,p}$}                          \\ 
AR(1) & 100 & 100 & CR & 0.834 & 0.892 & 0.918 & 0.906 & 0.938 & 0.924 & 0.902 & 0.950 & 0.960 & 0.966 & 0.976 & 0.966 & 0.922 & 0.964 & 0.964 & 0.970 & 0.978 & 0.976 \\
  &     &     & AL & 1.471 & 1.694 & 1.776 & 1.730 & 2.285 & 2.298 & 1.735 & 1.944 & 2.029 & 2.001 & 2.507 & 2.516 & 2.016 & 2.250 & 2.351 & 2.312 & 2.859 & 2.873 \\
  &     & 200 & CR & 0.870 & 0.936 & 0.934 & 0.918 & 0.942 & 0.942 & 0.942 & 0.962 & 0.958 & 0.964 & 0.976 & 0.976 & 0.948 & 0.964 & 0.966 & 0.966 & 0.972 & 0.984 \\
  &     &     & AL & 1.174 & 1.378 & 1.418 & 1.291 & 1.654 & 1.661 & 1.321 & 1.511 & 1.557 & 1.442 & 1.772 & 1.779 & 1.490 & 1.685 & 1.733 & 1.615 & 1.958 & 3.203 \\
  & 200 & 100 & CR & 0.730 & 0.814 & 0.862 & 0.874 & 0.930 & 0.918 & 0.906 & 0.922 & 0.924 & 0.964 & 0.960 & 0.950 & 0.898 & 0.958 & 0.970 & 0.972 & 0.984 & 0.982 \\
  &     &     & AL & 1.321 & 1.583 & 1.624 & 1.600 & 2.216 & 2.233 & 1.683 & 2.317 & 2.014 & 2.011 & 2.573 & 2.596 & 1.955 & 2.588 & 2.342 & 2.339 & 2.980 & 3.006 \\
  &     & 200 & CR & 0.816 & 0.902 & 0.912 & 0.876 & 0.944 & 0.936 & 0.928 & 0.948 & 0.960 & 0.962 & 0.964 & 0.962 & 0.960 & 0.982 & 0.984 & 0.984 & 0.984 & 0.972 \\
  &     &     & AL & 1.101 & 1.302 & 1.350 & 1.240 & 1.637 & 1.643 & 1.315 & 1.500 & 1.556 & 1.456 & 1.812 & 1.814 & 1.489 & 1.693 & 1.758 & 1.646 & 2.030 & 3.014 \\
Band & 100 & 100 & CR & 0.782 & 0.848 & 0.868 & 0.890 & 0.906 & 0.914 & 0.916 & 0.946 & 0.942 & 0.956 & 0.968 & 0.970 & 0.976 & 0.990 & 0.986 & 0.994 & 0.980 & 0.974 \\
  &     &     & AL & 1.647 & 2.003 & 2.102 & 1.927 & 2.491 & 2.526 & 1.925 & 2.282 & 2.385 & 2.215 & 2.710 & 2.741 & 2.281 & 2.716 & 2.846 & 2.609 & 3.172 & 1.964 \\
  &     & 200 & CR & 0.902 & 0.922 & 0.908 & 0.942 & 0.942 & 0.958 & 0.970 & 0.958 & 0.956 & 0.964 & 0.954 & 0.954 & 0.986 & 0.986 & 0.990 & 0.992 & 0.976 & 0.972 \\
  &     &     & AL & 1.467 & 1.817 & 1.854 & 1.596 & 2.021 & 2.042 & 1.616 & 1.933 & 1.972 & 1.741 & 2.116 & 2.123 & 1.865 & 2.214 & 2.261 & 2.011 & 2.407 & 2.418 \\
  & 200 & 100 & CR & 0.644 & 0.756 & 0.816 & 0.812 & 0.892 & 0.886 & 0.876 & 0.908 & 0.930 & 0.948 & 0.950 & 0.940 & 0.916 & 0.960 & 0.964 & 0.978 & 0.978 & 0.976 \\
  &     &     & AL & 1.398 & 1.675 & 1.751 & 1.643 & 2.221 & 2.244 & 1.774 & 2.058 & 2.139 & 2.035 & 2.555 & 2.575 & 2.098 & 2.455 & 2.531 & 2.400 & 2.992 & 2.038 \\
  &     & 200 & CR & 0.840 & 0.908 & 0.906 & 0.918 & 0.942 & 0.946 & 0.962 & 0.970 & 0.964 & 0.970 & 0.944 & 0.948 & 0.968 & 0.982 & 0.978 & 0.980 & 0.968 & 0.968 \\
  &     &     & AL & 1.273 & 1.600 & 1.663 & 1.426 & 1.902 & 1.913 & 1.536 & 1.836 & 1.904 & 1.682 & 2.100 & 2.105 & 1.800 & 2.148 & 2.223 & 1.965 & 2.451 & 2.458 \\ \hline
\end{tabularx}}
\caption{Coverage rates (CR) and average lengths (AL) of simultaneous confidence intervals in logistic regression}
\label{sim: simulation-inference_log}
\end{minipage}}
\end{sidewaystable}
From \cref{sim: simulation-inference_lin}, in linear regression, we can see that the coverage rates of the proposed HDMA methods perform better than the other competing methods in most cases. The coverage rates of Lasso are significantly lower than the pre-specified confidence level, especially when $n=100$ for the index set $\mc{G}=\set{1,2,\dots,5}$. For the average length of the simultaneous confidence intervals, the HDMA methods also yield narrower intervals compared with other methods. From \cref{sim: simulation-inference_log}, in logistic regression, although the average length of our methods is sometimes wider than that of other competing methods, our method provides more reasonable coverage rates, especially when the set of parameters we are interested in is $\mc{G}=\set{1,2\dots,5}$.

\section{Real Data Analysis}\label{sec: real data}
\begin{figure}[htbp]
  \centering
  \begin{subfigure}[b]{0.32\textwidth}
  \caption{HDMA (Lasso)}
    \includegraphics[width=\textwidth]{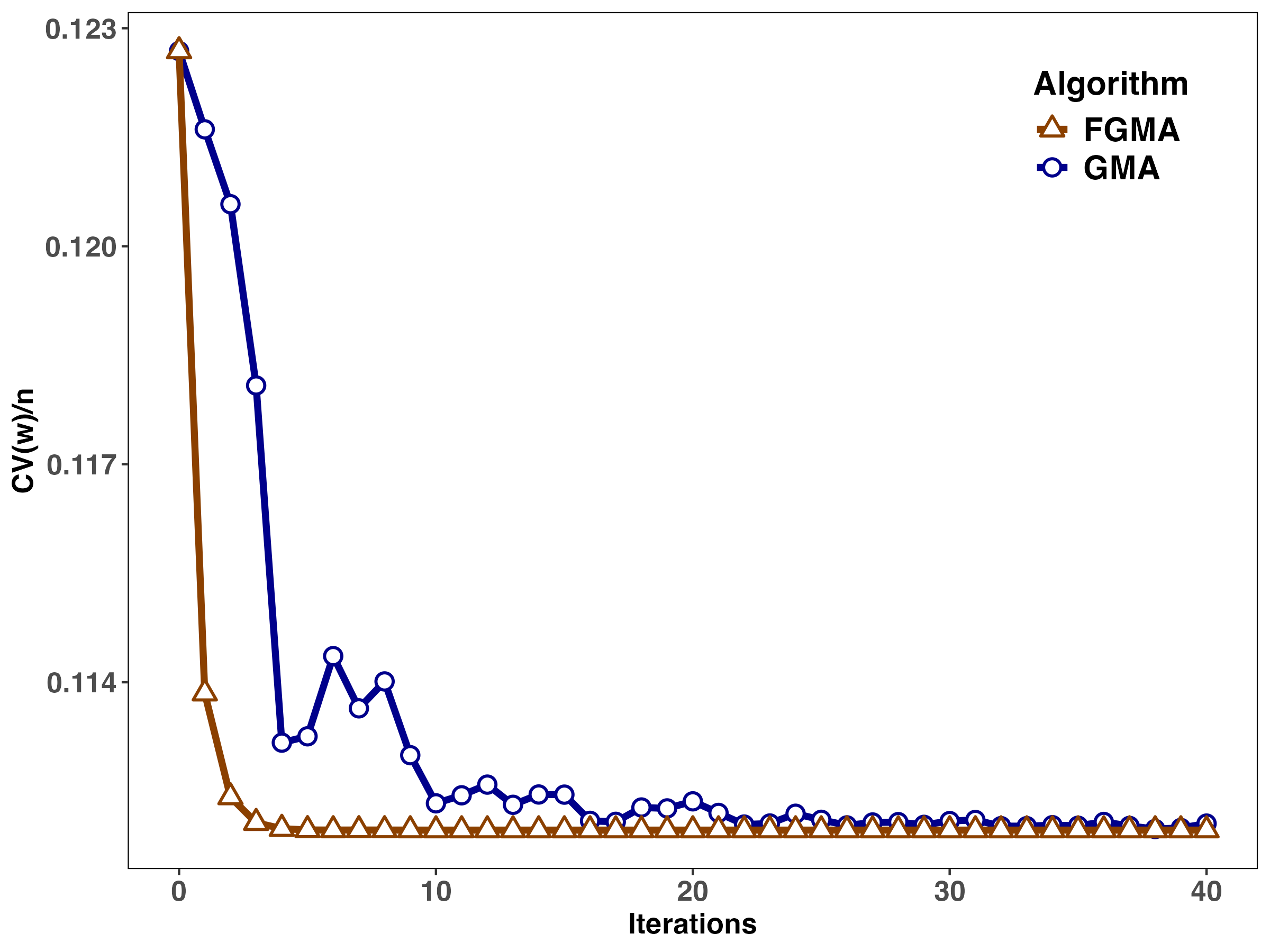}
  \end{subfigure}
  \hfill
  \begin{subfigure}[b]{0.32\textwidth}
  \caption{HDMA (SCAD)}
    \includegraphics[width=\textwidth]{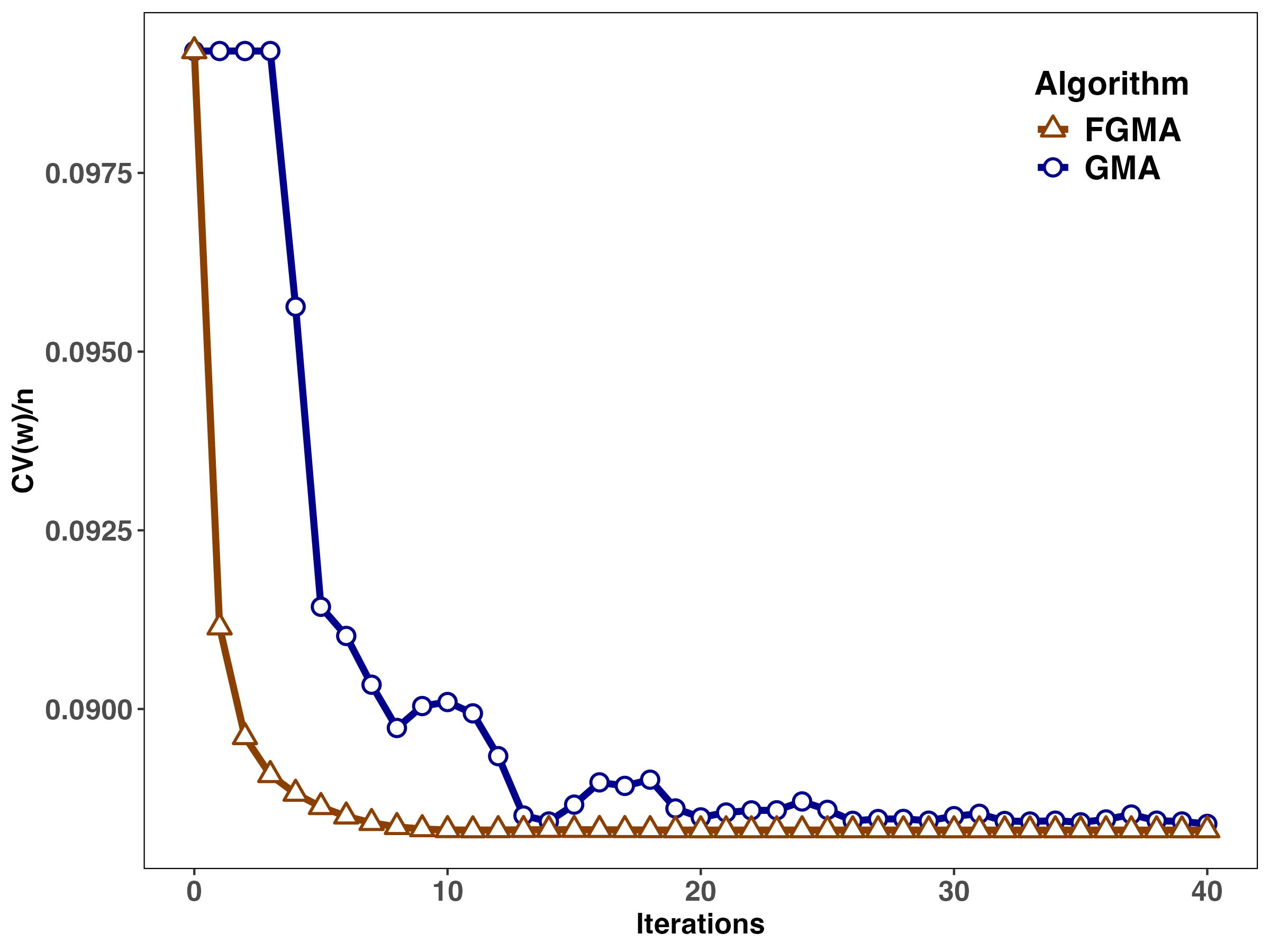}
  \end{subfigure}
  \hfill
  \begin{subfigure}[b]{0.32\textwidth}
  \caption{HDMA (MCP)}
    \includegraphics[width=\textwidth]{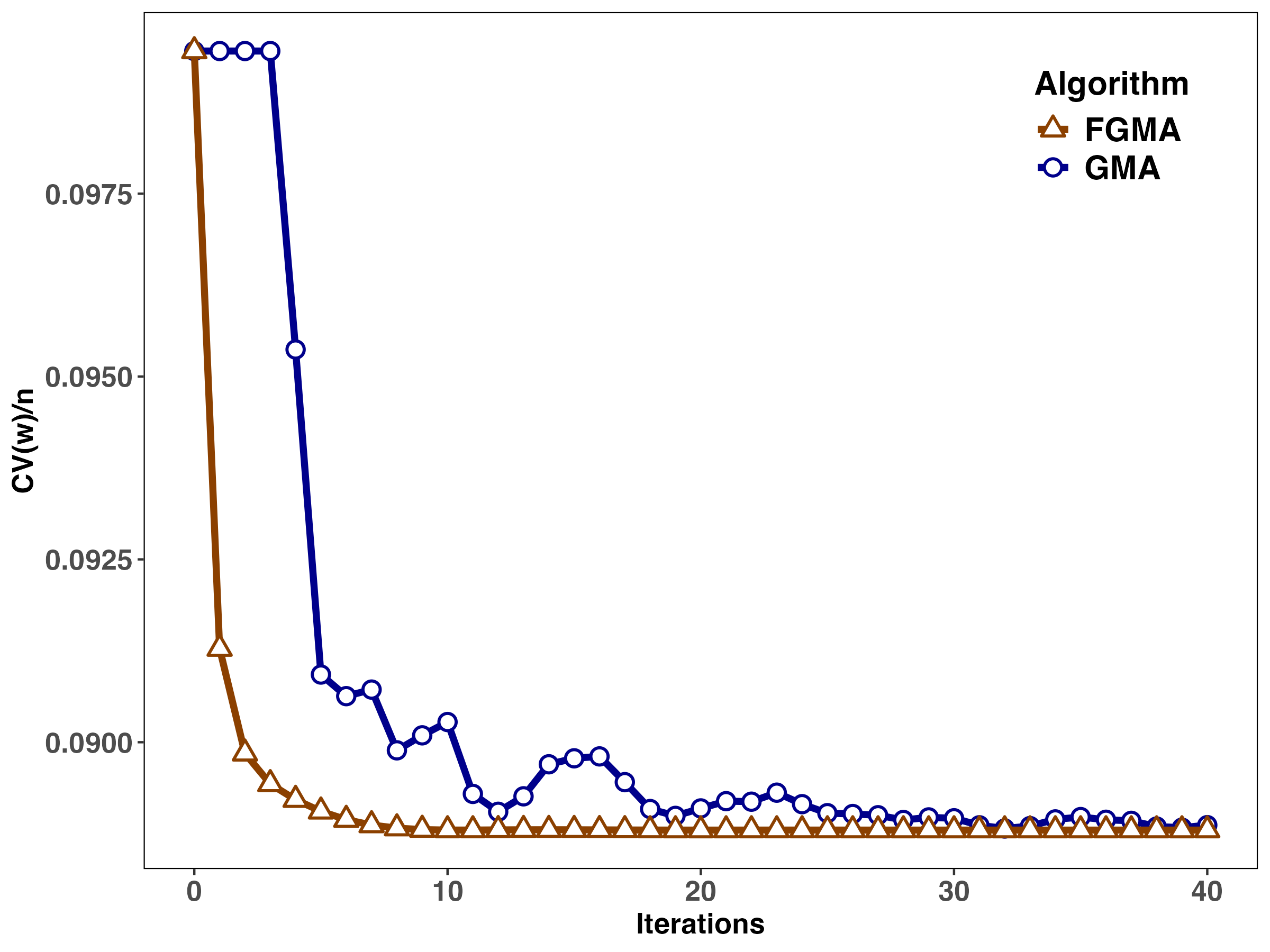}
  \end{subfigure}
  \caption{Cross-validation criterion (divided by $n$) in the first 40 iterations of the proposed FGMA algorithm and the original GMA algorithm for riboflavin data}
  \label{fig:Real-Data-Lin_FGMA}
\end{figure}

In this section, we apply our proposed method to the riboflavin data set concerning riboflavin production rate, which is publicly available in \cite{buhlmann2014high}. This data set contains $n=71$ observations and $p=4,088$ covariates corresponding to 4,088 genes. Following \cite{buhlmann2014high}, we model the riboflavin production rate as a linear model with the covariates. We first apply our proposed HDMA methods to the whole data set. The values of the cross-validation criterion (divided by $n$) in the first 40 iterations of the proposed FGMA \cref{alg:FGMA} and the original GMA \cref{alg:GMA} (in \cref{greedy model averaging}) are presented in \cref{fig:Real-Data-Lin_FGMA}. We can see that our proposed FGMA algorithm has a descent property and has a faster convergence rate than the original GMA algorithm. Note that the initial point in each algorithm corresponds to the model selection by $J$-fold cross-validation. This indicates that model averaging can result in a lower cross-validation criterion than model selection, so it delivers better out-of-sample performance.

\begin{figure}[tbp]
    \centering
    \includegraphics[width=0.7\linewidth]{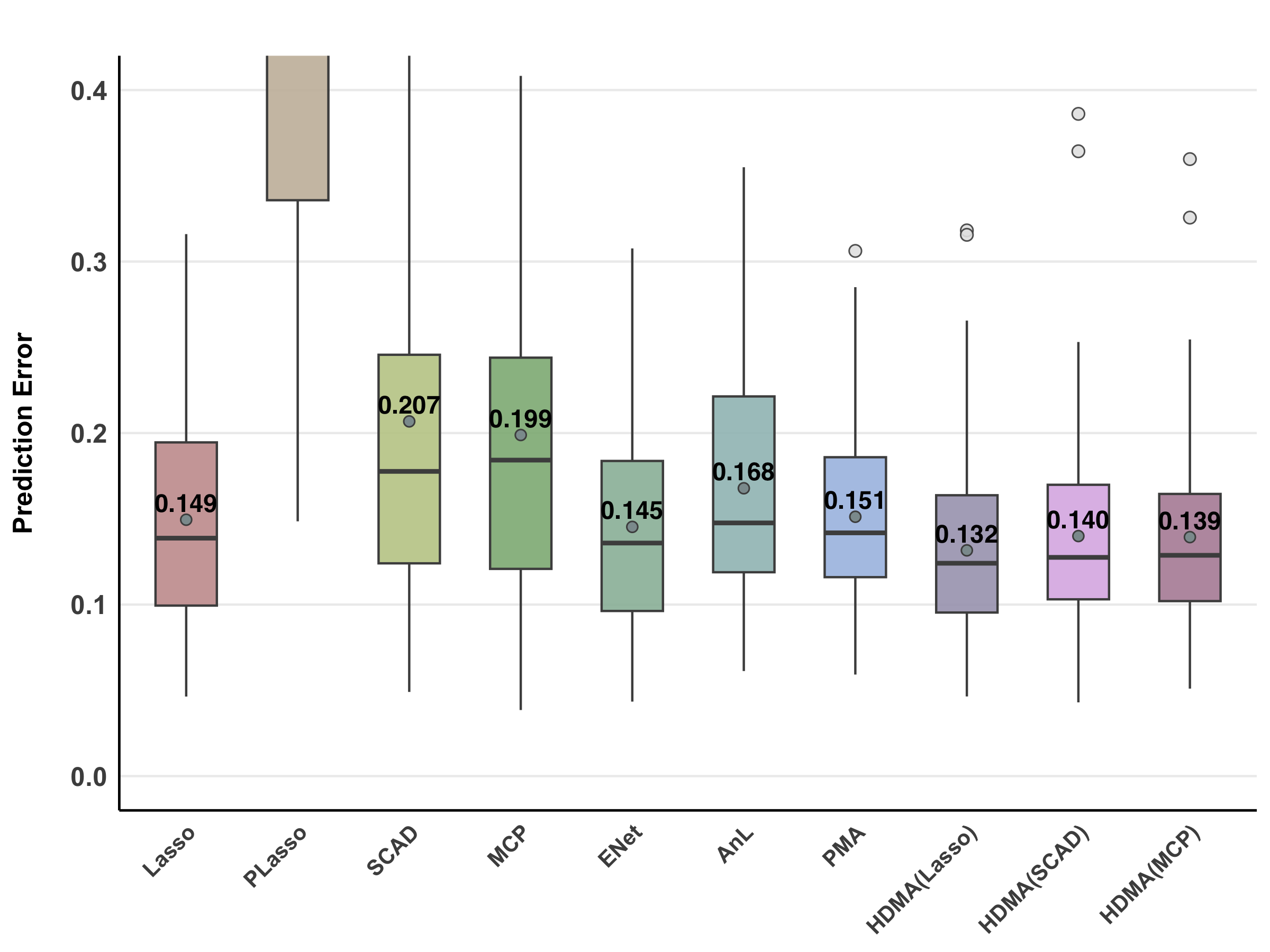}
    \caption{Out-of-sample prediction errors of different methods for riboflavin data (Bold numbers are the means of prediction errors)}
    \label{fig:Real-Data-Lin}
\end{figure}

To investigate the predictive power of the proposed approach, we randomly split the data set into a training set of 50 observations and a test set of 21 observations. We repeat the above data splitting process 100 times. \cref{fig:Real-Data-Lin} displays the prediction errors ($\text{PE}_1$) defined in \eqref{eq:prediction error1} for different approaches. We can observe that the HDMA methods have lower prediction errors than any other method. In addition, the performance of PLasso is poor in the real data. One possible explanation is that when the correlation between the covariates in the active set and other variables is strong, OLS Post-Lasso results in an overfitting model and leads to a terrible prediction. By aggregating a large number of penalized estimators, the HDMA estimators is more robust than other model averaging methods and penalized methods.

\begin{figure}[tbp]
    \centering
    \includegraphics[width=1\linewidth]{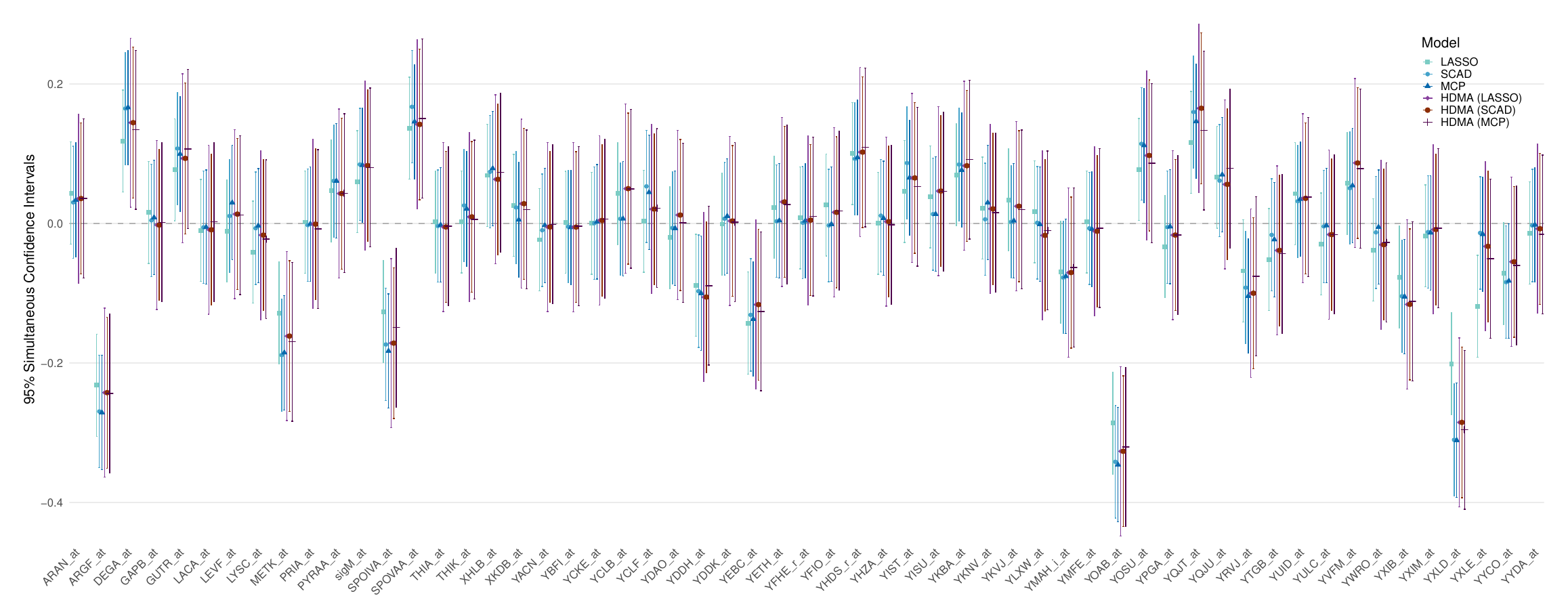}
    \caption{95\% simultaneous confidence intervals of difference methods for riboflavin data}
    \label{fig:Simultaneous confidence intervals}
\end{figure}

To make inference, similar to \cite{javanmard2014confidence}, we make a variable selection using Lasso at first. After the variable selection, we obtain 57 important variables. Then, we implement different methods to construct simultaneous confidence intervals. These intervals are shown in \cref{fig:Simultaneous confidence intervals}. We can see that compared with standard penalized methods, our approach identifies more conservative significant variables. There are nine significant variables which can be identified by all methods simultaneously: "ARGF\_at", "DEGA\_at", "METK\_at", "SPOIVA\_at", "SPOVAA\_at", "YEBC\_at", "YOAB\_at", "YQJT\_at", "YXLD\_at". The HDMA methods with Lasso and MCP penalties exactly find these variables. The variables "YXLD\_at" and "YOAB\_at" are found in \cite{yan2023confidence} and \cite{cai2025statistical}, and the former one is also identified in \cite{meinshausen2009p} and \cite{javanmard2014confidence}. This indicates that our proposed methods can effectively identify significant variables in high-dimensional data.

\section{Concluding Remark}\label{sec: concluding remarks}
In this paper, we present a unified high-dimensional model averaging framework for prediction and constructing post-averaging simultaneous confidence intervals with a large range of loss functions and penalties. The non-asymptotic bounds and asymptotic optimality of the proposed weight estimator are derived in the scenario where all candidate models are misspecified. When there is at least one correct model, we show that our method assigns all weights on correct models asymptotically and establish the Gaussian and bootstrap approximation results. Based on these results, we can construct post-averaging simultaneous confidence intervals by Gaussian multiplier bootstrap technique. To address the simplex-constrained optimization problem, we develop a fast greedy model averaging algorithm, which has a faster convergence rate and has a more stable performance in practice than the original greedy model averaging algorithm introduced in \cite{dai2012deviation} and \cite{he2023rank}. The convergence analysis of the FGMA algorithm matches the optimal rate of the first-order algorithms.

The basic $J$-fold cross-validation is used in our method. In recent years, many variants of cross-validation method have been proposed. For example, \cite{rad2020scalable} proposed a computationally efficient approximate leave-one-out (ALO) method to reduce the bias coming from K-fold cross-validation, which performs as well as leave-one-out cross-validation. \cite{bates2024cross} proposed a nested cross-validation (NCV) method to estimate the prediction error and empirically showed that NCV has a more accurate predictive power than the traditional cross-validation method. How to extend these new cross-validation methods to model averaging remains for future research.

% Acknowledgements and Disclosure of Funding should go at the end, before appendices and references

%\acks{This work was partially supported by the National Key R \& D Program of China (2021YFA10
%
%\noindent 00100, 2021YFA1000101, 2021YFA1000300), and the National Natural Science Foundation of China (72331005).}


\newpage

\appendix
\counterwithin{equation}{section}
\begin{center}
    \LARGE\bfseries Appendix
\end{center}

\cref{appendix:discussions} contains some further discussions on the results of weight consistency and asymptotic optimality. \cref{greedy model averaging} presents a greedy model averaging algorithm. Some useful lemmas are given in \cref{appendix:lemma}. The proof of propositions, the proof of theorems, the proof of lemmas are collected in \cref{appendix:proof-of-proposition,appendix:proof-of-theorem,appendix:proof-of-lemma}, respectively.

\section{Some Discussions}\label{appendix:discussions}
\subsection{Discussion on the Support of $w^*$}\label{remark:support of w^*}
Recall that $w^*$ is a minimizer of $L(Y_1,\<B^*w,X_1\>)$ over the simplex $\mc{W}^K$. By Karush-Kuhn-Tucher (KKT) condition, $w^*$ satisfies
\begin{align*}
    \E B^{*\top} X_1\frac{\partial L(Y_1,\<B^*w^*,X_1\>)}{\partial \mu}-u+v\mathbbm{1}_K=0,
\end{align*}
where $\mathbbm{1}_K$ is a $K$-dimensional vector of 1's and  $v\in\R, u\in\R^K, u_k\ge 0, u_kw^*_k=0$ for $k\in\brac{K}$. Then the support of $w^*$ is contained in an equicorrelation set: 
\begin{align}
    \text{supp}(w^*)\subseteq\setBig{k\in\brac{K}\Big|\E\beta^{*\top}_{(k)} X_1\bracBig{-\frac{\partial L(Y_1,\<B^*w^*,X_1\>)}{\partial \mu}}=v}\overset{\triangle}{=}\mc{J},\label{eq:support-of-w}
\end{align}
and for any $k\notin \mc{J}$, $\E\beta^{*\top}_{(k)} X_1\bracBig{-\frac{\partial L(Y_1,\<B^*w^*,X_1\>)}{\partial \mu}}<v$. 

For the $k$-th model, the quantity $\E\beta^{*\top}_{(k)} X_1\bracBig{-\frac{\partial L(Y_1,\<B^*w^*,X_1\>)}{\partial \mu}}$ can be interpreted as the regression coefficient between the prediction value $\beta^{*\top}_{(k)} X_1$ and the pseudo-residual (negative gradient) $-\frac{\partial L(Y_1,\<B^*w^*,X_1\>)}{\partial \mu}$  at the population level. Then positive weight values are assigned only to the candidate models with the strongest positive correlations between the prediction value and the pseudo-residual. 

Let $v$ be the quantity in \eqref{eq:support-of-w}. Then $\abs{v}$ can be bounded by $C(1-\linf{w^*})$, where $C$ is a positive constant independent of $n$. To see this, we define $h(w)=\E L(Y,\<B^*w,X_1\>)$ and let $e^{(k^*)}$ be the $k^*$-th coordinate unit vector in $\R^K$, where $k^*=\argmax_k{w^*_k}$. By \eqref{eq:support-of-w} and $\<e^{(k^*)},\nabla h(e^{(k^*)})\>=0$, we have
\begin{align*}
    \abs{v}=&~\abs{\<e^{(k^*)},\nabla h(w^*)\>}\\
    =&~\abs{\<e^{(k^*)},\nabla h(e^{(k^*)})-\nabla h(w^*)\>}\\
    =&~\abs{\<e^{(k^*)},\int_0^1\nabla^2 h(w^*+t(e^{(k^*)}-w^*))(e^{(k^*)}-w^*)\text{d} t\>}\\
    \leq&~L_ue^{(k^*)\top}B^{*\top}(\E XX^\top)B^*(e^{(k^*)}-w^*)\\
    \leq&~C(1-\linf{w^*}),
\end{align*}
where we use $\ltwo{\beta^*_{(k)}}\leq K_\beta$ and \cref{lasso:lossfunc,lasso:covariate}. This indicates that $\abs{v}$ measures the difference between model selection and model averaging. Moreover, if the weighted pseudo-true value $B^*w^*$ equals to one of $\beta^*_{(k)}$'s, then $v=0$.

%%%%%%%%%%%%%%%%%%%%%%%%%%%%%%%%%%%%%%%%%%%%
\subsection{Some Examples for the Convergence Rate of $\hat{w}$}\label{remark:rate of w}
To better understand \cref{MA}, we consider a linear model. Suppose that $Y_1=\sum_{j=1}^\infty\beta_j \mathrm{x}_{1j}+\varepsilon_1$, $\E(\varepsilon_1|\mathbf{x}_1)=0$, $\Var\paren{\varepsilon_1}=\sigma^2$ and $L(y,\mu)=\frac{1}{2}(y-\mu)^2$. In addition, we assume that $\Psi_{n,q}=T_{n,q}.$
\begin{enumerate}
    \item \textit{The largest candidate model in the nested part can minimize the prediction risk}, that is, $\beta^*_{(K_{ne})}=\argmin_{\beta\in\R^p}\E L(Y_1,\<\beta,X_1\>)$. We assume that $\beta^*_{(K_{ne})}\neq\beta^*_{(k)}$ for any $k\neq K_{ne}$. By choosing $S_*=\set{K_{ne}}$, we have $\lone{w^*_{S_*^c}}=0$, $\abs{S_*}=1$ and $\phi^2(S_*)\ge\lambdamin{B^{*\top}_{ne} B^*_{ne}}$, where $B^{*}_{ne}=(\beta^*_{(1)},\dots,\beta^*_{(K_{ne})})$. Consequently, with high probability, $\lone{\hat{w}-w^*}\lesssim \lambdamin{B^{*\top}_{ne} B^*_{ne}}^{-1}T_{n,q}$. When the minimum singular value of the matrix $B^*_{ne}$ is bounded away from zero, $\lone{\hat{w}-w^*}=O_p(T_{n,q}).$
    \item \textit{The covariates are uncorrelated}. That is, $\E \mathrm{x}_{1,j}=0$ and $\E \mathrm{x}_{1,j_1}\mathrm{x}_{1,j_2}=\delta_{j_1,j_2}$ for $j, j_1, j_2\in\N$, where $\delta_{j_1,j_2}$ is the Kronecker delta. We assume that $\set{\abs{\beta_j}}_{j=1}^\infty$ is a nonincreasing positive sequence, which is the same as Assumption 2 in \cite{peng2024optimality} and Assumption 1 in \cite{peng2022improvability}. For simplicity, suppose that $\mc{A}_{K_{ne}}=\brac{p_0}$ and $\ltwo{\beta^*_{(K_{ne})}}\ge\cdots\ge\ltwo{\beta^*_{(K)}}$. In this case, we must have $w^*_k=0$ for any $k\in\set{1,\dots,K_{ne}-1}$ and
    \begin{align}
        \E L(Y_1,\<B^*w^*,X_1\>)=\frac{1}{2}\sum_{k=K_{ne}}^{K}(1-w_k^*)^2\ltwo{\beta^*_{(k)}}^2+\frac{1}{2}\sum_{j=p+1}^\infty \abs{\beta_j}^2+\frac{1}{2}\sigma^2.\label{eq: uncorrelated loss}
    \end{align} Note that 
\begin{align}
    \E\beta^{*\top}_{(k)} X_1\bracBig{-\frac{\partial L(Y_1,\<B^*w^*,X_1\>)}{\partial \mu}}=\begin{cases}
    (1-w^*_{K_{ne}})\ltwo{\beta^*_{(k)}}^2,&~ \text{if }k=1,\dots,K_{ne}-1\\
        (1-w^*_k)\ltwo{\beta^*_{(k)}}^2,&~ \text{if }k=K_{ne},\dots,K
    \end{cases},\label{eq:equicorrelation}
\end{align}
which is nonnegative for each $k\in\brac{K}$ and can be maximized when $k=K_{ne}$ by \eqref{eq:support-of-w}. Thus, $w^*_{K_{ne}}\ge w^*_{K_{ne}+1}\ge\cdots\ge w^*_K$ and $\abs{\text{supp}(w^*)}\ge 2$. The following proposition shows that when $\ltwo{\beta^*_{(k)}}$'s are dominated by $\ltwo{\beta^*_{(K_{ne})}}$, $w^*$ is sparse.
\begin{proposition}\label{prop:sparse}
    Suppose that the covariates are uncorrelated with unit variance and $\set{\abs{\beta_j}}_{j=1}^\infty$ is a nonincreasing positive sequence. If $\frac{\ltwo{\beta^*_{(K_{ne}+2)}}^2}{\ltwo{\beta^*_{(K_{ne})}}^2}+\frac{\ltwo{\beta^*_{(K_{ne}+2)}}^2}{\ltwo{\beta^*_{(K_{ne}+1)}}^2}<1$ and $\ltwo{\beta^*_{(K_{ne})}}\ge\cdots\ge\ltwo{\beta^*_{(K)}}$, then $\abs{\text{supp}(w^*)}=2.$
\end{proposition}

To control the term $\abs{S_*}T_{n,q}/\phi^2(S_*)$ in \eqref{eq:oracle-inequality}, the size of the largest nested model could not diverge too fast. Let $\rho_n=\linf{\E\parenBig{\frac{\partial L(Y_1,\<B^*w^*,X_1\>)}{\partial \mu}B^{*\top}X_1}}=(1-w^*_{K_{ne}})\ltwo{\beta^*_{(K_{ne})}}^2$,  $p_0=\max\set{j\in\brac{p}:\abs{\beta_j}^2>T_{n,q}^\iota}$ and $d_1=\lfloor \nu p_0\rfloor$, where $\iota,\nu\in(0,1]$ are some constants. Suppose that $\ltwo{\beta^*_{(K_{ne}+1)}}^2\asymp T_{n,q}^\iota$ and $\text{supp}(w^*)=\set{K_{ne},K_{ne}+1}$. Then we can obtain that $ \rho_n=(1-w^*_{(K_{ne}+1)})\ltwo{\beta^*_{(K_{ne}+1)}}^2\leq\ltwo{\beta^*_{(K_{ne}+1)}}^2\asymp T_{n,q}^\iota$ and
\begin{align*}
    \phi^2(\set{K_{ne}})\ge&~ \ltwo{\beta^*_{(K_{ne})}}^2-\ltwo{\beta^*_{(K_{ne}-1)}}^2\ge d_1 T_{n,q}^\iota,\\
    \lone{w^*_{\set{K_{ne}}^c}}=&~1-w^*_{K_{ne}}=\rho_n/\ltwo{\beta^*_{(K_{ne})}}^2\lesssim 1/p_0,\\
    \phi^2(\set{K_{ne},K_{ne}+1})\ge&~\ltwo{\beta^*_{(K_{ne}+1)}}^2\asymp T_{n,q}^\iota,\\
  \lone{w^*_{\set{K_{ne},K_{ne}+1}^c}}=&~0. 
\end{align*}
Consequently, $\lone{\hat{w}-w^*}=O_p(\frac{1}{p_0}\wedge T_{n,q}^{1-\iota})$. Consider the following two different settings on the coefficients.
\begin{itemize}
    \item {\it Polynomially decaying coefficients.} When $\abs{\beta_j}=c_1j^{-\alpha_1}, \alpha_1>1/2$, we have $p_0\asymp T_{n,q}^{-\iota/(2\alpha_1)}$ and $\lone{\hat{w}-w^*}=O_p(T_{n,q}^{\iota/(2\alpha_1)}\wedge T_{n,q}^{1-\iota})$. If $\frac{1}{2}<\alpha_1\leq 1$, then $s_u= p_0$, $p_0\asymp(\frac{\log p}{n})^{-\iota/(4\alpha_1+\iota)}$ and $\lone{\hat{w}-w^*}=O_p((\frac{\log p}{n})^{\iota/(4\alpha_1+\iota)}\wedge(\frac{\log p}{n})^{2\alpha_1(1-\iota)/(4\alpha_1+\iota)})$ when we select $q=0$ in \cref{MA}. If $1<\alpha_1<\infty$, then $h_{q,u}^q<\sum_{j=1}^\infty\abs{\beta_j}^q<\infty$, $p_0\asymp (\frac{\log p}{n})^{-\iota(2\alpha_1-1-\tau)/(8\alpha_1^2)}$ and $\lone{\hat{w}-w^*}=O_p((\frac{\log p}{n})^{\iota(2\alpha_1-1-\tau)/(8\alpha_1^2)}\wedge (\frac{\log p}{n})^{(1-\iota)(2\alpha_1-1-\tau)/(4\alpha_1)})$ when we select $q=\frac{1+\tau}{\alpha_1}$ with some constant $\tau\in(0,\alpha_1-1]$. 
    \item {\it Exponentially decaying coefficients.} When $\abs{\beta_j}=c_1\exp(-c_2j^{\alpha_2}), \alpha_2>0$, we have $p_0\asymp\brac{\log(T_{n,q}^{-\iota/(2c_2)})}^{1/\alpha_2}$ and $\lone{\hat{w}-w^*}=O_p(\brac{\log(T_{n,q}^{-\iota/(2c_2)})}^{-1/\alpha_2}\wedge T_{n,1}^{1-\iota})$. Choosing $q\in(0,1]$, we have $h_{q,u}^q<\infty$ and $\lone{\hat{w}-w^*}=O_p((\log n-\log\log p)^{-1/\alpha_2}\wedge(\frac{\log p}{n})^{(1-\iota)(2-q)/4}).$
\end{itemize}

\end{enumerate}
%%%%%%%%%%%%%%%%%%%%%%%%%%%%%%%%
\subsection{Discussion on the Uniqueness of $w^*$}
\label{remark:uniqueness}
Clearly, if the loss $L(y,\mu)$ is strictly convex in $\mu$, $\E X_1X_1^\top$ is positive definite and the matrix $B^*$ is of full column rank, then $w^*$ is unique. However, it is possible that $w^*$ is not unique in some cases. Fortunately, \cref{MA} do not need the uniqueness of $w^*$ and it holds for any $w^*\in\argmin_{w\in\mc{W}^K}\E L(Y_1,\<B^*w,X_1\>)$. For example, consider the linear model with uncorrelated covariates in \cref{remark:rate of w} and assume $\ltwo{\beta^*_{(K_{ne}+1)}}=\ltwo{\beta^*_{(K_{ne}+2)}}$. By \eqref{eq: uncorrelated loss}, for any positive constants $a, b$ that satisfy $a+b=w^*_{K_{ne}+1}+w^*_{K_{ne}+2}$ we have 
\begin{align*}
    w^*(a,b)\overset{\triangle}{=}(w^*_1,\dots,w^*_{K_{ne}},a,b,w^*_{K_{ne}+3},\dots,w^*_K)\in\argmin_{w\in\mc{W}^K}\E L(Y_i,\<B^*w,X_i\>).
\end{align*} In this case, weights consistency holds for any $w^*(a,b)\in\argmin_{w\in\mc{W}^K}\E L(Y_1,\<B^*w,X_1\>)$. 

Obviously, if there are two candidate models which can minimize the risk function $\E L(Y_1,\<\beta,X_1\>)$ over $\R^p$, the weights consistency fails. For example, consider the linear model $Y_1=\sum_{j=1}^\infty\beta_j\mathrm{x}_{1,j}+\eps_1$, where $\mathrm{x}_{1,p+1}=f^0(\mathrm{x}_1)$ for some function $f^0(\cdot)$, $\mathrm{x}_{1,j}=0$ for $j\neq p+1$. By Proposition 3 in \cite{buhlmann2015high}, if $(\mathrm{x}_{1,1},\dots,\mathrm{x}_{1,p})^\top$ has a joint Gaussian distribution, then all pseudo-true values of the candidate models in the nested part are sparse and can minimize the risk $\E (Y_1-\<\beta,X_1\>)^2$ over $\R^p$, since $(c,0,\dots,0)^\top\in\argmin_{\beta\in\R^p}\E(Y_1-\<\beta,X_1\>)^2$ for some $c\in\R$. Then $e^{(1)},e^{(2)}\in\argmin_{w\in\mc{W}^K}\E (Y_1-\<B^*w,X_1\>)^2$ if $\set{1}\subseteq\mc{A}_1\subseteq\mc{A}_2$, where $e^{(k)}$ is the $k$-th coordinate unit vector in $\R^K$. If we take $S_*=\set{1}$, then $\phi(S_*)=0$ and the right-hand side of the oracle inequality \eqref{eq:oracle-inequality} is infinite. In fact, $\hat{w}$ cannot be simultaneously close to $e^{(1)}$ and $e^{(2)}$. To handle this case, we can assume that there are $K_0-1$ candidate models that cannot minimize the risk for some positive integer $K_0\leq K$. Using the same arguments as the proof of \cref{weights2one}, the sum of weights assigned to candidate models that can minimize the risk $\E L(Y_1,\langle \beta,X_1\rangle)$ converges to one with high probability. 
%%%%%%%%%%%%%%%%%%%%%%%%%%%%%%%%%%%%%%%%%%%%%%%%%%%%%%%
\subsection{A Specific Example for the Convergence Rate of $\xi_n$ in \cref{HD-MA:xi}}\label{remark:xi}

In this example, we consider the linear model with uncorrelated covariates in \cref{remark:rate of w}. By \eqref{eq: uncorrelated loss} and \eqref{eq:equicorrelation}, we can obtain that $\abs{\text{supp}(w^*)}\ge2$ and $\xi_n$ can be bounded by $\rho_n/2$ from below:
\begin{equation}
    \begin{aligned}
    \xi_n=&~\E L(Y_1,\<B^*w^*,X_1\>)-\E L(Y_1,\sum_{j=1}^\infty\beta_j\mathrm{x}_{1j})\\
    \ge&~\frac{1}{2}(1-w^*_{(K_{ne})})^2\ltwo{\beta_{(K_{ne}+1)}}^2+\frac{1}{2}(1-w^*_{(K_{ne}+1)})^2\ltwo{\beta_{(K_{ne}+1)}}^2+\frac{1}{2}\sum_{j=p+1}^\infty\abs{\beta_j}^2\\
    \ge&~ \frac{1}{2}\rho_n+\frac{1}{2}\sum_{j=p+1}^\infty\abs{\beta_j}^2\\
    \ge&~ \frac{1}{2}\rho_n,
\end{aligned}\label{eq:xi bounded below}
\end{equation}
where $\rho_n=\linf{\E\parenBig{\frac{\partial L(Y_1,\<B^*w^*,X_1\>)}{\partial \mu}B^{*\top}X_1}}$.

If $\rho_n$ is bounded below (for example, the class of candidate models is fixed), then $\hat{w}$ is asymptotically optimal provided that $T_{n,q}\to0$. Since $\rho_n=(1-w^*_{K_{ne}})\ltwo{\beta^*_{(K_{ne})}}^2=(1-w^*_{K_{ne}+1})\ltwo{\beta^*_{(K_{ne}+1)}}^2$, we have $1-w^*_{K_{ne}}\leq\ltwo{\beta^*_{(K_{ne}+1)}}^2/\ltwo{\beta^*_{(K_{ne})}}^2$. Suppose that $p_0=\max\set{j\in\brac{p}:\abs{\beta_j}^2>T_{n,q}^\iota}-1$, where $\iota\in(0,1)$. Then, we have $w^*_{K_{ne}}\to1$ as $p_0\to\infty$ and $2\xi_n\ge\rho_n\ge\frac{1}{2}\ltwo{\beta^*_{(K_{ne}+1)}}^2\gtrsim T_{n,q}^\iota$ by \eqref{eq:xi bounded below}. Then $\hat{w}$ is asymptotically optimal provided that $T_{n,q}^{1-\iota}\to0$. Thus, we have the following results for this special case.
\begin{itemize}
    \item When $p_0$ diverges at a slow rate ($\iota\in[0,1)$), $\hat{w}$ enjoys both consistency and asymptotic optimality.
    \item When $p_0$ diverges at a faster rate ($\iota\ge1$ and $T_{n,q}/\phi^2(\set{K_{ne}})\to0$), $\hat{w}$ only enjoys consistency.
\end{itemize}

%%%%%%%%%%%%%%%%%%%%%%%%%%%%%%%%%%%%%%%%%%%%%%%%%%%%%%%%%%%%
\section{Greedy Model Averaging Algorithm}\label{greedy model averaging}

To minimize the convex objective function $CV(w)$ over a simplex $\mc{W}^K$, we can use the following greedy model averaging algorithm described in \cref{alg:GMA}, which is similar to the GMA-0 algorithm of \cite{dai2012deviation}. This algorithm leads to a $N$-sparse approximate weights estimator after the $N$-th iteration and only uses zero-order information in each iteration. We show the convergence result in \cref{GMA}, which is similar to Theorem 4 in \cite{he2020functional} and only requires that the convex objective function $CV(w)/n$ has a Hessian matrix with bounded eigenvalues.

As \cref{alg:GMA} only requires zero-order information in each step, in practice, we can apply this algorithm (although the condition in \cref{GMA} cannot be satisfied) when the loss function is not smooth (e.g. quantile check loss function). To avoid using first-order information, we suggest the algorithm stopping criterion as $\alpha_N<\eps_1$ and $\linf{\hat{w}^{(N)}-\hat{w}^{(N-1)}}<\eps_2$ for some positive constants $\eps_1$ (say, $10^{-2}$) and $\eps_2$ (say, $10^{-3}$). The former condition ensures that the objective function can decrease, as the value of the objective function could remain unchanged if $\alpha_N$ is not sufficiently small.
\begin{algorithm}
\caption{Greedy model averaging}
\label{alg:GMA}
\begin{algorithmic}[1]
\INPUT{Initial value $\hat{w}^{(0)}=\hat{w}^{(0)}=e^{(\varsigma)}$ (where $\varsigma=\argmin_{k\in\brac{K}}CV(e^{(k)})$), constants $\eps_1>0$ and $\eps_2>0$.}
\FORALL{$N=1,2,\dots$ until $\|\hat{w}^{(l)} - \hat{w}^{(l+1)}\|_\infty\leq\eps_1$ and $\alpha_N<\eps_2$}
\STATE Set $\alpha_N=\frac{2}{N+2}$.
\STATE Compute $\varsigma=\argmin_{k\in\brac{K}}CV(\hat{w}^{(N-1)}+\alpha_N(e^{(k)}-\hat{w}^{(N-1)}))$, where $e^{(k)}$ is the $k$-th coordinate unit vector in $\R^K$.
\STATE Set $\hat{w}^{(N)}=\hat{w}^{(N-1)}+\alpha_N(e^{(\varsigma)}-\hat{w}^{(N-1)})$.
\ENDFOR
\OUTPUT{Approximate weight estimator $\hat{w}^{(N)}$.}
\end{algorithmic}
\end{algorithm}

\begin{theorem}\label{GMA}
    Under \cref{lasso:lossfunc}, for any $N\in\N$, we have
    \begin{align*}
        \frac{CV(\hat{w}^{(N+1)})}{n}\leq \frac{CV(\hat{w})}{n}+\frac{4L^\prime_{CV}}{N+4},
    \end{align*}
    where $L^\prime_{CV}=\sup_{w\in\mc{W}^K}\matrixnorm{\nabla^2 CV(w)/n}_2$.
\end{theorem}

\section{Proof of Proposition}\label{appendix:proof-of-proposition}

%%%%%%%%%%%%%%%%%%%%
\vspace{0.2cm}
\subsection{Proof of \cref{lasso}}
\begin{proof} It suffices to show the desired result for the case $k=1$. Without loss of generality, we assume that $\set{X_i}_{i\in\brac{n}}$, $\beta^*_{(1)}$ and $\hat{\beta}_{(1)}$ select only the elements from the set $\mc{A}_1 \subseteq \brac{p}$. Denote $L_n(Y,\langle \beta,X\rangle)=\dn\sumn L(Y_i,\langle \beta,X_i\rangle)$ and
\begin{align*}
    \delta\mc{L}\paren{\Delta,\beta}\overset{\triangle}{=}&~L_n(Y,\langle \beta+\Delta,X\rangle)-L_n(Y,\langle\beta,X\rangle)-\langle \grad_\beta L_n(Y,\langle \beta,X\rangle),\Delta\rangle,
\end{align*}
where $\beta\in\R^{p_1}.$ Applying \cref{lemma:RSC} with $\beta_0=\beta^*_{(1)}$ and $R_1=1,R_2=0$, if $n\ge c_1$, the following event $\mc{B}_1$ holds with probability at least $1-2\exp\paren{-c_2n}$:
\begin{align*}
    \mc{B}_1\overset{\triangle}{=}\setBig{\delta\mc{L}\paren{\Delta,\beta^*}\ge C_1\ltwo{\Delta}\parenBig{\ltwo{\Delta}-C_2\sqrt{\frac{\log p_1}{n}}\lone{\Delta}}\text{ for any }\Delta \in \setBig{\Delta\in\R^{p_1}\Big|\ltwo{\Delta}\leq1}},
\end{align*}
 where $c_1,c_2,C_1,C_2$ are some positive constants not depending on $n$. 

Note that
\begin{align*}
    \grad_\beta L_n(Y,\langle\beta^*_{(1)},X\rangle)=\dn\sumn\frac{\partial L(Y_i,\langle\beta^*_{(1)},X_i\>)}{\partial \mu}X_i
\end{align*}
and for each $j\in \brac{p_1}$
\begin{align*}
    \norm{\frac{\partial L(Y_i,\langle\beta^*_{(1)},X_i\>)}{\partial \mu}X_{ij}}_{\psi_1}
   \leq \norm{\frac{\partial L(Y_i,0)}{\partial \mu}X_{ij}}_{\psi_1}+L_u\norm{\<\beta^*_{(1)},X_i\>X_{ij}}_{\psi_1}\leq K_yK_x+L_uK_\beta K_x^2,
\end{align*}
where we use Lemma 2.7.7 of \cite{vershynin2018high}.

Since $\E\grad_\beta L_n(Y,\<\beta^*_{(1)},X\>)=0$, by the Bernstein's inequality, for any $j\in \brac{p_1}$ we can obtain
\begin{align*}
    \P\parenBig{\absBig{\dn\sumn\frac{\partial L(Y_i,\<\beta^*_{(1)},X_i\>)}{\partial \mu}X_{ij}}>Ct}\leq 2\exp\paren{-cnt^2},
\end{align*}
where $t\in[0,1)$. Let $t=\sqrt{\frac{1}{c}}\parenBig{\delta+\sqrt{\frac{\log p_1}{n}}}$. By a union bound, we have
\begin{align}
    \P\parenBig{\linf{\grad_\beta L_n(Y,\<\beta^*_{(1)},X\>)}>C\paren{\delta+\sqrt{\frac{\log p_1}{n}}}}\leq 2\exp\paren{-n\delta^2}.\label{eq:infnormbound}
\end{align}
Set $\lambda_{n,1}=C\paren{\delta+\sqrt{\frac{\log p_1}{n}}}$. Then the event 
\begin{align*}
\mc{B}_2=\setBig{\lambda_{n,1}\ge 2\linf{\grad_\beta L_n(Y,\<\beta^*_{(1)},X\>)}} 
\end{align*}
holds with probability at least $1-2\exp\paren{-n\delta^2}$.  

We first show the desired results for the Lasso penalty. Let $\tau$ be a positive constant that we will choose later. Under the events $\mc{B}_1,\mc{B}_2$, for any $S\subseteq \brac{p_1}$ and any $\Delta \in \setBig{\Delta\in\R^{p_1}\Big|\ltwo{\Delta}\leq \tau}$, by triangle inequality, we have
\begin{equation}\label{proof:prop1-1}
    \begin{aligned}
    \mc{F}\paren{\Delta}\overset{\triangle}{=}&~L_n\paren{Y,\<\beta^*_{(1)}+\Delta,X\>}-L_n\paren{Y,\<\beta^*_{(1)},X\>}+\lambda_{n,1}\parenBig{\lone{\beta^*_{(1)}+\Delta}-\lone{\beta^*_{(1)}}}\\
    \ge&~\<\Delta,\grad_\beta L_n\paren{Y,\<\beta^*_{(1)},X\>}\>+C_1\ltwo{\Delta}\parenBig{\ltwo{\Delta}-C_2\sqrt{\frac{\log p_1}{n}}\lone{\Delta}}\\
    &~+\lambda_{n,1}\parenBig{\lone{\beta^*_{(1),S}+\Delta_S}+\lone{\Delta_{S^c}}-\lone{\beta^*_{(1),S}}-2\lone{\beta^*_{(1),S^c}}}\\
    \ge&~-\lone{\Delta}\linf{\grad_\beta L_n\paren{Y,\<\beta^*_{(1)},X\>}}+C_1\ltwo{\Delta}\parenBig{\ltwo{\Delta}-C_2\sqrt{\frac{\log p_1}{n}}\lone{\Delta}}\\
    &~+\lambda_{n,1}\parenBig{\lone{\Delta_{S^c}}-\lone{\Delta_S}-2\lone{\beta^*_{(1),S^c}}}\\
    \ge&~\parenBig{\frac{\lambda_{n,1}}{2}-C_1C_2\tau\sqrt{\frac{\log p_1}{n}}}\lone{\Delta_{S^c}}-\parenBig{\frac{3\lambda_{n,1}}{2}+C_1C_2\tau\sqrt{\frac{\log p_1}{n}}}\lone{\Delta_{S}}\\
    &~+C_1\ltwo{\Delta}^2-2\lambda_{n,1}\lone{\beta^*_{(1),S^c}}\\
    \ge&~-2\lambda_{n,1}\lone{\Delta_{S}}+C_1\ltwo{\Delta}^2-2\lambda_{n,1}\lone{\beta^*_{(1),S^c}}\\
    \ge&~-2\sqrt{\abs{S}}\lambda_{n,1}\ltwo{\Delta}+C_1\ltwo{\Delta}^2-2\lambda_{n,1}\lone{\beta^*_{(1),S^c}},
\end{aligned}
\end{equation}
where we choose $\tau\in(0,1]$ such that $C_1C_2\tau\sqrt{\frac{\log p_1}{n}}\leq \frac{\lambda_{n,1}}{2}$.
If $\lambda_{n,1}<\frac{C_1}{2\sqrt{\abs{S}}\tau+2\lone{\beta^*_{(1),S^c}}}\tau^2$, then under the events $\mc{B}_1,\mc{B}_2$, $\mc{F}\paren{\Delta}>0$ for any $\Delta \in \setBig{\Delta\in\R^{p_1}\Big|\ltwo{\Delta}= \tau}$. We claim that under $\mc{B}_1$ and $\mc{B}_2$, $\ltwo{\hat{\Delta}}\overset{\triangle}{=}\ltwo{\hat{\beta}_{(1)}-\beta^*_{(1)}}\leq \tau\leq1$. In fact, if $\ltwo{\hat{\Delta}}>\tau$, then by the convexity of $\mc{F}$ and $\mc{F}\paren{0}=0$, we have
\begin{align*}
    0<\mc{F}\paren{\frac{\tau}{\ltwo{\hat{\Delta}}}\hat{\Delta}+ \paren{1-\frac{\tau}{\ltwo{\hat{\Delta}}}}\times 0}\leq \frac{\tau}{\ltwo{\hat{\Delta}}}\mc{F}\paren{\hat{\Delta}},
\end{align*}
which is contradictory with the fact that $\mc{F}\paren{\hat{\Delta}}\leq 0$.
Thus under the events $\mc{B}_1$ and $\mc{B}_2$, by \eqref{proof:prop1-1}, we have
\begin{align*}
   0\ge\mc{F}\paren{\hat{\Delta}}\ge-2\sqrt{\abs{S}}\lambda_{n,1}\ltwo{\hat{\Delta}}+C_1\ltwo{\hat{\Delta}}^2-2\lambda_{n,1}\lone{\beta^*_{(1),S^c}},
\end{align*}
which implies that (since the right-hand side is a quadratic function of $\ltwo{\hat{\Delta}}$)
\begin{align}\label{proof:prop1-2}
\ltwo{\hat{\Delta}}\leq\frac{\sqrt{\abs{S}}\lambda_{n,1}+\sqrt{\abs{S}\lambda_{n,1}^2+2C_1\lone{\beta^*_{(1),S^c}}\lambda_{n,1}}}{C_1}\leq C_3\sqrt{\abs{S}}\lambda_{n,1}+C_3\sqrt{\lambda_{n,1}\lone{\beta^*_{(1),S^c}}}
\end{align} holds with probability at least $1-2\exp\paren{-n\delta^2}-2\exp\paren{-cn}$. By the convexity of the loss function, under the event $\mc{B}_2$, we have
\begin{align*}
    0\ge&~L_n\paren{Y,\<\beta^*_{(1)}+\hat{\Delta},X\>}-L_n\paren{Y,\<\beta^*_{(1)},X\>}+\lambda_{n,1}\parenBig{\lone{\beta^*_{(1)}+\hat{\Delta}}-\lone{\beta^*_{(1)}}}\\
    \ge&~\<\hat{\Delta},\grad_\beta L_n\paren{Y,\<\beta^*_{(1)},X\>}\>+\lambda_{n,1}\parenBig{\lone{\beta^*_{(1),S}+\hat{\Delta}_S}+\lone{\hat{\Delta}_{S^c}}-\lone{\beta^*_{(1),S}}-2\lone{\beta^*_{(1),S^c}}}\\
    \ge&~-\lone{\hat{\Delta}}\linf{\grad_\beta L_n\paren{Y,\<\beta^*_{(1)},X\>}}+\lambda_{n,1}\parenBig{\lone{\hat{\Delta}_{S^c}}-\lone{\hat{\Delta}_S}-2\lone{\beta^*_{(1),S^c}}}\\
    \ge&~\frac{1}{2}\lambda_{n,1}\lone{\hat{\Delta}_{S^c}}-\frac{3}{2}\lambda_{n,1}\lone{\hat{\Delta}_S}-2\lambda_{n,1}\lone{\beta^*_{(1),S^c}},
\end{align*}
which implies $\hat{\Delta}\in\setBig{\Delta\in\R^{p_1}\Big| ~3\lone{\Delta_S}+4\lone{\beta^*_{(1),S^c}}\ge\lone{\Delta_{S^c}}}$. Consequently, the $\ell_1$ bound follows from the inequality
\begin{align}\label{proof:prop1-3}
    \lone{\hat{\Delta}}\leq 4\lone{\hat{\Delta}_S}+4\lone{\beta^*_{(1),S^c}}\leq 4\sqrt{\abs{S}}\ltwo{\hat{\Delta}}+4\lone{\beta^*_{(1),S^c}},
\end{align}
under the event $\mc{B}_2$. 

To derive the result, we consider the following two cases.
\begin{enumerate}
\item Letting $S=\text{supp}(\beta^*_{(1)})$. Then $\lone{\beta^*_{(1),S^c}}=0$ and by \eqref{proof:prop1-2} and \eqref{proof:prop1-3} we have
\begin{align*}
    \ltwo{\hat{\Delta}}\leq C\sqrt{s_1}\lambda_{n,1}, \lone{\hat{\Delta}}\leq Cs_1\lambda_{n,1}
\end{align*}
with probability at least $1-2\exp\paren{-n\delta^2}-2\exp\paren{-cn}$ if $\sqrt{s_1}\lambda_{n,1}\leq C$, where $s_1=\lzero{\beta^*_{(1)}}$.
    \item Letting $S=\set{j\in\brac{p_1}|~\abs{\beta_{(1),j}^*}\ge\lambda_{n,1}}$. Then $\abs{S}\lambda_{n,1}^q\leq h_{q,1}^q$ and $\lone{\beta^*_{(1),S^c}}/\lambda_n\leq h_{q,1}^q/\lambda_{n,1}^q$. Thus, by \eqref{proof:prop1-2} and \eqref{proof:prop1-3}
    \begin{align*}
        \ltwo{\hat{\Delta}}\leq C\lambda_{n,1}^{1-\frac{q}{2}}h_{q,1}^{\frac{q}{2}}, \lone{\hat{\Delta}}\leq C\lambda^{1-q}_{n,1}h_{q,1}^q
    \end{align*}
    with probability at least $1-2\exp\paren{-n\delta^2}-2\exp\paren{-cn}$ if $\lambda_{n,1}^{1-\frac{q}{2}}h_{q,1}^{\frac{q}{2}}\leq C$, where $h_{q,1}=\norm{\beta^*_{(1)}}_q$.
\end{enumerate}

Now, we derive the result for the general regularizers that satisfy \cref{con:regularizer}. Similar to \eqref{proof:prop1-1}, under the events $\mc{B}_1, \mc{B}_2$, for any $S\subseteq\brac{p_1}$ and $\Delta \in \setBig{\Delta\in\R^{p_1}\Big|\ltwo{\Delta}\leq \tau}$, by triangle inequality and \cref{lemma:regularizer-to-lone}, we have
\begin{equation}\label{proof:prop1-4}
    \begin{aligned}
    \mc{F}^\prime(\Delta)\overset{\triangle}{=}&~L_n\paren{Y,\<\beta^*_{(1)}+\Delta,X\>}-L_n\paren{Y,\<\beta^*_{(1)},X\>}+r_{\lambda_{n,1}}(\beta^*_{(1)}+\Delta)-r_{\lambda_{n,1}}(\beta^*_{(1)})\\
    \ge&~-\frac{\lambda_{n,1}}{2}\lone{\Delta}+C_1\ltwo{\Delta}\parenBig{\ltwo{\Delta}-C_2\sqrt{\frac{\log p_1}{n}}\lone{\Delta}}\\
    &~+\lambda_{n,1}\lone{\Delta_{S^c}}-\lambda_{n,1}\lone{\Delta_S}-\frac{\kappa_r}{2}\ltwo{\Delta_{S^c}}^2-2\lambda_{n,1}\lone{\beta^*_{(1),S^c}}\\
    \ge&~-2\sqrt{\abs{S}}\lambda_{n,1}\ltwo{\Delta}+(C_1-\frac{\kappa_r}{2})\ltwo{\Delta}^2-2\lambda_{n,1}\lone{\beta^*_{(1),S^c}}\\
    \ge&~-2\sqrt{\abs{S}}\lambda_{n,1}\ltwo{\Delta}+\frac{C_1}{2}\ltwo{\Delta}^2-2\lambda_{n,1}\lone{\beta^*_{(1),S^c}},
\end{aligned}
\end{equation}
where the last inequality uses $C_1\ge c\kappa_l\ge \kappa_r$ from the proof of \cref{lemma:RSC} and \cref{con:regularizer} (iii) for some $c>0$ only depending on $K_\beta$ and the loss function. We claim that under $\mc{B}_1, \mc{B}_2$, $\ltwo{\hat{\Delta}}\leq\tau$. If not, taking $t=\frac{\tau}{\ltwo{\hat{\Delta}}}<1$, by $\lone{\hat{\Delta}}\leq R$, Lemma 4 in \cite{loh2015regularized} and the convexity of the loss function, we have
\begin{align*}
    0\ge&~L_n\paren{Y,\<\beta^*_{(1)}+\hat{\Delta},X\>}-L_n\paren{Y,\<\beta^*_{(1)},X\>}+r_{\lambda_{n,1}}(\beta^*_{(1)}+\hat{\Delta})-r_{\lambda_{n,1}}(\beta^*_{(1)})\\
    \ge&~\frac{1}{t}\parenBig{L_n\paren{Y,\<\beta^*_{(1)}+t\hat{\Delta},X\>}-L_n\paren{Y,\<\beta^*_{(1)},X\>}}-\lambda_{n,1}\lone{\hat{\Delta}}\\
    \ge&~-2\lambda_{n,1}\lone{\hat{\Delta}}+C_1\tau \ltwo{\hat{\Delta}}-\tau C_1C_2\sqrt{\frac{\log p_1}{n}}\lone{\hat{\Delta}}\\
    \ge&~C_1\tau\ltwo{\hat{\Delta}}-(2\lambda_{n,1}+\tau  C_1C_2\sqrt{\frac{\log p_1}{n}})R,
\end{align*}
which is impossible when $(2\lambda_{n,1}+\tau C_1C_2\sqrt{\frac{\log p_1}{n}})R\leq C_1\tau^2$. Then, by \eqref{proof:prop1-4}, under $\mc{B}_1,\mc{B}_2$, we have
\begin{align}\label{proof:prop1-5}
    \ltwo{\hat{\Delta}}\leq C_3\sqrt{\abs{S}}\lambda_{n,1}+C_3\sqrt{\lambda_{n,1}\lone{\beta^*_{(1),S^c}}}.
\end{align}
We also note that by \eqref{proof:prop1-4}, under $\mc{B}_1,\mc{B}_2$, choosing $\tau\in(0,1]$ such that $C_1C_2\tau\sqrt{\frac{\log p_1}{n}}\leq \frac{\lambda_{n,1}}{4}$, we have
\begin{align*}
    0\ge&~-\frac{\lambda_{n,1}}{2}\lone{\hat{\Delta}}+C_1\ltwo{\hat{\Delta}}\parenBig{\ltwo{\hat{\Delta}}-C_2\sqrt{\frac{\log p_1}{n}}\lone{\hat{\Delta}}}\\
    &~+\lambda_{n,1}\lone{\hat{\Delta}_{S^c}}-\lambda_{n,1}\lone{\hat{\Delta}_S}-\frac{\kappa_r}{2}\ltwo{\hat{\Delta}_{S^c}}^2-2\lambda_{n,1}\lone{\beta^*_{(1),S^c}}\\
    \ge&~\frac{1}{4}\lambda_{n,1}\lone{\hat{\Delta}_{S^c}}-\frac{7}{4}\lambda_{n,1}\lone{\hat{\Delta}_S}+\frac{C_1}{2}\ltwo{\hat{\Delta}}^2-2\lambda_{n,1}\lone{\beta^*_{(1),S^c}}\\
    \ge&~\frac{1}{4}\lambda_{n,1}\lone{\hat{\Delta}_{S^c}}-\frac{7}{4}\lambda_{n,1}\lone{\hat{\Delta}_S}-2\lambda_{n,1}\lone{\beta^*_{(1),S^c}},
\end{align*}
which implies that $\hat{\Delta}\in\setBig{\Delta\in\R^{p_1}\Big| ~7\lone{\Delta_S}+8\lone{\beta^*_{(1),S^c}}\ge\lone{\Delta_{S^c}}}$. Consequently, the $\ell_1$ bound follows from the inequality
\begin{align}\label{proof:prop1-6}
    \lone{\hat{\Delta}}\leq 8\lone{\hat{\Delta}_S}+8\lone{\beta^*_{(1),S^c}}\leq 8\sqrt{\abs{S}}\ltwo{\hat{\Delta}}+8\lone{\beta^*_{(1),S^c}},
\end{align}
under the event $\mc{B}_1, \mc{B}_2$. Combining \eqref{proof:prop1-5} and \eqref{proof:prop1-6}, the rest of the proof is similar to the proof for the Lasso penalty.

\end{proof}
%%%%%%%%%%%%%%%%%%%%%%%%%%%%%%%%%%%%%%%%%%%%%%%%%%%%%%%%%%%
\subsection{Proof of \cref{unpenalized}}
\begin{proof} Without loss of generality, we assume $1\in\mc{A}^d$ and prove the result for the case $k=1$. Analogous to the proof of \cref{lasso}, we assume that $\set{X_i}_{i\in\brac{n}}$, $\beta$, $\beta^*_{(1)}$ and $\hat{\beta}_{(1)}$ select only the elements from the set $\mc{A}_1 \subseteq \brac{p}$.  Since $\lone{\hat{\beta}_{(1)}-\beta^*_{(1)}}\leq \sqrt{d}\ltwo{\hat{\beta}_{(1)}-\beta^*_{(1)}}$, it suffices to show that there exists a sufficiently large constant $C>0$ such that
\begin{align*}
    &~\P\parenBig{\inf_{\ltwo{\beta-\beta^*_{(1)}}=C\sqrt{d}(\delta+\sqrt{\frac{\log d}{n}})} L_n(Y,\<\beta,X\>)-L_n(Y,\<\beta^*_{(1)},X\>)>0}\\
    >&~1-2\exp\paren{-n\delta^2}-2\exp\paren{-c_2n},
\end{align*}
where $\delta$ is a small positive constant such that $C\sqrt{d}(\delta+\sqrt{\frac{\log d}{n}})\leq 1$. Following the proof of the inequality \eqref{eq:infnormbound} in \cref{lasso} and \cref{lemma:RSC}, with probability $1-2\exp\paren{-n\delta^2}-2\exp\paren{-c_2n}$,
\begin{align*}
    &~L_n\paren{Y,\<\beta^*_{(1)}+\Delta,X\>}-L_n\paren{Y,\<\beta^*_{(1)},X\>}\\
    \ge&~\<\Delta,\grad_\beta L_n\paren{Y,\<\beta^*_{(1)},X\>}\>+C_1\ltwo{\Delta}^2-C_2\sqrt{\frac{\log d}{n}}\lone{\Delta}\ltwo{\Delta}\\
    \ge&~-\linf{\grad_\beta L_n\paren{Y,\<\beta^*_{(1)},X\>}}\lone{\Delta}+\frac{C_1}{2}\ltwo{\Delta}^2\\
    \ge&~-C_3\sqrt{d}(\sqrt{\frac{\log d}{n}}+\delta)\ltwo{\Delta}+\frac{C_1}{2}\ltwo{\Delta}^2
\end{align*}
holds uniformly for any $\Delta\in\set{\Delta\in\R^{\abs{\mc{A}_1}}|~\ltwo{\Delta}\leq 1}$ , where we take large enough $n$ such that $ 2C_2\sqrt{\frac{d\log d}{n}}\leq C_1$ and $C_1, C_2, C_3$ are some constants not related to $n$ and $d$. By choosing $C=4C_3/C_1$, when $C\sqrt{d}(\delta+\sqrt{\frac{\log d}{n}})\leq 1$, we have
\begin{align*}
    \inf_{\ltwo{\beta-\beta^*_{(1)}}=C\sqrt{d}(\sqrt{\frac{\log d}{n}}+\delta)} L_n(Y,\<\beta,X\>)-L_n(Y,\<\beta^*_{(1)},X\>)
    \ge\frac{4C_3^2}{C_1}d(\sqrt{\frac{\log d}{n}}+\delta)^2
    >0,
\end{align*}
with probability $1-2\exp\paren{-n\delta^2}-2\exp\paren{-c_2n}$, which completes the proof.
\end{proof}
%%%%%%%%%%%%%%%%%%%%%%%%%%%%%%%%%%%%%%%%%%%%%%%%%%%%%%%%%%
\subsection{Proof of \cref{prop:sparse} in \cref{remark:rate of w}}
\begin{proof} Let $s_w=\abs{\text{supp}(w^*)}$. By \eqref{eq:support-of-w} and \eqref{eq:equicorrelation}, for some $v>0$ we have
\begin{align*}
    (1-w^*_{K_{ne}})\ltwo{\beta_{(K_{ne})}}^2=\cdots=(1-w^*_{K_{ne}+s_w-1})\ltwo{\beta_{(K_{ne}+s_w-1)}}^2=v.
\end{align*}
Since $\lone{w^*}=1$, we have 
\begin{align*}
    v=\frac{s_w-1}{\sum_{k=K_{ne}}^{K_{ne}+s_w-1}\frac{1}{\ltwo{\beta^*_{(k)}}^2}}=(1-w^*_{K_{ne}+s_w-1})\ltwo{\beta^*_{(K_{ne}+s_w-1)}}^2\leq\ltwo{\beta^*_{(K_{ne}+s_w-1)}}^2.
\end{align*}
Assume that $s_w\ge3$. By the condition $\frac{\ltwo{\beta^*_{(K_{ne}+2)}}^2}{\ltwo{\beta^*_{(K_{ne})}}^2}+\frac{\ltwo{\beta^*_{(K_{ne}+2)}}^2}{\ltwo{\beta^*_{(K_{ne}+1)}}^2}<1$ and $\ltwo{\beta^*_{(K_{ne})}}\ge\cdots\ge\ltwo{\beta^*_{(K)}}$, we can obtain that
\begin{align*}
    s_w\leq1+\sum_{k=K_{ne}}^{K_{ne}+s_w-1}\frac{\ltwo{\beta^*_{(K_{ne}+s_w-1)}}^2}{\ltwo{\beta^*_{(k)}}^2}\leq \frac{\ltwo{\beta^*_{(K_{ne}+2)}}^2}{\ltwo{\beta^*_{(K_{ne})}}^2}+\frac{\ltwo{\beta^*_{(K_{ne}+2)}}^2}{\ltwo{\beta^*_{(K_{ne}+1)}}^2}+s_w-1<s_w,
\end{align*}
which is impossible. Therefore, $s_w=2.$
\end{proof}
\section{Proof of Theorem}\label{appendix:proof-of-theorem}

\subsection{Proof of \cref{MA}}
\begin{proof} Denote 
\begin{align*}
    T_{n,q}^{(3)}(\delta)=\begin{cases}
        2K_\beta\sqrt{s_u}&~~\text{if }q=0\\
        K_\beta h_{q,u}^{\frac{q}{2}}(\delta+\sqrt{\frac{\log p_u}{n}})^{-\frac{q}{2}}+(K_\beta+1)h_{q,u}^{\frac{2}{2-q}}&~~\text{if }q\in(0,1)\\
    \end{cases}
\end{align*}and 
    \begin{align*}
        \mc{B}_3=\setBig{\max_{m\in\brac{J}}\max_{k\in\brac{K}}\ltwo{\hat{\beta}_{(k)}^{\brac{-m}}-\beta_{(k)}^*}\leq T_{n,q}^{(2)}(\delta)}\bigcap\setBig{\max_{m\in\brac{J}}\max_{k\in\brac{K}}\lone{\hat{\beta}_{(k)}^{\brac{-m}}-\beta_{(k)}^*}\leq T_{n,q}^{(1)}(\delta)}.
    \end{align*}
    By \cref{betarate} and \cref{HD-MA:beta-deviation}, if $\delta\in(0,c]$, $T_{n,q}^{(2)}(\delta)\leq K_\beta$, then under $\mc{B}_3$, for any $\lone{\Delta}\leq 2$, we have 
    \begin{align*}
\max_{m\in\brac{J}}\lone{\hat{B}^{\brac{-m}}\Delta}&\leq T_{n,q}^{(3)}(\delta)\lone{\Delta}\\
\max_{m\in\brac{J}}\ltwo{\hat{B}^{\brac{-m}}\Delta}&\leq 2K_\beta\lone{\Delta}\leq 4K_\beta.
\end{align*} and $\mc{B}_3$ holds with probability at least $1-4JK\exp\paren{-n_{-1}\delta^2}$.

 For each $m\in\brac{J}$, we will prove that the following event  $\mc{B}_{4,m}$
{\footnotesize \begin{align*}
    &~\frac{1}{n_1}\sum_{i\in\mc{I}_m}\Big\{L\paren{Y_i,\<\hat{B}^{\brac{-m}}\paren{w^*+\Delta},X_i\>}-L\paren{Y_i,\<\hat{B}^{\brac{-m}}w^*,X_i\>}-\frac{\partial L\paren{Y_i,\<\hat{B}^{\brac{-m}}w^*,X_i\>}}{\partial \mu}\<\hat{B}^{\brac{-m}}\Delta,X_i\>\Big\}\\
    \ge&~C_1\ltwo{\hat{B}^{\brac{-m}}\Delta}\parenBig{\ltwo{\hat{B}^{\brac{-m}}\Delta}-C_2\sqrt{\frac{\log p}{n_1}}\lone{\hat{B}^{\brac{-m}}\Delta}}\text{\ \ \ for any\ }\Delta\in\setBig{\Delta\in\R^K\Big|\lone{\Delta}\leq 2}
\end{align*}}
holds with high probability. By \cref{lemma:RSC} with $\beta_0=\hat{B}^{\brac{-m}}w^*, R_1=4K_\beta, R_2=0$, if $n_1\ge c_1$ and $T_{n,q}^{(2)}(\delta)\leq K_\beta$, we have 
    \begin{align*}
    \E\parenBig{\P\paren{\mc{B}_{4,m}^c|\hat{B}^{\brac{-m}}}\Big|\mc{B}_{3,m}}\leq 2\exp\paren{-c_2n_1},
\end{align*}
where
\begin{align}\label{proof:weight-B3m}
    \mc{B}_{3,m}=\setBig{\max_{k\in\brac{K}}\ltwo{\hat{\beta}_{(k)}^{\brac{-m}}-\beta_{(k)}^*}\leq T_{n,q}^{(2)}(\delta)}\bigcap\setBig{\max_{k\in\brac{K}}\lone{\hat{\beta}_{(k)}^{\brac{-m}}-\beta_{(k)}^*}\leq T_{n,q}^{(1)}(\delta)}.
\end{align}
Since $\mc{B}_3\subseteq\mc{B}_{3,m}$, then $\P(\mc{B}_{3,m}^c)\leq\P(\mc{B}_3^c)$. Thus,
    \begin{align*}
        \P\paren{\mc{B}_{4,m}^c}=&~\E\parenBig{\P\paren{\mc{B}_{4,m}^c|\hat{B}^{\brac{-m}}}}\leq\E\parenBig{\P\paren{\mc{B}_{4,m}^c|\hat{B}^{\brac{-m}}}\Big|\mc{B}_{3,m}}+\P\parenBig{\mc{B}_3^c}\\
        \leq&~ 4JK\exp\paren{-n_{-1}\delta^2}+2\exp\paren{-c_2n_1}.
    \end{align*}
Consequently, the following event $\mc{B}_4\supseteq \cap_{m\in\brac{J}}\mc{B}_{4,m}$
{\footnotesize\begin{align*}
    &~\dn\sum_{m=1}^J\sum_{i\in\mc{I}_m}\setBig{L\paren{Y_i,\<\hat{B}^{\brac{-m}}\paren{w^*+\Delta},X_i\>}-L\paren{Y_i,\<\hat{B}^{\brac{-m}}w^*,X_i\>}-\frac{\partial L\paren{Y_i,\<\hat{B}^{\brac{-m}}w^*,X_i\>}}{\partial \mu}\<\hat{B}^{\brac{-m}}\Delta,X_i\>}\\
    \ge&~\frac{1}{J}\sum_{m=1}^JC_1\ltwo{\hat{B}^{\brac{-m}}\Delta}\parenBig{\ltwo{\hat{B}^{\brac{-m}}\Delta}-C_2\sqrt{\frac{\log p}{n_1}}\lone{\hat{B}^{\brac{-m}}\Delta}}\text{\ \ \ for any\ }\Delta\in\setBig{\Delta\in\R^K\Big|\lone{\Delta}\leq 2}
\end{align*}}
 holds with probability at least $1-4J^2K\exp\paren{-n_{-1}\delta^2}-2J\exp\paren{-c_2n_1}$.
 
Denote 
\begin{align*}
    V_{i,k_1}^{\brac{-m}}\overset{\triangle}{=}\frac{\partial L\paren{Y_i,\<\hat{B}^{\brac{-m}}w^*,X_i\>}}{\partial \mu}\<\hat{\beta}_{(k_1)}^{\brac{-m}},X_i\>-\E\parenBig{\frac{\partial L\paren{Y_i,\<B^*w^*,X_i\>}}{\partial \mu}\<\beta^*_{(k_1)},X_i\>},
\end{align*}
where $i\in\mc{I}_m$ and $k_1\in\brac{K}$. To control $\absBig{\dn\sum_{m\in\brac{J}}\sum_{i\in\mc{I}_m}V_{i,k_1}^{\brac{-m}}}$ for fixed $k_1\in\brac{K}$, note that
{\small\begin{align*}
    &~\absBig{\dn\sum_{m\in\brac{J}}\sum_{i\in\mc{I}_m}V_{i,k_1}^{\brac{-m}}}\\
    \leq&~\frac{1}{J}\sum_{m\in\brac{J}}\sup_{w\in\mc{W}^K}\absBig{\frac{1}{n_1}\sum_{i\in\mc{I}_m}\frac{\partial L\paren{Y_i,\<\hat{B}^{\brac{-m}}w,X_i\>}}{\partial \mu}\<\hat{\beta}_{(k_1)}^{\brac{-m}},X_i\>-\frac{\partial L\paren{Y_i,\<B^*w,X_i\>}}{\partial \mu}\<\hat{\beta}_{(k_1)}^{\brac{-m}},X_i\>}\\
    &~+\frac{1}{J}\sum_{m\in\brac{J}}\absBig{\frac{1}{n_1}\sum_{i\in\mc{I}_m}\frac{\partial L\paren{Y_i,\<B^*w^*,X_i\>}}{\partial \mu}\<\hat{\beta}_{(k_1)}^{\brac{-m}},X_i\>-\E\parenBig{\frac{\partial L\paren{Y_i,\<B^*w^*,X_i\>}}{\partial \mu}\<\hat{\beta}_{(k_1)}^{\brac{-m}},X_i\>\Big|\hat{B}^{\brac{-m}}}}\\
    &~+ \frac{1}{J}\sum_{m\in\brac{J}}\absBig{\E\parenBig{\frac{\partial L\paren{Y_i,\<B^*w^*,X_i\>}}{\partial \mu}\<\hat{\beta}_{(k_1)}^{\brac{-m}},X_i\>\Big|\hat{B}^{\brac{-m}}}-\E\parenBig{\frac{\partial L\paren{Y_i,\<B^*w^*,X_i\>}}{\partial \mu}\<\beta^*_{(k_1)},X_i\>}}\\
    \overset{\triangle}{=}&~\frac{1}{J}\sum_{m\in\brac{J}}(I_{1,m,k_1}+I_{2,m,k_1}+I_{3,m,k_1}).
\end{align*}}
For any $k_1\in\brac{K}$, we have
{\footnotesize\begin{align*}
    I_{1,m,k_1}\leq&~\sup_{w\in\mc{W}^K}\absBig{\frac{1}{n_1}\sum_{i\in\mc{I}_m}\frac{\partial L\paren{Y_i,\<\hat{B}^{\brac{-m}}w,X_i\>}}{\partial \mu}\<\hat{\beta}_{(k_1)}^{\brac{-m}},X_i\>-\frac{\partial L\paren{Y_i,\<B^*w,X_i\>}}{\partial \mu}\<\hat{\beta}_{(k_1)}^{\brac{-m}},X_i\>}\\
    \leq &~\sup_{w\in\mc{W}^K}\frac{1}{n_1}\sum_{i\in\mc{I}_m}L_u\abs{\<\hat{\beta}_{(k_1)}^{\brac{-m}},X_i\>}\abs{\<\hat{B}^{\brac{-m}}w-B^*w,X_i\>}\\
    \leq&~~\sup_{k_2\in\brac{K}}\frac{1}{n_1}\sum_{i\in\mc{I}_m}L_u\abs{\<\hat{\beta}_{(k_1)}^{\brac{-m}},X_i\>}\abs{\<\hat{\beta}_{(k_2)}^{\brac{-m}}-\beta^*_{(k_2)},X_i\>}\\
    \leq&~\sup_{k_2\in\brac{K}}\absBig{\frac{1}{n_1}\sum_{i\in\mc{I}_m}L_u\abs{\<\hat{\beta}_{(k_1)}^{\brac{-m}},X_i\>}\abs{\<\hat{\beta}_{(k_2)}^{\brac{-m}}-\beta^*_{(k_2)},X_i\>}-\E \parenBig{L_u\abs{\<\hat{\beta}_{(k_1)}^{\brac{-m}},X_i\>}\abs{\<\hat{\beta}_{(k_2)}^{\brac{-m}}-\beta^*_{(k_2)},X_i\>}\Big|\hat{B}^{\brac{-m}}}}\\
    &~+\sup_{k_2\in\brac{K}}\E \parenBig{L_u\abs{\<\hat{\beta}_{(k_1)}^{\brac{-m}},X_i\>}\abs{\<\hat{\beta}_{(k_2)}^{\brac{-m}}-\beta^*_{(k_2)},X_i\>}\Big|\hat{B}^{\brac{-m}}}\\
    \overset{\triangle}{=}&~I_{4,m,k_1}+I_{5,m,k_1}.
\end{align*}}
Denote
\begin{align*}
    U_{i,k_1,k_2}^{\brac{-m}}\overset{\triangle}{=}L_u\abs{\<\hat{\beta}_{(k_1)}^{\brac{-m}},X_i\>}\abs{\<\hat{\beta}_{(k_2)}^{\brac{-m}}-\beta^*_{(k_2)},X_i\>},
\end{align*}
where $i\in\mc{I}_m$ and $k_1,k_2\in\brac{K}$. Conditioning on $\hat{B}^{\brac{-m}}$, by Lemma 2.7.7 of \cite{vershynin2018high} we have
\begin{align*}
    \norm{\frac{U_{i,k_1,k_2}^{\brac{-m}}}{L_uK_x^2\ltwo{\hat{\beta}_{(k_1)}^{\brac{-m}}}\sup_{k_2\in\brac{K}}\ltwo{\hat{\beta}_{(k_2)}^{\brac{-m}}-\beta^*_{(k_2)}}}}_{\psi_1}\leq 1
\end{align*}
and then by Proposition 2.7.1 of \cite{vershynin2018high}
\begin{align}
   I_{5,m,k_1}= \sup_{k_2\in\brac{K}}\E \parenBig{U_{i,k_1,k_2}^{\brac{-m}}\Big|\hat{B}^{\brac{-m}}}\leq C\ltwo{\hat{\beta}_{(k_1)}^{\brac{-m}}}\sup_{k_2\in\brac{K}}\ltwo{\hat{\beta}_{(k_2)}^{\brac{-m}}-\beta^*_{(k_2)}}.\label{proof:consistency-I5}
\end{align}
By Bernstein's inequality, \cref{centering} and taking expectation, for any $t_1\in(0,1)$,
\begin{align}
    \P\parenBig{I_{4,m,k_1}>C\ltwo{\hat{\beta}_{(k_1)}^{\brac{-m}}}\sup_{k_2\in\brac{K}}\ltwo{\hat{\beta}_{(k_2)}^{\brac{-m}}-\beta^*_{(k_2)}}t_1}\leq2K\exp(-cn_1t_1^2).\label{proof:consistency-I4}
\end{align}
Similarly, for any $t_2\in(0,1)$, $k_1\in\brac{K}$, we have
\begin{align}
   \P \parenBig{I_{2,m,k_1}>C\ltwo{\hat{\beta}_{(k_1)}^{\brac{-m}}}t_2}\leq2\exp(-cn_1t_2^2),\label{proof:consistency-I2}
\end{align}
since
\begin{align*}
    \norm{\frac{\partial L\paren{Y_i,\<B^*w^*,X_i\>}}{\partial \mu}}_{\psi_2}\leq L_uK_\beta K_x+K_y.
\end{align*}
By Lemma 2.7.7 and Proposition 2.7.1 in \cite{vershynin2018high}, we have
\begin{align}
    I_{3,m,k_1}\leq C\ltwo{\hat{\beta}_{(k_1)}^{\brac{-m}}-\beta^*_{(k_1)}}.\label{proof:consistency-I3}
\end{align}
Combined with \eqref{proof:consistency-I5}, \eqref{proof:consistency-I4}, \eqref{proof:consistency-I2} and \eqref{proof:consistency-I3}, with probability at least $1-2JK^2\exp(-cn_1t_1^2)+2JK\exp(-cn_1t_2^2)$, we have
\begin{align*}
&~\sup_{k_1\in\brac{K}}\absBig{\dn\sum_{m\in\brac{J}}\sum_{i\in\mc{I}_m}V_{i,k_1}^{\brac{-m}}}\\
    \leq&~C\sup_{m\in\brac{J}}\Big\{\sup_{k_1\in\brac{K}}\ltwo{\hat{\beta}_{(k_1)}^{\brac{-m}}}\parenBig{\sup_{k_2\in\brac{K}}\ltwo{\hat{\beta}_{(k_2)}^{\brac{-m}}-\beta^*_{(k_2)}}(t_1+1)+t_2}+\sup_{k_1\in\brac{K}}\ltwo{\hat{\beta}_{(k_1)}^{\brac{-m}}-\beta^*_{(k_1)}}\Big\}.
\end{align*}
Note that under the event $\mc{B}_3$, if $T_{n,q}^{(2)}(\delta)\leq K_\beta$, we have
 \begin{align*}
     &~\sup_{m\in\brac{J}}\setBig{\sup_{k_1\in\brac{K}}\ltwo{\hat{\beta}_{(k_1)}^{\brac{-m}}}\parenBig{\sup_{k_2\in\brac{K}}\ltwo{\hat{\beta}_{(k_2)}^{\brac{-m}}-\beta^*_{(k_2)}}(t_1+1)+t_2}+\sup_{k_1\in\brac{K}}\ltwo{\hat{\beta}_{(k_1)}^{\brac{-m}}-\beta^*_{(k_1)}}}\\
     \leq&~ C\parenBig{(t_1+1)T_{n,q}^{(2)}(\delta)+t_2}.
 \end{align*}
Let $\lambda_n^\prime=C\parenBig{(t_1+1)T_{n,q}^{(2)}(\delta)+t_2}$. Then the following event $\mc{B}_5$
\begin{align*}
    \setBig{\lambda_n^\prime\ge\linf{\dn\sum_{m=1}^J\sum_{i\in \mc{I}_m}\frac{\partial L\paren{Y_i,\<\hat{B}^{\brac{-m}}w^*,X_i\>}}{\partial \mu}\hat{B}^{\brac{-m}\top}X_i-\E\parenBig{\frac{\partial L\paren{Y_i,\<B^*w^*,X_i\>}}{\partial\mu}B^{*\top}X_i}}}
\end{align*} holds with probability at least $1-2JK^2\exp(-cn_1t_1^2)-2JK\exp(-cn_1t_2^2)-4JK\exp\paren{-n_{-1}\delta^2}$. 

Let $\hat{\Delta}=\hat{w}-w^*$. By the optimality condition, we have $\<\E\parenBig{\frac{\partial L\paren{Y_i,\<B^*w^*,X_i\>}}{\partial\mu}B^{*\top}X_i},\hat{\Delta}\>\ge 0$. Note that under the events $\mc{B}_5$, for any $\varrho>1$ and $S_*\subseteq \brac{K}$, using the property $\lone{\hat{w}}=\lone{w^*}=1$ we have
\begin{align*}
    0\ge&~\dn\sum_{m=1}^J\sum_{i\in\mc{I}_m}\setBig{L\paren{Y_i,\<\hat{B}^{\brac{-m}}\paren{w^*+\hat{\Delta}},X_i\>}-L\paren{Y_i,\<\hat{B}^{\brac{-m}}w^*,X_i\>}}\\
    &~+\varrho\lambda_n^\prime\parenBig{\lone{w^*+\hat{\Delta}}-\lone{w^*}}\\
    \ge&~\<\dn\sum_{m=1}^J\sum_{i\in \mc{I}_m}\frac{\partial L\paren{Y_i,\<\hat{B}^{\brac{-m}}w^*,X_i\>}}{\partial \mu}\hat{B}^{\brac{-m}\top}X_i-\E\parenBig{\frac{\partial L\paren{Y_i,\<B^*w^*,X_i\>}}{\partial\mu}B^{*\top}X_i},\hat{\Delta}\>\\
    &~+\varrho\lambda_n^\prime\parenBig{\lone{w^*+\hat{\Delta}}-\lone{w^*}}\\
    \ge&~-\lambda_n^\prime\lone{\hat{\Delta}}+\varrho\lambda_n^\prime\parenBig{\lone{w^*_{S_*}+\hat{\Delta}_{S_*}}+\lone
    {\hat{\Delta}_{S_*^c}}-\lone{w^*_{S_*}}-2\lone{w^*_{S_*^c}}}\\
    \ge&~-\lambda_n^\prime \lone{\hat{\Delta}}+\varrho\lambda_n^\prime\parenBig{\lone{\hat{\Delta}_{S_*^c}}-\lone{\hat{\Delta}_{S_*}}-2\lone{w^*_{S_*^c}}},
\end{align*}
which implies $\lone{\hat{\Delta}_{S_*}}\ge\frac{\varrho-1}{\varrho+1}\lone{\hat{\Delta}_{S_*^c}}-2\frac{\varrho}{\varrho+1}\lone{w^*_{S_*^c}}$. Letting $\varrho\to+\infty$, we have 
\begin{align}
\lone{\hat{\Delta}_{S_*}}\ge\lone{\hat{\Delta}_{S_*^c}}-2\lone{w^*_{S_*^c}}\label{proof:consistency-cone}
\end{align}
under $\mc{B}_5$. Consider the following two cases.
\begin{enumerate}
    \item $\lone{\hat{\Delta}_{S_*}}\leq\lone{w^*_{S_*^c}}$. Then by \eqref{proof:consistency-cone}
    \begin{align}
        \lone{\hat{\Delta}}\leq4\lone{w^*_{S_*^c}}.\label{proof:consistency-case1}
    \end{align}
    \item $\lone{\hat{\Delta}_{S_*}}>\lone{w^*_{S_*^c}}$. Recall the definition of $\phi^2(S_*)$. Under the events $\mc{B}_3$ and $\mc{B}_5$, if $4T_{n,q}^{(2)}(\delta)\sqrt{\abs{S_*}}\leq \phi(S_*)$, by \eqref{proof:consistency-cone} we have
\begin{align*}
\frac{\phi^2(S_*)(\lone{\hat{\Delta}_{S_*}}-\lone{w^*_{S_*^c}})^2}{4\abs{S_*}}
    \leq&~\parenBig{\frac{\phi(S_*)\lone{\hat{\Delta}_{S_*}}}{\sqrt{\abs{S_*}}}-2T_{n,q}^{(2)}(\delta)\lone{\hat{\Delta}_{S_*}}-2T_{n,q}^{(2)}(\delta)\lone{w^*_{S_*^c}}}^2\\
    \leq&~\parenBig{\ltwo{B^*\hat{\Delta}}-T_{n,q}^{(2)}(\delta)\lone{\hat{\Delta}}}^2\\
    \leq&~\frac{1}{J}\sum_{m=1}^J\parenBig{\ltwo{B^*\hat{\Delta}}-\ltwo{\parenBig{\hat{B}^{\brac{-m}}-B^*}\hat{\Delta}}}^2\\
    \leq&~\frac{1}{J}\sum_{m=1}^J\ltwo{\hat{B}^{\brac{-m}}\hat{\Delta}}^2.
\end{align*}
Under the events $\mc{B}_3,\mc{B}_4,\mc{B}_5$, by \eqref{proof:consistency-cone} and some basic calculation, if $T_{n,q}^{(2)}(\delta)\leq K_\beta$, we have
\begin{equation}\label{proof:consistency-RSC}
\begin{aligned}
    0\ge&~\dn\sum_{m=1}^J\sum_{i\in\mc{I}_m}\setBig{L\paren{Y_i,\<\hat{B}^{\brac{-m}}\paren{w^*+\hat{\Delta}},X_i\>}-L\paren{Y_i,\<\hat{B}^{\brac{-m}}w^*,X_i\>}}\\
    &~+\lambda_n^\prime\parenBig{\lone{w^*+\hat{\Delta}}-\lone{w^*}}\\
    \ge&~-\lambda_n^\prime\lone{\hat{\Delta}}+\lambda_n^\prime\parenBig{\lone{\hat{\Delta}_{S_*^c}}-\lone{\hat{\Delta}_{S_*}}-2\lone{w^*_{S_*^c}}}\\
    &~+\frac{1}{J}\sum_{m=1}^JC_1\ltwo{\hat{B}^{\brac{-m}}\hat{\Delta}}\parenBig{\ltwo{\hat{B}^{\brac{-m}}\hat{\Delta}}-C_2\sqrt{\frac{\log p}{n_1}}\lone{\hat{B}^{\brac{-m}}\hat{\Delta}}}\\
    \ge&~-\lambda_n^\prime \lone{\hat{\Delta}}+\lambda_n^\prime\parenBig{\lone{\hat{\Delta}_{S_*^c}}-\lone{\hat{\Delta}_{S_*}}-2\lone{w^*_{S_*^c}}}+\frac{1}{J}\sum_{m=1}^JC_1\ltwo{\hat{B}^{\brac{-m}}\hat{\Delta}}^2\\
    &~-C_3T_{n,q}^{(3)}(\delta)\sqrt{\frac{\log p}{n_1}}\lone{\hat{\Delta}}^2\\
    \ge&~-\lambda_n^\prime\lone{\hat{\Delta}}+\lambda_n^\prime\parenBig{\lone{\hat{\Delta}_{S_*^c}}-\lone{\hat{\Delta}_{S_*}}-2\lone{w^*_{S_*^c}}}-32C_3T_{n,q}^{(3)}(\delta)\sqrt{\frac{\log p}{n_1}}\lone{w_{S_*^c}^*}\\
    &~+(\frac{\phi^2(S_*)}{4\abs{S_*}}-4C_3T^{(3)}_{n,q}(\delta)\sqrt{\frac{\log p}{n_1}})(\lone{\hat{\Delta}_{S_*}}-\lone{w^*_{S_*^c}})^2\\
    \ge&~-2\lambda_n^\prime\absBig{\lone{\hat{\Delta}_{S_*}}-\lone{w^*_{S_*^c}}}-\parenBig{4\lambda_n^\prime+32C_3T_{n,q}^{(3)}(\delta)\sqrt{\frac{\log p}{n_1}}}\lone{w^*_{S_*^c}}\\
    &~+(\frac{\phi^2(S_*)}{4\abs{S_*}}-4C_3T^{(3)}_{n,q}(\delta)\sqrt{\frac{\log p}{n_1}})(\lone{\hat{\Delta}_{S_*}}-\lone{w^*_{S_*^c}})^2\\
    \ge&~-2\lambda_n^\prime\absBig{\lone{\hat{\Delta}_{S_*}}-\lone{w^*_{S_*^c}}}-\parenBig{4\lambda_n^\prime+32C_3T_{n,q}^{(3)}(\delta)\sqrt{\frac{\log p}{n_1}}}\lone{w^*_{S_*^c}}\\
    &~+\frac{\phi^2(S_*)}{8\abs{S_*}}(\lone{\hat{\Delta}_{S_*}}-\lone{w^*_{S_*^c}})^2
\end{aligned}
\end{equation}
provided that $4C_3T^{(3)}_{n,q}(\delta)\sqrt{\frac{\log p}{n_1}}\leq \frac{\phi^2(S_*)}{8\abs{S_*}}$. Note that the right-hand side of \eqref{proof:consistency-RSC} is a quadratic function of $\absBig{\lone{\hat{\Delta}_{S_*}}-\lone{w^*_{S_*^c}}}$, which implies \begin{align}
\lone{\hat{\Delta}_{S_*}}\leq&~\frac{2\lambda_n^\prime+\sqrt{(2\lambda_n^\prime)^2+4\lone{w^*_{S_*^c}}\parenBig{4\lambda_n^\prime+32C_3T_{n,q}^{(3)}(\delta)\sqrt{\frac{\log p}{n_1}}}\frac{\phi^2(S_*)}{8\abs{S_*}}}}{2\times\frac{\phi^2(S_*)}{8\abs{S_*}}}+\lone{w^*_{S_*^c}}\\
\leq&~\frac{16\abs{S_*}\lambda_n^\prime}{\phi^2(S_*)}+\sqrt{\frac{8\abs{S_*}}{\phi^2(S_*)}\parenBig{4\lambda_n^\prime+32C_3T_{n,q}^{(3)}(\delta)\sqrt{\frac{\log p}{n_1}}}\lone{w^*_{S_*^c}}}+\lone{w^*_{S_*^c}}.\label{proof:consistency-case2}
\end{align}
\end{enumerate}

Recall that under $\mc{B}_5$ we have $\lone{\hat{\Delta}_{S_*}}\ge\lone{\hat{\Delta}_{S_*^c}}-2\lone{w^*_{S_*^c}}$ for any $S_*\subseteq \brac{K}$. Combining \eqref{proof:consistency-case1} and \eqref{proof:consistency-case2}, under the events $\mc{B}_3,\mc{B}_4,\mc{B}_5$, if $T_{n,q}^{(2)}(\delta)\leq K_\beta$, $4T_{n,q}^{(2)}(\delta)\sqrt{\abs{S_*}}\leq \phi(S_*)$ and $4C_3T^{(3)}_{n,q}(\delta)\sqrt{\frac{\log p}{n_1}}\leq \frac{\phi^2(S_*)}{8\abs{S_*}}$, we have
\begin{align}
\lone{\hat{\Delta}}\leq\frac{32\abs{S_*}\lambda_n^\prime}{\phi^2(S_*)}+4\lone{w^*_{S_*^c}}+2\sqrt{\frac{8\abs{S_*}}{\phi^2(S_*)}\parenBig{4\lambda_n^\prime+32C_3T_{n,q}^{(3)}(\delta)\sqrt{\frac{\log p}{n_1}}}\lone{w^*_{S_*^c}}}.
\end{align}
Setting $t_1=\sqrt{\frac{2\log (pK)}{cn_1}}$, $t_2=\sqrt{\frac{\log (pK)}{cn_1}}$ and $\delta=\sqrt{\frac{\log (pK)}{n_{-1}}}$, then $T_{n,q}^{(2)}(\delta)\asymp T_{n,q}\lesssim\Psi_{n,q}$ and $T_{n,q}^{(3)}(\delta)\sqrt{\frac{\log p}{n}}\lesssim\Psi_{n,q}$. Consequently, we have
\begin{align*}
    \lambda_n^\prime=C\parenBig{(t_1+1)T_{n,q}^{(2)}(\delta)+t_2}\asymp T_{n,q},
\end{align*}
and then 
\begin{align*}
    \lone{\hat{\Delta}}\leq C_1\frac{\abs{S_*}}{\phi^2(S_*)}T_{n,q}+C_1\lone{w^*_{S_*^c}}+C_1\sqrt{\frac{\abs{S_*}}{\phi^2(S_*)}\Psi_{n,q}\lone{w^*_{S_*^c}}}
\end{align*}
uniformly holds for any $S_*\subseteq \brac{K}$ with $\abs{S_*}\leq c_2\phi^2(S_*)\Psi_{n,q}^{-1}$ with probability $1-C_2(pK)^{-1}-C_2e^{-cn}$ since $J$ is fixed. 
\end{proof}
%%%%%%%%%%%%%%%%%%%%%%%%%%%%%%%%%%%%%%%%%%%%%%%%%%%%%%%%
\subsection{Proof of \cref{HD-MA:minimax}}
\begin{proof} By \cref{VGBound} (Varshamov-Gilbert lemma), there exists a subset $\set{w^{(1)},\dots,w^{(M)}}\subseteq\mc{H}(K-1,s_w-1)$ such that $\lone{w^{(j_1)}-w^{(j_2)}}\ge (s_w-1)/2$ for any $j_1, j_2\in\brac{M}$, where $M$ is a positive integer that satisfies $\log M\ge\frac{s_w-1}{2}\log\frac{2(K-1)}{5(s_w-1)}$ and $\mc{H}(d,s)=\set{w\in\set{0,1}^d|~\lzero{w}=s}$. Denote 
\begin{align*}
    \tilde{w}^{(j)}\overset{\triangle}{=}\parenBig{\frac{4\delta_1}{s_w-1}w^{(j)\top},1-4\delta_1}^\top\in\mc{W}^K,
\end{align*}
where $\delta_1\leq 1/4$. Then $\lone{\tilde{w}^{(j)}}=1$ and $\lzero{\tilde{w}^{(j)}}\leq s_w$ for $j\in\brac{M}$. Denote $S^{(j)}=\text{supp}(\tilde{w}^{(j)})$. We set $\mc{A}_k=\set{(k-1)\lfloor\frac{p}{K}\rfloor,\dots,k\lfloor\frac{p}{K}\rfloor}, k=1,\dots,K$ and consider the following form of $\bm{\beta}^{(j)}$:
\begin{align*}
    \bm{\beta}^{(j)}=(\underbrace{z_1^{(j)},0,\dots,0}_{\mc{A}_1},\underbrace{z_2^{(j)},0,\dots,0}_{\mc{A}_2},\dots,\underbrace{z_K^{(j)},0,\dots,0}_{\mc{A}_K},0,\dots),
\end{align*}
where
\begin{align*}
    z_k^{(j)}=\begin{cases}
        2\sqrt{\delta_1}\delta_2,~&k\notin S^{(j)}\\
        2\sqrt{\delta_1}\delta_2\sqrt{\frac{1}{1-\frac{4\delta_1}{s_w-1}}},~&k\in S^{(j)}\backslash\set{K}\\
        \delta_2,~&k=K
    \end{cases}
\end{align*}
and $\delta_2>0$ for $j\in\brac{M}$ and $k\in\brac{K}$. Since $X_1\sim N(0,I_p)$, we have
\begin{align*}
    \beta^{*(j)}_{(k)}=(0,\dots,0,\underbrace{z_k^{(j)},0,\dots,0}_{\mc{A}_k},0,\dots,0)^\top.
\end{align*}
Let $B^{(j)}=(\beta^{*(j)}_{(1)},\beta^{*(j)}_{(2)},\dots,\beta^{*(j)}_{(K)})$. Note that for any $w\in\mc{W}^K$, by some simple calculation,
\begin{align*}
    &~\<\E_{\bm{\beta}^{(j)}}\frac{\partial L(Y_1,\<B^{(j)}\tilde{w}^{(j)},X_1\>)}{\partial \mu}B^{(j)\top}X_1,w-\tilde{w}^{(j)}\>\\
    =&~\<B^{(j)}\tilde{w}^{(j)}-\bm{\beta}^{(j)}_{\brac{p}},B^{(j)}(w-\tilde{w}^{(j)})\>\\
    =&~\sum_{k=1}^Kw_k(\tilde{w}_k^{(j)}-1)\abs{z_k^{(j)}}^2-\parenBig{\sum_{k=1}^K(\tilde{w}_k^{(j)})^2\abs{z_k^{(j)}}^2-\sum_{k=1}^K\tilde{w}_k^{(j)}\abs{z_k^{(j)}}^2}\\
    =&~0,
\end{align*}
which implies $\tilde{w}^{(j)}=\argmin_{w\in\mc{W}^K}\E_{\bm{\beta}^{(j)}}L(Y_1,\<B^{(j)}w,X_1\>)$ by optimality condition. Since $B^{(j)}$ is of full column rank, $\tilde{w}^{(j)}$ is the unique minimizer for each $j\in\brac{M}$. Thus, when
\begin{align}\label{proof:lower-bound}
   \delta_2=\max_{k\in\brac{K}}\ltwo{\beta_{(k)}^{(j)}}\leq C\text{\ \ \ and\ \ \ }\phi^2(S^{(j)})\ge4\delta_1\delta_2^2\ge\phi_n^2,
\end{align}
we have $(\tilde{w}^{(j)},B^{(j)})\in\Theta$ for each $j\in\brac{M}$.

Note that for any $j_1,j_2\in\brac{M}, j_1\neq j_2$
\begin{align*}
    \lone{\tilde{w}^{(j_1)}-\tilde{w}^{(j_2)}}=\frac{4\delta_1}{s_w-1}\parenBig{\lone{w^{(j_1)}-w^{(j_2)}}}\ge2\delta_1
\end{align*}
and the KL divergence
{\small\begin{align*}
    D_{KL}^{(n)}\parenBig{f_{\bm{\beta}^{(j_2)}}\Big|\Big|f_{\bm{\beta}^{(j_1)}}}\overset{\triangle}{=}&~\E_{\bm{\beta}^{(j_2)}}\bracBig{-\log\parenBig{\frac{\prod_{i=1}^nf_{Y,X}(Y_i,X_i|\bm{\beta}^{(j_1)})}{\prod_{i=1}^nf_{Y,X}(Y_i,X_i|\bm{\beta}^{(j_2)})}}}\\
    =&~\E_{\bm{\beta}^{(j_2)}}\parenBig{\sum_{i=1}^n \setBig{L(Y_1,\<\bm{\beta}^{(j_1)}_{\brac{p}},X_i\>)-L(Y_i,\<\bm{\beta}^{(j_2)}_{\brac{p}},X_i\>)}}\\
    \leq&~\E_{\bm{\beta}^{(j_2)}}\parenBig{\sum_{i=1}^n\setBig{\frac{\partial L(Y_i,\<\bm{\beta}^{(j_2)}_{\brac{p}},X_i\>)}{\partial\mu}\<\bm{\beta}^{(j_1)}_{\brac{p}}-\bm{\beta}^{(j_2)}_{\brac{p}},X_i\>+\<\bm{\beta}^{(j_1)}_{\brac{p}}-\bm{\beta}^{(j_2)}_{\brac{p}},X_i\>^2}}\\
    =&~n\E\<\bm{\beta}^{(j_1)}_{\brac{p}}-\bm{\beta}^{(j_2)}_{\brac{p}},X_i\>^2\\
    \leq&~n\ltwo{\bm{\beta}^{(j_1)}_{\brac{p}}-\bm{\beta}^{(j_2)}_{\brac{p}}}^2\\
    \leq &~2n(s_w-1)\times 4\delta_1\delta^2_2\parenBig{1-\sqrt{\frac{1}{1-\frac{4\delta_1}{s_w-1}}}}^2\\
    \leq&~\frac{Cn\delta_1^3\delta_2^2}{s_w-1},
\end{align*}}
where we use $\E_{\bm{\beta}^{(j_2)}}\bracBig{\frac{\partial L(Y_i,\<\bm{\beta}^{(j_2)}_{\brac{p}},X_i\>)}{\partial\mu}X_i}=0$ and the inequality $(1-(1-x)^{-1/2})^2\leq Cx^2, x\in[0,c)$ for some constant $c, C>0$. 
In order to apply Fano's inequality \citep[see, for example, Proposition 15.12 in][]{wainwright2019high}, we need to find $\delta_1, \delta_2>0$ such that the inequalities in \eqref{proof:lower-bound} and the following inequality
\begin{align*}
    1-\frac{\log 2+\frac{Cn\delta_1^3\delta_2^2}{s_w-1}}{\frac{s_w-1}{2}\log\frac{2(K-1)}{5(s_w-1)}}\ge \frac{1}{2}
\end{align*}
hold simultaneously. Choosing $\delta_1^2\asymp \frac{s_w^2\log(K/s_w)}{n\phi_n^2}$ and $\delta_2^2\asymp\frac{\phi_n^2}{\delta_1}$, when $\phi_n^2\lesssim\parenBig{s_w^2\frac{\log(K/s_w)}{n}}^{1/3}$, by Fano's inequality, we have
\begin{align*}
\inf_{\hat{w}}\sup_{Q\in\mc{P}_{Y,X}^{1:n}}\P_{(Y_i,X_i)\sim Q}\parenBig{\lone{\hat{w}-w^*}\ge Cs_w\sqrt{\frac{\log(K/s_w)}{n\phi_n^2}}}\ge\frac{1}{2}.
    \end{align*}

\end{proof}
%%%%%%%%%%%%%%%%%%%%%%%%%%%%%%%%%%%%%%%%%%%%%%%%%%%%%%%%%%

\subsection{Proof of \cref{HD-MA:asymptotical-optimality}}
\begin{proof} By \cref{HD-MA:beta-deviation},
\begin{align*}
    \sup_{k\in\brac{K}}\ltwo{\hat{\beta}_{(k)}-\beta_{(k)}^*}=O_p(T_{n,q})
\end{align*}
and similarly we have
\begin{align*}
    \sup_{m\in\brac{J}}\sup_{k\in\brac{K}}\ltwo{\hat{\beta}^{\brac{-m}}_{(k)}-\beta_{(k)}^*}=O_p(T_{n,q}),
\end{align*}
since $J$ is fixed.
Since $\hat{w}$ is the minimizer of $CV(w)$, by \cref{AOPrate}, it suffices to show that
\begin{align*}
    \xi_n^{-1}\sup_{w\in\mc{W}^K}\abs{R(w)-R^*(w)}=O_p(\xi_n^{-1}T_{n,q})
\end{align*}
and
\begin{align*}
    \xi_n^{-1}\sup_{w\in\mc{W}^K}\abs{R^*(w)-CV(w)/n}=O_p(\xi_n^{-1}T_{n,q}).
\end{align*}

Note that
\begin{align*}
    &~\xi_n^{-1}\sup_{w\in\mc{W}^K}\abs{R(w)-R^*(w)}\\
    =&~\xi_n^{-1}\sup_{w\in\mc{W}^K}\absBig{\E\brac{L\paren{Y_{n+1},\<\hat{B}w,X_{n+1}\>}|\set{(Y_i,X_i)}_{i=1}^n}-\E L\paren{Y_{n+1},\<B^*w,X_{n+1}\>}}
\end{align*}
and
\begin{align*}
    &~\xi_n^{-1}\sup_{w\in\mc{W}^K}\abs{R^*(w)-CV(w)/n}\\
\leq&~\xi_n^{-1}\sup_{w\in\mc{W}^K}\absBig{\dn\sum_{m=1}^J\sum_{i\in\mc{I}_m}\setBig{L\paren{Y_{i},\<\hat{B}^{\brac{-m}}w,X_{i}\>}-L\paren{Y_{i},\<B^*w,X_{i}\>}}}\\
&~+\xi_n^{-1}\sup_{w\in\mc{W}^K}\absBig{\dn\sumn L\paren{Y_{i},\<B^*w,X_{i}\>}-\E L\paren{Y_{i},\<B^*w,X_{i}\>}}.
\end{align*}
Thus, we just need to show that
\begin{align*}
  \xi_n^{-1}\sup_{w\in\mc{W}^K}\absBig{\E\brac{L\paren{Y_{n+1},\<\hat{B}w,X_{n+1}\>}|\set{Y_i,X_i}_{i=1}^n}-\E L\paren{Y_{n+1},\<B^*w,X_{n+1}\>}}=&~O_p(\xi_n^{-1}T_{n,q})\\
    \xi_n^{-1}\sup_{w\in\mc{W}^K}\absBig{\dn\sum_{m=1}^J\sum_{i\in\mc{I}_m}\setBig{L\paren{Y_{i},\<\hat{B}^{\brac{-m}}w,X_{i}\>}-L\paren{Y_{i},\<B^*w,X_{i}\>}}}=&~O_p(\xi_n^{-1}T_{n,q})\\
    \xi_n^{-1}\sup_{w\in\mc{W}^K}\absBig{\dn\sumn L\paren{Y_{i},\<B^*w,X_{i}\>}-\E L\paren{Y_{i},\<B^*w,X_{i}\>}}=&~O_p(\xi_n^{-1}T_{n,q}).
\end{align*}
\begin{enumerate}
    \item Suppose $L\paren{y,\mu}=\frac{1}{2}\paren{y-\mu}^2$. Conditioning on $\set{(Y_i,X_i)}_{i=1}^n$, for any $w\in\mc{W}^K$, we have, 
\begin{align*}
    \norm{\frac{\<\hat{B}w-B^*w, X_{n+1}Y_{n+1}\>}{\ltwo{\hat{B}w-B^*w}}}_{\psi_1}\leq&~ K_xK_y\\
    \norm{\frac{\<\hat{B}w-B^*w,X_{n+1}\>\<\hat{B}w+B^*w,X_{n+1}\>}{\ltwo{\hat{B}w-B^*w}\ltwo{\hat{B}w+B^*w}}}_{\psi_1}\leq&~ K_x^2.
\end{align*} By Proposition 2.7.1 in \cite{vershynin2018high} we have
\begin{equation}
\label{proof:AOP-case1-1}
\begin{aligned}
    &~ \sup_{w\in\mc{W}^K}\absBig{\E\brac{L\paren{Y_{n+1},\<\hat{B}w,X_{n+1}\>}|\set{(Y_i,X_i)}_{i=1}^n}-\E L\paren{Y_{n+1},\<B^*w,X_{n+1}\>}}\\
    \leq&~\sup_{w\in\mc{W}^K}\E\bracBig{\abs{\frac{\<\hat{B}w-B^*w, X_{n+1}Y_{n+1}\>}{\ltwo{\hat{B}w-B^*w}}}\Big|\set{(Y_i,X_i)}_{i=1}^n}\sup_{w\in\mc{W}^K}\ltwo{\hat{B}w-B^*w}\\
    &~+\sup_{w\in\mc{W}^K}\frac{1}{2}\E\bracBig{\abs{\frac{\<\hat{B}w-B^*w,X_{n+1}\>\<\hat{B}w+B^*w,X_{n+1}\>}{\ltwo{\hat{B}w-B^*w}\ltwo{\hat{B}w+B^*w}}}\Big|\set{(Y_i,X_i)}_{i=1}^n}\\
    &~\ \ \ \ ~\times\sup_{w\in\mc{W}^K}\ltwo{\hat{B}w-B^*w}\ltwo{\hat{B}w+B^*w}\\
    \leq&~CK_xK_y\sup_{k\in\brac{K}}\ltwo{\hat{\beta}_{(k)}-\beta^*_{(k)}}+CK_x^2\parenBig{\sup_{k\in\brac{K}}\ltwo{\hat{\beta}_{(k)}-\beta^*_{(k)}}}\parenBig{\sup_{k\in\brac{K}}\ltwo{\hat{\beta}_{(k)}+\beta^*_{(k)}}}\\
    =&~O_p(T_{n,q}).
\end{aligned}
\end{equation}
For any $m\in\brac{J}$ we have
\begin{equation}
\label{proof:AOP-case1-2}
\begin{aligned}
    &~\sup_{w\in\mc{W}^K}\absBig{\frac{1}{n_1}\sum_{i\in\mc
    {I}_m}L\paren{Y_{i},\<\hat{B}^{\brac{-m}}w,X_{i}\>}-L\paren{Y_{i},\<B^*w,X_{i}\>}}\\
    \leq&~\sup_{w\in\mc{W}^K}\absBig{\frac{1}{n_1}\sum_{i\in\mc
    {I}_m}\<\hat{B}^{\brac{-m}}w-B^*w, X_{i}Y_{i}\>}\\
    &~+\sup_{w\in\mc{W}^K}\frac{1}{2}\absBig{\frac{1}{n_1}\sum_{i\in\mc
    {I}_m}\<\hat{B}^{\brac{-m}}w-B^*w,X_{i}\>\<\hat{B}^{\brac{-m}}w+B^*w,X_{i}\>}\\
    \leq&~\sup_{k\in\brac{K}}\absBig{\frac{1}{n_1}\sum_{i\in\mc
    {I}_m}\<\hat{\beta}^{\brac{-m}}_{(k)}-\beta^*_{(k)}, X_{i}\>Y_i}\\
    &~+\frac{1}{2}\sup_{k_1,k_2\in\brac{K}}\absBig{\frac{1}{n_1}\sum_{i\in\mc
    {I}_m}\<\hat{\beta}^{\brac{-m}}_{(k_1)}-\beta^*_{(k_1)},X_{i}\>\<\hat{\beta}^{\brac{-m}}_{(k_2)}+\beta^*_{(k_2)},X_{i}\>}.
\end{aligned}
\end{equation}
By \cref{indnorm} and Lemma 2.7.7 in \cite{vershynin2018high}, for any $m\in\brac{J}$, $i\in\mc{I}_m$ and $k\in\brac{K}$, we have
\begin{align*}
    \norm{\frac{\<\hat{\beta}^{\brac{-m}}_{(k)}-\beta^*_{(k)}, X_{i}\>Y_i}{\ltwo{\hat{\beta}_{(k)}^{\brac{-m}}-\beta^*_{(k)}}}}_{\psi_1}\leq K_xK_y,
\end{align*}
and then we have
\begin{align*}
   \sup_{k\in\brac{K}}\E\absBig{\frac{\<\hat{\beta}^{\brac{-m}}_{(k)}-\beta^*_{(k)}, X_{i}\>Y_i}{\ltwo{\hat{\beta}^{\brac{-m}}_{(k)}-\beta^*_{(k)}}}}\leq CK_xK_y .
\end{align*}
By \cref{centering,supsumsubexp} and Jensen's inequality we have
\begin{equation}
\label{proof:AOP-case1-3}
\begin{aligned}
    &~\sup_{k\in\brac{K}}\absBig{\frac{1}{n_1}\sum_{i\in\mc
    {I}_m}\<\hat{\beta}^{\brac{-m}}_{(k)}-\beta^*_{(k)}, X_{i}\>Y_i}\\
    \leq&~\sup_{k\in\brac{K}}\absBig{\frac{1}{n_1}\sum_{i\in\mc
    {I}_m}\frac{\<\hat{\beta}^{\brac{-m}}_{(k)}-\beta^*_{(k)}, X_{i}\>Y_i}{\ltwo{\hat{\beta}^{\brac{-m}}_{(k)}-\beta^*_{(k)}}}-\E\parenBig{\frac{\<\hat{\beta}^{\brac{-m}}_{(k)}-\beta^*_{(k)}, X_{i}\>Y_i}{\ltwo{\hat{\beta}^{\brac{-m}}_{(k)}-\beta^*_{(k)}}}}}\sup_{k\in\brac{K}}\ltwo{\hat{\beta}^{\brac{-m}}_{(k)}-\beta^*_{(k)}}\\
    &~+\sup_{k\in\brac{K}}\E\absBig{\frac{\<\hat{\beta}^{\brac{-m}}_{(k)}-\beta^*_{(k)}, X_{i}\>Y_i}{\ltwo{\hat{\beta}^{\brac{-m}}_{(k)}-\beta^*_{(k)}}}}\sup_{k\in\brac{K}}\ltwo{\hat{\beta}^{\brac{-m}}_{(k)}-\beta^*_{(k)}}\\
    =&~\parenBig{O_p(\sqrt{\frac{\log K}{n}})+CK_xK_y}O_p(T_{n,q})\\
    =&~O_p(T_{n,q}),
\end{aligned}
\end{equation}
since $\sqrt{\frac{\log K}{n}}\to 0$ as $n\to\infty$. 
By \cref{indnorm} and Lemma 2.7.7 of \cite{vershynin2018high}, for any $m\in\brac{J}$, $i\in\mc{I}_m$ and $k_1,k_2\in\brac{K}$ we have
\begin{align*}
    \norm{\frac{\<\hat{\beta}^{\brac{-m}}_{(k_1)}-\beta^*_{(k_1)},X_{i}\>\<\hat{\beta}^{\brac{-m}}_{(k_2)}+\beta^*_{(k_2)},X_{i}\>}{\ltwo{\hat{\beta}_{(k_1)}^{\brac{-m}}-\beta^*_{(k_1)}}\ltwo{\hat{\beta}_{(k_2)}^{\brac{-m}}+\beta^*_{(k_2)}}}}_{\psi_1}\leq K_x^2,
\end{align*}
Similarly to the previous proof, we have
\begin{equation}\label{proof:AOP-case1-4}
\begin{aligned}
    \frac{1}{2}\sup_{k_1,k_2\in\brac{K}}\absBig{\frac{1}{n_1}\sum_{i\in\mc
    {I}_m}\<\hat{\beta}^{\brac{-m}}_{(k_1)}-\beta^*_{(k_1)},X_{i}\>\<\hat{\beta}^{\brac{-m}}_{(k_2)}+\beta^*_{(k_2)},X_{i}\>}=O_p(T_{n,q}).
\end{aligned}
\end{equation}
Combining \eqref{proof:AOP-case1-2}, \eqref{proof:AOP-case1-3} and \eqref{proof:AOP-case1-4}, since $J$ is fixed, by Jensen's inequality, we have
\begin{equation}\label{proof:AOP-case1-5}
\begin{aligned}
     &~\sup_{w\in\mc{W}^K}\absBig{\dn\sum_{m\in\brac{J}}\sum_{i\in\mc{I}_m}L\paren{Y_{i},\<\hat{B}^{\brac{-m}}w,X_{i}\>}-L\paren{Y_{i},\<B^*w,X_{i}\>}}\\
     \leq &~\frac{1}{J}\sum_{m\in\brac{J}}\sup_{w\in\mc{W}^K}\absBig{\frac{1}{n_1}\sum_{i\in\mc
    {I}_m}L\paren{Y_{i},\<\hat{B}^{\brac{-m}}w,X_{i}\>}-L\paren{Y_{i},\<B^*w,X_{i}\>}}\\
     =&~O_p(T_{n,q}).
\end{aligned}
\end{equation}

By symmetrization, Holder's inequality and \cref{supsumsubexp}, we have
{\small\begin{equation}\label{proof:AOP-case1-6}
\begin{aligned}
    &~\E\sup_{w\in\mc{W}^K}\absBig{\dn\sumn L\paren{Y_{i},\<B^*w,X_{i}\>}-\E L\paren{Y_{i},\<B^*w,X_{i}\>}}\\
    \leq&~2\E\sup_{w\in\mc{W}^K}\absBig{\dn\sumn \eps_iL\paren{Y_{i},\<B^*w,X_{i}\>}}\\
    \leq&~\E\absBig{\dn\sumn \eps_iY_i^2}+2\E\sup_{w\in\mc{W}^K}\absBig{\dn\sumn\eps_iY_i\<w,B^{*\top}X_i\>}+\E\sup_{w\in\mc{W}^K}\absBig{\dn\sumn\eps_i\<w^{\otimes 2}
    ,(B^{*\top}X_i)^{\otimes 2}\>}\\
    \leq&~\parenBig{\E\bracBig{\dn\sumn \eps_iY_i^2}^2}^{\frac{1}{2}}+2\E\sup_{k\in\brac{K}}\absBig{\dn\sumn\eps_iY_i\<\beta_{k}^*,X_i\>}+\E\sup_{k_1,k_2\in\brac{K}}\absBig{\dn\sumn\eps_i\<\beta_{k_1}^*,X_i\>\<\beta_{k_2}^*,X_i\>}\\
    \leq&~\parenBig{\dn\E\bracBig{ \eps_1Y_1^2}^2}^{\frac{1}{2}}+2\E\sup_{k\in\brac{K}}\absBig{\dn\sumn\eps_iY_i\<\beta_{k}^*,X_i\>}+\E\sup_{k_1,k_2\in\brac{K}}\absBig{\dn\sumn\eps_i\<\beta_{k_1}^*,X_i\>\<\beta_{k_2}^*,X_i\>}\\
    =&~O(\sqrt{\frac{1}{n}})+O(\sqrt{\frac{\log K}{n}})+O(\sqrt{\frac{\log K}{n}})\\
    =&~O(T_{n,q}),
\end{aligned}
\end{equation}}
where $\set{\eps_i}_{i=1}^n$ is a sequence of independent Rademacher variables and is independent of $\set{(Y_i,X_i)}_{i=1}^n$.
By \eqref{proof:AOP-case1-1}, \eqref{proof:AOP-case1-5} and \eqref{proof:AOP-case1-6}, the result holds under the first condition in \cref{HD-MA:lossfunc}.
\item Suppose $\abs{L\paren{y,\mu_1}-L\paren{y,\mu_2}}\leq L_u^\prime\abs{\mu_1-\mu_2}$ for any $\mu_1,\mu_2\in\R$, $y\in\mc{T}$. Conditioning on $\set{(Y_i,X_i)}_{i=1}^n$, for any $w\in\mc{W}^K$ we have
\begin{align*}
    \norm{\frac{\abs{\<\hat{B}w-B^*w,X_{n+1}\>}}{\ltwo{\hat{B}w-B^*w}}}_{\psi_2}\leq K_x,
\end{align*}
which means $\sup_{w\in\mc{W}^K}\E\parenBig{\frac{\abs{\<\hat{B}w-B^*w,X_{n+1}\>}}{\ltwo{\hat{B}w-B^*w}}|\set{(Y_i,X_i)}_{i=1}^n}\leq CK_x$ for some absolute constant $C>0$ by Proposition 2.5.2 of \cite{vershynin2018high}. Then we have
\begin{equation}\label{proof:AOP-case2-1}
\begin{aligned}
    &~\sup_{w\in\mc{W}}\absBig{\E\brac{L\paren{Y_{n+1},\<\hat{B}w,X_{n+1}\>}|\set{(Y_i,X_i)}_{i=1}^n}-\E L\paren{Y_{n+1},\<B^*w,X_{n+1}\>}}\\
    \leq&~L_u^\prime\sup_{w\in\mc{W}^K}\absBig{\E\brac{\<\hat{B}w-B^*w,X_{n+1}\>|\set{(Y_i,X_i)}_{i=1}^n}}\\
    \leq&~L_u^\prime\sup_{w\in\mc{W}}\E\bracBig{\abs{\frac{\<\hat{B}w-B^*w, X_{n+1}\>}{\ltwo{\hat{B}w-B^*w}}}\Big|\set{(Y_i,X_i)}_{i=1}^n}\sup_{w\in\mc{W}^K}\ltwo{\hat{B}w-B^*w}\\
    \leq&~CL_u^\prime K_x\max_{k\in\brac{K}}\ltwo{\hat{\beta}_{(k)}-\beta^*_{(k)}}\\
    =&~O_p(T_{n,q}).
\end{aligned}
\end{equation}

By \cref{indnorm}, for any $m\in\brac{J}$, $i\in\mc{I}_m$ and $k\in\brac{K}$, we have
\begin{align*}
    \norm{\frac{\<\hat{\beta}^{\brac{-m}}_{(k)}-\beta^*_{(k)}, X_{i}\>}{\ltwo{\hat{\beta}_{(k)}^{\brac{-m}}-\beta^*_{(k)}}}}_{\psi_2}\leq K_x
\end{align*}
and then we have
\begin{align*}
   \sup_{k\in\brac{K}}\E\absBig{\frac{\<\hat{\beta}^{\brac{-m}}_{(k)}-\beta^*_{(k)}, X_{i}\>}{\ltwo{\hat{\beta}^{\brac{-m}}_{(k)}-\beta^*_{(k)}}}}\leq CK_x .
\end{align*}
 By \cref{centering,supsumsubexp} and Jensen's inequality, for any $m\in\brac{J}$, we have
\begin{align*}
    &~\sup_{w\in\mc{W}^K}\absBig{\frac{1}{n_1}\sum_{i\in\mc{I}_m}L\paren{Y_{i},\<\hat{B}^{\brac{-m}}w,X_{i}\>}- L\paren{Y_{i},\<B^*w,X_{i}\>}}\\
    \leq&~L_u^\prime\sup_{w\in\mc{W}^K}\absBig{\frac{1}{n_1}\sum_{i\in\mc{I}_m}\<\hat{B}w-B^*w,X_{i}\>}\\
    \leq&~L_u^\prime\sup_{k\in\brac{K}}\absBig{\frac{1}{n_1}\sum_{i\in\mc{I}_m}\<\hat{\beta}^{\brac{-m}}_{(k)}-\beta^*_{(k)},X_{i}\>}\\
    \leq&~L_u^\prime\sup_{k\in\brac{K}}\absBig{\frac{1}{n_1}\sum_{i\in\mc{I}_m}\frac{\<\hat{\beta}^{\brac{-m}}_{(k)}-\beta^*_{(k)}, X_{i}\>}{\ltwo{\hat{\beta}^{\brac{-m}}_{(k)}-\beta^*_{(k)}}}-\E\frac{\<\hat{\beta}^{\brac{-m}}_{(k)}-\beta^*_{(k)}, X_{i}\>}{\ltwo{\hat{\beta}^{\brac{-m}}_{(k)}-\beta^*_{(k)}}}}\sup_{k\in\brac{K}}\ltwo{\hat{\beta}^{\brac{-m}}_{(k)}-\beta^*_{(k)}}\\
    &~+L_u^\prime\sup_{k\in\brac{K}}\E\absBig{\frac{\<\hat{\beta}^{\brac{-m}}_{(k)}-\beta^*_{(k)}, X_{i}\>}{\ltwo{\hat{\beta}^{\brac{-m}}_{(k)}-\beta^*_{(k)}}}}\sup_{k\in\brac{K}}\ltwo{\hat{\beta}^{\brac{-m}}_{(k)}-\beta^*_{(k)}}\\
    \leq&~L_u^\prime\parenBig{O_p(\sqrt{\frac{\log K}{n}})+CK_x}O_p(T_{n,q})\\
    =&~O_p(T_{n,q}),
\end{align*}
which implies
\begin{equation}\label{proof:AOP-case2-2}
\begin{aligned}
\sup_{w\in\mc{W}^K}\absBig{\dn\sum_{m=1}^J\sum_{i\in\mc{I}_m}\setBig{L\paren{Y_{i},\<\hat{B}^{\brac{-m}}w,X_{i}\>}-L\paren{Y_{i},\<B^*w,X_{i}\>}}}=O_p\paren{T_{n,q}}
\end{aligned}
\end{equation}
since $J$ is fixed.
By symmetrization, contraction inequality from Theorem 4.12 of \cite{ledoux1991probability} and \cref{supsumsubexp}, we have
\begin{equation}\label{proof:AOP-case2-3}
\begin{aligned}
    &~\E\sup_{w\in\mc{W}^K}\absBig{\dn\sumn L\paren{Y_{i},\<B^*w,X_{i}\>}-\E L\paren{Y_{i},\<B^*w,X_{i}\>}}\\
    \leq&~2\E\sup_{w\in\mc{W}^K}\absBig{\dn\sumn \eps_iL\paren{Y_{i},\<B^*w,X_{i}\>}}\\
    \leq&~4L_u^\prime\E\sup_{w\in\mc{W}^K}\absBig{\dn\sumn \eps_i\<B^*w,X_i\>}\\
    \leq&~4L_u^\prime\E\sup_{k\in\brac{K}}\absBig{\dn\sumn \eps_i\<\beta_{(k)}^*,X_i\>}\\
    =&~O\paren{\sqrt{\frac{\log K}{n}}}\\
    =&~O(T_{n,q}),
\end{aligned}
\end{equation}
where $\set{\eps_i}_{i=1}^n$ is a sequence of independent Rademacher variables and is independent of $\set{(Y_i,X_i)}_{i=1}^n$. Thus, the result follows from \eqref{proof:AOP-case2-1}, \eqref{proof:AOP-case2-2} and \eqref{proof:AOP-case2-3}.
\end{enumerate}

\end{proof}
%%%%%%%%%%%%%%%%%%%%%%%%%%%%%%%%%%%%%%%%%%%%%%%%%%%%%%%%%%%%%%%%%

\subsection{Proof of \cref{weights2one}}
\begin{proof} The proof of this theorem is similar to the proof of \cref{MA}. We need to use the notation introduced in the proof of \cref{MA}. Denote $\mc{W}^*:=\set{w\in\mc{W}^K|\sum_{k=K^*}^{K_{ne}}w_k=1}$. Then we have
\begin{align}
    B^*w_0=\beta^*\label{proof:correctMA-w0}
\end{align}
for any $w_0\in\mc{W}^*$.  To show that the following event  $\mc{B}_{6,m}$ 
{\footnotesize\begin{align*}
    &~\Big\{\frac{1}{n_1}\sum_{i\in\mc{I}_m}\setBig{L\paren{Y_i,\<\hat{B}^{\brac{-m}}\paren{w_0+\Delta},X_i\>}-L\paren{Y_i,\<\hat{B}^{\brac{-m}}w_0,X_i\>}-\frac{\partial L\paren{Y_i,\<\hat{B}^{\brac{-m}}w_0,X_i\>}}{\partial \mu}\<\hat{B}^{\brac{-m}}\Delta,X_i\>}\\
    \ge&~C_1\ltwo{\hat{B}^{\brac{-m}}\Delta}\parenBig{\ltwo{\hat{B}^{\brac{-m}}\Delta}-C_2\sqrt{\frac{\log p}{n_1}}\lone{\hat{B}^{\brac{-m}}\Delta}}\text{ for any }\Delta\in\setBig{\Delta\in\R^K\Big|\lone{\Delta}\leq 2}\text{ and }w_0\in\mc{W}^*\Big\}
\end{align*}}holds with high probability for each $m\in\brac{J}$, we use \cref{lemma:RSC}.
\begin{enumerate}
    \item Assume that $\inf_{\mu\in\R}\inf_{y\in\R}\frac{\partial^2L(y,\mu)}{\partial\mu^2}\ge\eta_l>0$. By the second part of \cref{lemma:RSC} with $R_1=4K_\beta$, if $n_1\ge c_1$ and $T_{n,q}^{(2)}(\delta)\leq K_\beta$, we have
    \begin{align*}
        \E\parenBig{\P\paren{\mc{B}_{6,m}^c|\hat{B}^{\brac{-m}}}\Big|\mc{B}_{3,m}}\leq 2\exp\paren{-c_2n_1},
    \end{align*}
    where $\mc{B}_{3,m}$ is defined in \eqref{proof:weight-B3m}.
    
\item Assume that $(1\vee\sqrt{\log p})T_{n,q}^{(1)}(\delta)\leq 1$. Under $\mc{B}_{3,m}$, by \eqref{proof:correctMA-w0} for any $w_0\in\mc{W}^*$ we have
\begin{align*}
        \lone{\hat{B}^{\brac{-m}}w_0-\beta^*}\leq T_{n,q}^{(1)}(\delta)\leq\frac{1}{1\vee\sqrt{\log p}}.
    \end{align*}
    By the first part of \cref{lemma:RSC} with $\beta_0=\beta^*, R_1=4K_\beta, R_2=1$, if $n_1\ge c_1$ and $T_{n,q}^{(2)}(\delta)\leq K_\beta$, we can obtain that
    \begin{align*}
        \E\parenBig{\P\paren{\mc{B}_{6,m}^c|\hat{B}^{\brac{-m}}}\Big|\mc{B}_{3,m}}\leq 2\exp\paren{-c_2n_1}.
    \end{align*}

\end{enumerate}
Thus, we have
    \begin{align*}
        &~\P\paren{\mc{B}_{6,m}^c}=\E\parenBig{\P\paren{\mc{B}_{6,m}^c|\hat{B}^{\brac{-m}}}}\leq \E\parenBig{\P\paren{\mc{B}_{6,m}^c|\hat{B}^{\brac{-m}}}\Big|\mc{B}_{3,m}}+\P\parenBig{\mc{B}_3^c}\\
        \leq&~ 4JK\exp\paren{-n_{-1}\delta^2}+2\exp\paren{-c_2n_1}.
    \end{align*}
Consequently, the following event $\mc{B}_6\supseteq\cap_{m\in\brac{J}}\mc{B}_{6,m}$
{\footnotesize\begin{align*}
    &~\Big\{\dn\sum_{m=1}^J\sum_{i\in\mc{I}_m}\setBig{L\paren{Y_i,\<\hat{B}^{\brac{-m}}\paren{w_0+\Delta},X_i\>}-L\paren{Y_i,\<\hat{B}^{\brac{-m}}w_0,X_i\>}-\frac{\partial L\paren{Y_i,\<\hat{B}^{\brac{-m}}w_0,X_i\>}}{\partial \mu}\<\hat{B}^{\brac{-m}}\Delta,X_i\>}\\
    \ge&~\frac{1}{J}\sum_{m=1}^JC_1\ltwo{\hat{B}^{\brac{-m}}\Delta}\parenBig{\ltwo{\hat{B}^{\brac{-m}}\Delta}-C_2\sqrt{\frac{\log p}{n_1}}\lone{\hat{B}^{\brac{-m}}\Delta}}\text{ for any }\Delta\in\setBig{\Delta\in\R^K\Big|\lone{\Delta}\leq 2}\text{ and }w_0\in\mc{W}^*\Big\}
\end{align*}}
 holds with probability at least $1-4J^2K\exp\paren{-n_{-1}\delta^2}-2J\exp\paren{-c_2n_1}$.

Denote
\begin{align*}
    \hat{w}_0=(0,\dots,0,\sum_{k\neq K^*+1,\dots,K_{ne}}\hat{w}_k,\hat{w}_{K^*+1},\dots,\hat{w}_{K_{ne}},0,\dots,0)^\top\in\mc{W}^*.
\end{align*}
Since $B^*\hat{w}_0=\beta^*$ is a deterministic quantity and $\linf{\E\parenBig{\frac{\partial L\paren{Y_i,\<\beta^*,X_i\>}}{\partial\mu}B^{*\top}X_i}}=0$, similar to the proof of \cref{MA}, the event 
\begin{align*}
    \mc{B}_5=\setBig{\lambda_n^\prime\ge\linf{\dn\sum_{m=1}^J\sum_{i\in \mc{I}_m}\frac{\partial L\paren{Y_i,\<\hat{B}^{\brac{-m}}w_0,X_i\>}}{\partial \mu}\hat{B}^{\brac{-m}\top}X_i}}
\end{align*} holds with probability at least $1-2JK^2\exp(-cn_1t_1^2)-2JK\exp(-cn_1t_2^2)-4JK\exp\paren{-n_{-1}\delta^2}$, where $\lambda_n^\prime=C\parenBig{(t_1+1)T_{n,q}^{(2)}(\delta)+t_2}$ and $t_1,t_2\in (0,1)$.

Denote $\hat{\Delta}=\hat{w}-\hat{w}_0$ and $\breve{\Delta}=(\hat{w}_1,\dots,\hat{w}_{K^*-1},\sum_{k=K^*}^{K_{ne}}\hat{w}_k-1,\hat{w}_{K_{ne}+1},\dots,\hat{w}_K)^\top\in\mc{W}^{K_0}-\tilde{w}^*$. Then we have $\hat{\Delta}_k=0$ for $k=K^*+1,\dots,K_{ne}$ and 
\begin{align*}
    \lone{\hat{\Delta}}=&~\lone{\breve{{\Delta}}}=2\sum_{k\neq K^*,\dots,K_{ne}}\hat{w}_k=2\abs{1-\sum_{k\in \set{K^*,\dots,K_{ne}}}\hat{w}_k},\\
\ltwo{\hat{\Delta}}=&~\ltwo{\breve{{\Delta}}}\ge\abs{\breve{{\Delta}}_{K^*}}=\abs{1-\sum_{k\in\set{ K^*,\dots,K_{ne}}}\hat{w}_k}=\frac{\lone{\hat{\Delta}}}{2}
\end{align*}
Recall the definition of $\tilde{\phi}^2$. Under the event $\mc{B}_3$, if $\tilde{\phi}\ge 4T_{n,q}^{(2)}(\delta)$, we have
\begin{align*}
&~\frac{\tilde{\phi}^2}{4}\abs{\breve{{\Delta}}_{K^*}}^2\\
\leq&~\parenBig{\tilde{\phi}\abs{\breve{{\Delta}}_{K^*}}-2T_{n,q}^{(2)}(\delta)\abs{\breve{{\Delta}}_{K^*}}}^2\\
    \leq&~\parenBig{\ltwo{\tilde{B}^*\breve{{\Delta}}}-T_{n,q}^{(2)}(\delta)\lone{\hat{\Delta}}}^2\\
    \leq&~\parenBig{\ltwo{B^*\hat{\Delta}}-T_{n,q}^{(2)}(\delta)\lone{\hat{\Delta}}}^2\\
    \leq&~\frac{1}{J}\sum_{m=1}^J\parenBig{\ltwo{B^*\hat{\Delta}}-\ltwo{\parenBig{\hat{B}^{\brac{-m}}-B^*}\hat{\Delta}}}^2\\
    \leq&~\frac{1}{J}\sum_{m=1}^J\ltwo{\hat{B}^{\brac{-m}}\hat{\Delta}}^2.
\end{align*}
Under the events $\mc{B}_3,\mc{B}_5,\mc{B}_6$, if $T_{n,q}^{(2)}(\delta)\leq K_\beta$, by Holder's inequality, we have
\begin{align*}
    0\ge&~\dn\sum_{m=1}^J\sum_{i\in\mc{I}_m}\setBig{L\paren{Y_i,\<\hat{B}^{\brac{-m}}\paren{\hat{w}_0+\hat{\Delta}},X_i\>}-L\paren{Y_i,\<\hat{B}^{\brac{-m}}\hat{w}_0,X_i\>}}\\
    \ge&~-\lambda_n^\prime\lone{\hat{\Delta}}+\frac{1}{J}\sum_{m=1}^JC_1\ltwo{\hat{B}^{\brac{-m}}\hat{\Delta}}\parenBig{\ltwo{\hat{B}^{\brac{-m}}\hat{\Delta}}-C_2\sqrt{\frac{\log p}{n_1}}\lone{\hat{B}^{\brac{-m}}\hat{\Delta}}}\\
    \ge&~-\lambda_n^\prime \lone{\hat{\Delta}}+\frac{1}{J}\sum_{m=1}^JC_1\ltwo{\hat{B}^{\brac{-m}}\hat{\Delta}}^2-C_3T_{n,q}^{(3)}(\delta)\sqrt{\frac{\log p}{n_1}}\lone{\hat{\Delta}}^2\\
    \ge&~-2\lambda_n^\prime\abs{\breve{{\Delta}}_{K^*}}+\parenBig{\frac{C_1\tilde{\phi}^2}{4}-4C_3T_{n,q}^{(3)}(\delta)\sqrt{\frac{\log p}{n_1}}}\abs{\breve{{\Delta}}_{K^*}}^2,
\end{align*}
which implies \begin{align*}
1-\sum_{k=K^*}^{K_{ne}}\hat{w}_k=\abs{\breve{{\Delta}}_{K^*}}\leq\frac{16\lambda_n^\prime}{C_1\tilde{\phi}^2}
\end{align*}
provided that $4C_3T_{n,q}^{(3)}(\delta)\sqrt{\frac{\log p}{n_1}}\leq \frac{C_1\tilde{\phi}^2}{8}$. Setting $t_1=\sqrt{\frac{2\log (pK)}{cn_1}}$, $t_2=\sqrt{\frac{\log (pK)}{cn_1}}$ and $\delta=\sqrt{\frac{\log (pK)}{n_{-1}}}$, then $T_{n,q}^{(2)}(\delta)\asymp T_{n,q}\lesssim\Psi_{n,q}$, $T_{n,q}^{(3)}(\delta)\sqrt{\frac{\log p}{n}}\lesssim\Psi_{n,q}$ and
\begin{align*}
    \lambda_n^\prime=C\parenBig{(t_1+1)T_{n,q}^{(2)}(\delta)+t_2}\asymp T_{n,q}.
\end{align*}
Thus, if $\Psi_{n,q}\lesssim\tilde{\phi}^2$, we have
$1-\sum_{k=K^*}^{K_{ne}}\hat{w}_k\leq C_1T_{n,q}/\tilde{\phi}^2$
holds with probability at least $1-C_2(pK)^{-1}-C_2e^{-cn}$ since $J$ is fixed.
\end{proof}
%%%%%%%%%%%%%%%%%%%%%%%%%%%%%%%%%%%%%%%%%%%%%%
\subsection{Proof of \cref{riskconsistency}}
\begin{proof} By \cref{weights2one} and $\E\frac{\partial L(Y_{n+1},\<\beta^*,X_{n+1})}{\partial \mu}X_{n+1}=0$, we have
\begin{align*}
        &~\E\bracBig{L\paren{Y_{n+1},\<\hat{B}\hat{w},X_{n+1}\>}\Big|\set{Y_i,X_i}_{i=1}^n}\\
        \leq&~R_{Bayes}+\absBig{\E\bracBig{L\paren{Y_{n+1},\<\hat{B}\hat{w},X_{n+1}\>}-L\paren{Y_{n+1},\<\beta^*,X_{n+1}\>}\Big|\set{Y_i,X_i}_{i=1}^n}}\\
        \leq&~R_{Bayes}+\frac{1}{2}L_u\E\bracBig{\absBig{\<\hat{B}\hat{w}-\beta^*,X_{n+1}\>}^2\Big|\set{Y_i,X_i}_{i=1}^n}\\
        =&~R_{Bayes}+O_p\parenBig{\ltwo{\sum_{k=1}^K\hat{w}_k\hat{\beta}_{(k)}-\beta^*}^2}\\
        =&~R_{Bayes}+O_p\parenBig{\frac{s\log p}{n}},
    \end{align*}
    which completes the proof.
\end{proof}
%%%%%%%%%%%%%%%%%%%%%%%%%%%%%%%%%%%%%%%%%%%%%%%%
\subsection{Proof of \cref{inference: bahadur representation}}
\begin{proof} Consider the following expansion:
\begin{align*}
    &~\tilde{\beta}(\hat{w})-\beta^*+J^{-1}\dn\sumn\frac{\partial L(Y_i,\<\beta^*,X_i\>)}{\partial \mu}X_i\\
    =&~(I_p-\hat{W}\dn\sumn\frac{\partial^2 L(Y_i,\<\hat{\beta}_{MA}(\hat{w}),X_i\>)}{\partial \mu^2}X_iX_i^\top)(\hat{\beta}_{MA}(\hat{w})-\beta^*)\\
    &~-(\hat{W}-J^{-1})\dn\sumn\frac{\partial L(Y_i,\<\beta^*,X_i\>)}{\partial \mu}X_i\\
    &~+\hat{W}\parenBig{\dn\sumn\frac{\partial^2 L(Y_i,\<\hat{\beta}_{MA}(\hat{w}),X_i\>)}{\partial \mu^2}X_iX_i^\top-Q_n}(\hat{\beta}_{MA}(\hat{w})-\beta^*),
\end{align*}
where
\begin{align*}
    Q_n=\int_0^1\dn\sumn\frac{\partial^2 L(Y_i,\<\beta^*+t(\hat{\beta}_{MA}(\hat{w})-\beta^*),X_i\>)}{\partial \mu^2}X_iX_i^\top\text{d}t.
\end{align*}
It suffices to show that
\begin{equation}\label{proof: bahadur representation}
    \begin{aligned}
    &~\linfBig{\dn\sumn\frac{\partial L(Y_i,\<\beta^*,X_i\>)}{\partial \mu}X_i}=O_p(\sqrt{\frac{\log p}{n}}),\\
    &~\matrixnorm{\dn\sumn\frac{\partial^2 L(Y_i,\<\hat{\beta}_{MA}(\hat{w}),X_i\>)}{\partial \mu^2}X_iX_i^\top-\dn\sumn\frac{\partial^2 L(Y_i,\<\beta^*,X_i\>)}{\partial \mu^2}X_iX_i^\top}_{\max}=O_p(\frac{s\log^2 (p)}{\sqrt{n}}),\\
    &~\matrixnorm{\dn\sumn\frac{\partial^2 L(Y_i,\<\beta^*,X_i\>)}{\partial \mu^2}X_iX_i^\top-Q_n}_{\max}=O_p(\frac{s\log^2 (p)}{\sqrt{n}}).
\end{aligned}
\end{equation}
In fact, by \cref{inference: hessian}, \cref{lemma: CLIME} and \eqref{eq: estimation error}, we have
{\small\begin{align*}
    &~\linfBig{\tilde{\beta}(\hat{w})-\beta^*+J^{-1}\dn\sumn\frac{\partial L(Y_i,\<\beta^*,X_i\>)}{\partial \mu}X_i}\\
    \leq&~\gamma_n\lone{\hat{\beta}_{MA}(\hat{w})-\beta^*}+\matrixnorm{\hat{W}-J^{-1}}_1\linf{\dn\sumn\frac{\partial L(Y_i,\<\beta^*,X_i\>)}{\partial \mu}X_i}\\
    &~+\matrixnorm{\hat{W}}_1\matrixnorm{\dn\sumn\frac{\partial^2 L(Y_i,\<\hat{\beta}_{MA}(\hat{w}),X_i\>)}{\partial \mu^2}X_iX_i^\top-\dn\sumn\frac{\partial^2 L(Y_i,\<\beta^*,X_i\>)}{\partial \mu^2}X_iX_i^\top}_{\max}\lone{\hat{\beta}_{MA}(\hat{w})-\beta^*}\\
    &~+\matrixnorm{\hat{W}}_1\matrixnorm{\dn\sumn\frac{\partial^2 L(Y_i,\<\beta^*,X_i\>)}{\partial \mu^2}X_iX_i^\top-Q_n}_{\max}\lone{\hat{\beta}_{MA}(\hat{w})-\beta^*}\\
    =&~O_p\parenBig{n^{-1/2}\parenBig{s_0\setBig{s\log^{(5-4q^\prime)/(2-2q^\prime)}(p)/\sqrt{n}}^{1-q^\prime}+s^2\log^{5/2}(p)/\sqrt{n}}},
\end{align*}}
which completes the proof.

Since $\E\brac{\frac{\partial L(Y_1,\<\beta^*,X_1\>}{\partial \mu}X_1}=0$, by \cref{lasso:covariate,lasso:lossfunc,supsumsubexp}, we have 
\begin{align*}
    \E\sup_{j\in\brac{p}}\abs{\dn\sumn\frac{\partial L(Y_i,\<\beta^*,X_i\>)}{\partial \mu}X_{i,j}}=O(\sqrt{\frac{\log p}{n}}),
\end{align*} 
which implies $\linf{\dn\sumn\frac{\partial L(Y_i,\<\beta^*,X_i\>)}{\partial \mu}X_i}=O_p(\sqrt{\frac{\log p}{n}}).$ 
Using the same arguments of the proof of \eqref{proof: lemma-CLIME2}, we have
\begin{align*}
&~\matrixnorm{\dn\sumn\frac{\partial^2 L(Y_i,\<\hat{\beta}_{MA}(\hat{w}),X_i\>)}{\partial \mu^2}X_iX_i^\top-\dn\sumn\frac{\partial^2 L(Y_i,\<\beta^*,X_i\>)}{\partial \mu^2}X_iX_i^\top}_{\max}=O_p(\frac{s\log^2 (p)}{\sqrt{n}}),\\
    &~\matrixnorm{\dn\sumn\frac{\partial^2 L(Y_i,\<\beta^*,X_i\>)}{\partial \mu^2}X_iX_i^\top-Q_n}_{\max}=O_p(\frac{s\log^2 (p)}{\sqrt{n}}),
\end{align*}
which implies \eqref{proof: bahadur representation}.
\end{proof}
%%%%%%%%%%%%%%%%%%%%%%%%%%%%%%%%%%%%%%%
\subsection{Proof of \cref{Gaussian approximation}}
\begin{proof} Since $\eta_n=o(\log^{-1/2}(p))$ and $\matrixnorm{S}_{\max}\lesssim\norm{X_1}_{\psi_2}^2$, by \eqref{proof: lemma-CLIME2}, \eqref{proof: lemma-CLIME3}, \cref{lasso:covariate} and \cref{inference: bahadur representation}, we have
\begin{equation}\label{proof:gaussian approximation1}
    \begin{aligned}
     \norm{R_n}_{\mc{G}}\overset{\triangle}{=}&~ \linf{\sqrt{n}\hat{S}(\tilde{\beta}(\hat{w})-\beta^*)+\sqrt{n}SJ^{-1}\dn\sumn\frac{\partial L(Y_i,\<\beta^*,X_i\>)}{\partial \mu}X_i}\\
    \leq&~\matrixnorm{\hat{S}}_{\max}\linf{\sqrt{n}(\tilde{\beta}(\hat{w})-\beta^*)+\sqrt{n}J^{-1}\dn\sumn\frac{\partial L(Y_i,\<\beta^*,X_i\>)}{\partial \mu}X_i}\\
    &~+\matrixnorm{\hat{S}-S}_{\max}\linf{\sqrt{n}J^{-1}\dn\sumn\frac{\partial L(Y_i,\<\beta^*,X_i\>)}{\partial \mu}X_i}\\
    =&~O_p(\eta_n)+O_p(\frac{s\log^{5/2}(p)}{\sqrt{n}})\\
    =&~o_p(\log^{-1/2}(p)),
\end{aligned}
\end{equation}
and then for any $\delta>0$, $\P(\norm{R_n}_{\mc{G}}>\delta\log^{-1/2}(p))\leq \eps/3$ for large enough $p$, where $\eps>0$ is a sufficiently small constant. By anti-concentration inequality \citep[see Corollary 1 and Comment 5 in][]{chernozhukov2015comparison}, we can obtain that
\begin{align}\label{proof:gaussian approximation2}
    \sup_{t\in\R}\P\parenBig{\abs{\norm{\mathbf{Z}}_{\mc{G}}-t}\leq\delta\log^{-1/2}(p)}\leq C\delta\sqrt{\frac{\log p+\log\log p-\log\delta}{\log p}.}
\end{align}
Then there exists a sufficiently small $\delta>0$ such that the right-hand side of \eqref{proof:gaussian approximation1} is less than $\eps/3$ for large enough $p$. Since $\log(np)=o(n^{-1/5})$, combining \eqref{eq: bahadur approximation}, \eqref{proof:gaussian approximation1} and \eqref{proof:gaussian approximation2}, by Lemma 7 in \cite{cai2025statistical}, we have
\begin{align*}
\sup_{t\ge0}\absBig{\P(\norm{\sqrt{n}\hat{S}(\tilde{\beta}(\hat{w})-\beta^*)}_{\mc{G}}\leq t)-\P(\norm{\mathbf{Z}}_{\mc{G}}\leq t)}\leq \eps
\end{align*}
for large enough $p$, which completes the proof.
\end{proof}
%%%%%%%%%%%%%%%%%%%%%%%%%%%%%%%%%%%%%%%%%%%%%%
\subsection{Proof of \cref{bootstrap approximation}}
\begin{proof} Denote $\mathbf{\Sigma}_0=SJ^{-1}\Lambda J^{-1}S$. Conditioning on the data, $\hat{\mathbf{Z}}^{(boot)}$ has the covariance $\hat{\mathbf{\Sigma}}_0\overset{\triangle}{=}\hat{S}\hat{W}\hat{\Lambda}\hat{W}\hat{S}$, where $\hat{\Lambda}=\dn\sumn(\frac{\partial L(Y_i,\<\hat{\beta},X_i\>)}{\partial \mu})^2X_iX_i^\top$. The first step of the proof is to control $\matrixnorm{\hat{\mathbf{\Sigma}}_0-\mathbf{\Sigma}_0}_{\max}$. Note that
\begin{equation}\label{proof:bootstrap-approximation1}
   \begin{aligned}
    \hat{\mathbf{\Sigma}}_0-\mathbf{\Sigma}_0=&~(\hat{S}\hat{W}-SJ^{-1})\hat{\Lambda}(\hat{W}\hat{S}-J^{-1}S)+2(\hat{S}\hat{W}-SJ^{-1})\hat{\Lambda}J^{-1}S\\
    &~+SJ^{-1}\hat{\Lambda}J^{-1}S-SJ^{-1}\Lambda J^{-1}S.
\end{aligned} 
\end{equation}
It suffices to bound the terms $\matrixnorm{\hat{S}\hat{W}-SJ^{-1}}_1$ and $\matrixnorm{\hat{\Lambda}-\Lambda}_{\max}$. By \eqref{proof: lemma-CLIME2}, \eqref{proof: lemma-CLIME3}, \cref{lemma: CLIME} and condition $\eta_n=o(\log^{-3/2}(p))$, we have
\begin{equation}\label{proof:bootstrap-approximation2}
    \begin{aligned}
    \matrixnorm{\hat{S}\hat{W}-SJ^{-1}}_1\leq&~\parenBig{\matrixnorm{\hat{S}-S}_{\max}+\matrixnorm{S}_{\max}}\matrixnorm{\hat{W}-J^{-1}}_1+\matrixnorm{\hat{S}-S}_{\max}\matrixnorm{J^{-1}}_1
   \\ =&~ O_p(s_0\gamma_n^{1-q^\prime}+\frac{s\log^2(p)}{\sqrt{n}})=O_p(\eta_n\log^{-1/2} (p))=o_p(\log^{-2} (p)).
\end{aligned}
\end{equation}
Note that 
\begin{equation}\label{proof:bootstrap-approximation3}
   \begin{aligned}
    \hat{\Lambda}-\Lambda=&~\dn\sumn(\frac{\partial L(Y_i,\<\hat{\beta}_{MA}(\hat{w}),X_i\>)}{\partial \mu}-\frac{\partial L(Y_i,\<\beta^*,X_i\>)}{\partial \mu})^2X_iX_i^\top\\
    &~+2\dn\sumn(\frac{\partial L(Y_i,\<\beta^*,X_i\>)}{\partial \mu})(\frac{\partial L(Y_i,\<\hat{\beta}_{MA}(\hat{w}),X_i\>)}{\partial \mu}-\frac{\partial L(Y_i,\<\beta^*,X_i\>)}{\partial \mu})X_iX_i^\top\\
    &~+\setBig{\dn\sumn(\frac{\partial L(Y_i,\<\beta^*,X_i\>)}{\partial \mu})^2X_iX_i^\top-\E\brac{(\frac{\partial L(Y_i,\<\beta^*,X_i\>)}{\partial \mu})^2X_iX_i^\top}}\\
    \overset{\triangle}{=}&~I_1+I_2+I_3.
\end{aligned} 
\end{equation}
By \cref{inference: hessian}, $\E\brac{\max_{j\in\brac{p}}\abs{X_{1,j}}^4}\lesssim\norm{\max_{j\in\brac{p}}\abs{X_{1j}}}_{\psi_2}^4=O(\log^2(p))$, we have
\begin{align}\label{proof:bootstrap-approximation4}
    \matrixnorm{I_1}_{\max}\leq C^2\lone{\hat{\beta}_{MA}(\hat{w})-\beta^*}^2\dn\sumn\max_{j\in\brac{p}}\abs{X_{i,j}}^4=O_p(\frac{s^2\log^3 (p)}{n}).
\end{align}
By Cauchy-Schwartz inequality, we have
\begin{align*}
    \E\brac{\max_{j\in\brac{p}}\abs{\frac{\partial L(Y_1,\<\beta^*,X_1\>)}{\partial \mu}X_{1,j}^3}}\leq\sqrt{\E\abs{\frac{\partial L(Y_1,\<\beta^*,X_1\>)}{\partial \mu}}^2\E\brac{\max_{j\in\brac{p}}\abs{X_{1,j}}^6}}=O(\log^{3/2}(p)),
\end{align*}
then by \cref{inference: hessian} we can obtain
\begin{align}\label{proof:bootstrap-approximation5}
    \matrixnorm{I_2}_{\max}\leq C\lone{\hat{\beta}_{MA}(\hat{w})-\beta^*}\dn\sumn\max_{j\in\brac{p}}\abs{\frac{\partial L(Y_1,\<\beta^*,X_1\>)}{\partial \mu}X_{i,j}^3}=O_p(\frac{s\log^2(p)}{\sqrt{n}}).
\end{align}
By Corollary 6.5 in \cite{zhang2020concentration}, the $\psi_{1/2}$-norm of quantity $(\frac{\partial L(Y_1,\<\beta^*,X_1\>)}{\partial \mu})^2X_{1,j}X_{1,k}$ is bounded by a constant for any $j,k\in\brac{p}$. By Corollary 6.2 in \cite{zhang2020concentration} and a general concentration inequality \cite[see, for example, Lemma K.3 in the supplementary materials of][]{ning2017general}, we have
\begin{align}\label{proof:bootstrap-approximation6}
    \P\parenBig{\matrixnorm{I_3}_{\max}\ge t}\leq4p^2\exp(-\frac{1}{8}n^{1/5}t^{2/5})+4np^2C_2\exp(-C_2n^{1/5}t^{2/5}),
\end{align}
with $t=C\sqrt{\frac{\log^5(np)}{n}}$ for some sufficiently large constant $C$. Combining \eqref{proof:bootstrap-approximation3}, \eqref{proof:bootstrap-approximation4}, \eqref{proof:bootstrap-approximation5} and \eqref{proof:bootstrap-approximation6}, we have $\matrixnorm{\hat{\Lambda}-\Lambda}_{\max}=O_p(\eta_n\log^{-1/2}(p)+\sqrt{\frac{\log^5(np)}{n}})=o_p(\log^{-2}(p))$ by the condition $\eta_n=o(\log^{-3/2}(p))$ and $\log^5(np)\log^4(p)=o(n).$ Thus, by \eqref{proof:bootstrap-approximation1} and \eqref{proof:bootstrap-approximation2}, we can obtain that
\begin{align*}
    \matrixnorm{\hat{\mathbf{\Sigma}}_0-\mathbf{\Sigma}}_{\max}\leq&~\matrixnorm{\hat{S}\hat{W}-SJ^{-1}}_1^2\matrixnorm{\hat{\Lambda}}_{\max}+2\matrixnorm{\hat{S}\hat{W}-SJ^{-1}}_1\matrixnorm{\hat{\Lambda}}_{\max}\matrixnorm{SJ^{-1}}_1\\
    &~+\matrixnorm{SJ^{-1}}_1^2\matrixnorm{\hat{\Lambda}-\Lambda}_{\max}\\
    =&~o_p(\log^{-2}(p)).
\end{align*}
By Gaussian comparison result \citep[see Proposition 2.1 in ][]{chernozhuokov2022improved}, we have
\begin{align*}
    \sup_{t\ge0}\absBig{\P\parenBig{\norm{\hat{\mathbf{Z}}^{(boot)}}_{\mc{G}}\leq t\Big|\set{(Y_i,X_i)}_{i=1}^n}-\P\parenBig{\norm{\mathbf{Z}}_{\mc{G}}\leq t}}=o_p(1).
\end{align*}
\end{proof}

%%%%%%%%%%%%%%%%%%%%%%%%%%%%%%%%%%%%%%%%%%%%%%%%%%%%%%%%%

\subsection{Proof of \cref{thm:FGMA}}
\begin{proof} Let $\mc{L}(w)=CV(w)/n$. First we show that \begin{align*}
    \frac{A_N^2}{L_N}\brac{\mc{L}(\hat{w}^{(N)})-\mc{L}(\hat{w})}+\frac{1}{2}\ltwo{A_N\hat{w}^{(N)}-(A_N-1)\hat{w}^{(N-1)}-\hat{w}}^2
\end{align*} is a Lyapunov energy function of the number of iterations $N$ for this algorithm, that is, a nonnegative quantity that decreases in each iterations:
\begin{equation}\label{proof:FGMA-lyapunov}
    \begin{aligned}
    &~\frac{A_{N+1}^2}{L_{N+1}}\bracBig{\mc{L}(\hat{w}^{(N+1)})-\mc{L}(\hat{w})}+\frac{1}{2}\ltwo{A_{N+1}\hat{w}^{(N+1)}-(A_{N+1}-1)\hat{w}^{(N)}-\hat{w}}^2\\
    \leq&~\frac{A_N^2}{L_N}\bracBig{\mc{L}(\hat{w}^{(N)})-\mc{L}(\hat{w})}+\frac{1}{2}\ltwo{A_N\hat{w}^{(N)}-(A_N-1)\hat{w}^{(N-1)}-\hat{w}}^2.
\end{aligned}
\end{equation}
Note that for any $w\in\mc{W}^K, l\in\N$, by optimality condition $\<\nabla \mc{L}(z^{(l)})+L_N(p_{L_N}(z^{(l)})-z^{(l)}),w-p_{L_N}(z^{(l)})\>\ge0$, we have
\begin{equation}\label{proof:FGMA-lower-bound}
    \begin{aligned}
    \mc{L}(w)-\mc{L}(p_{L_l}(z^{(l)}))\ge&~\mc{L}(z^{(l)})+\<\nabla \mc{L}(z^{(l)}), w-z^{(l)}\>\\
    &~- \setBig{\mc{L}(z^{(l)})+\<\nabla \mc{L}(z^{(l)}), p_{L_l}(z^{(l)})-z^{(l)}\>+\frac{L_l}{2}\ltwo{p_{L_l}(z^{(l)})-z^{(l)}}^2}\\
    \ge&~L_l\<p_{L_l}(z^{(l)})-z^{(l)}, p_{L_l}(z^{(l)})-w\>-\frac{L_l}{2}\ltwo{p_{L_l}(z^{(l)})-z^{(l)}}^2\\
    \ge&~L_l\<p_{L_l}(z^{(l)})-z^{(l)}, z^{(l)}-w\>+\frac{L_l}{2}\ltwo{p_{L_l}(z^{(l)})-z^{(l)}}^2.
\end{aligned}
\end{equation}
Following \eqref{proof:FGMA-lower-bound}, we have
\begin{align*}
    \mc{L}(\hat{w}^{(N)})-\mc{L}(\hat{w}^{(N+1)})\ge&~ L_{N+1}\<\hat{w}^{(N+1)}-z^{(N+1)}, z^{(N+1)}-\hat{w}^{(N)}\>+\frac{L_{N+1}}{2}\ltwo{\hat{w}^{(N+1)}-z^{(N+1)}}^2,\\
    \mc{L}(\hat{w})-\mc{L}(\hat{w}^{(N+1)})\ge&~ L_{N+1}\<\hat{w}^{(N+1)}-z^{(N+1)}, z^{(N+1)}-\hat{w}\>+\frac{L_{N+1}}{2}\ltwo{\hat{w}^{(N+1)}-z^{(N+1)}}^2.
\end{align*}
Consider a weighted sum of the above inequalities:
\begin{align*}
   &~A_N^2\bracBig{\mc{L}(\hat{w}^{(N)})-\mc{L}(\hat{w}^{(N+1)})}+A_{N+1}\bracBig{\mc{L}(\hat{w})-\mc{L}(\hat{w}^{(N+1)})}\\
    \ge&~\frac{L_{N+1}}{2}\ltwo{A_{N+1}(\hat{w}^{(N+1)}-z^{(N+1)})}^2\\
    &~+L_{N+1}\<A_{N+1}(\hat{w}^{(N+1)}-z^{(N+1)}),A_{N+1}z^{(N+1)}-(A_{N+1}-1)\hat{w}^{(N)}-\hat{w}\>,
\end{align*}
where we use the relation $A_N^2=A_{N+1}^2-A_{N+1}$. By the relation $z^{(N+1)}=\hat{w}^{(N)}+\frac{A_N-1}{A_{N+1}}(\hat{w}^{(N)}-\hat{w}^{(N-1)})$, $L_N\leq L_{N+1}$ and some basic algebra, the above inequality can be reorganized as \eqref{proof:FGMA-lyapunov}. By \eqref{proof:FGMA-lower-bound} again, we have
\begin{align*}
    &~\frac{1}{L_1}\setBig{\mc{L}(\hat{w}^{(1)})-\mc{L}(\hat{w})}+\frac{1}{2}\ltwo{\hat{w}^{(1)}-\hat{w}}^2\\
    \leq&~\frac{1}{2}\ltwo{\hat{w}^{(1)}-\hat{w}}^2-\frac{1}{2}\ltwo{\hat{w}^{(1)}-z^{(1)}}^2-\<\hat{w}^{(1)}-z^{(1)},z^{(1)}-\hat{w}\>\\
    =&~\frac{1}{2}\ltwo{z^{(1)}-\hat{w}}^2.
\end{align*}
By \eqref{proof:FGMA-lyapunov} and induction, we obtain
\begin{align*}
    \mc{L}(\hat{w}^{(N+1)})-\mc{L}(\hat{w})\leq\frac{L_{N+1}}{2A_{N+1}^2}\ltwo{\hat{w}^{(0)}-\hat{w}}^2.
\end{align*}
Note that $A_{N+1}\ge \frac{1}{2}+A_N\ge\cdots\ge\frac{N}{2}+1$ and $L_{N+1}\leq \gamma L_{CV}$. Then
\begin{align*}
    \mc{L}(\hat{w}^{(N+1)})-\mc{L}(\hat{w})\leq\frac{2\gamma L_{CV}}{(N+2)^2}\ltwo{\hat{w}^{(0)}-\hat{w}}^2,
\end{align*}
which completes the proof.
\end{proof}

%%%%%%%%%%%%%%%%%%%%%%%%%%%%%%%%%%%%%%%%%%%

\subsection{Proof of \cref{GMA} in \cref{greedy model averaging}}
\begin{proof} Let $\mc{L}(w)=CV(w)/n$. 
For any $w, w^\prime\in\mc{W}^K$ and $N\in\N$, by Taylor's expansion and \cref{lasso:lossfunc} we have
\begin{align*}
    \mc{L}(w)\leq\mc{L}(w^\prime)+\<\grad_w\mc{L}(w^\prime),w-w^\prime\>+\frac{L_{CV}}{2}\ltwo{w-w^\prime}^2,
\end{align*}
where $L_{CV}^\prime=\inf_{w\in\mc{W}^K}\matrixnorm{\nabla^2\mc{L}(w)}_2$.
Setting $w=\hat{w}^{(N-1)}+\alpha_{N}(e^{(k)}-\hat{w}^{(N-1)}), w^\prime=\hat{w}^{(N-1)}$, we have
\begin{align*}
    &~\mc{L}(\hat{w}^{(N-1)}+\alpha_{N+1}(e^{(k)}-\hat{w}^{(N-1))})\\
    \leq&~\mc{L}(\hat{w}^{(N-1)})+\alpha_{N}\<\grad_w\mc{L}(\hat{w}^{(N-1)}), e^{(k)}-\hat{w}^{(N-1)}\>+\frac{\alpha_{N}^2L_{CV}^\prime}{2}\ltwo{e^{(k)}-\hat{w}^{(N-1)}}^2.
\end{align*}
By the convexity of $\mc{L}(w)$ and $\mc{L}(\hat{w}^{(N)})\leq\sum_{k=1}^K\hat{w}_k\mc{L}(\hat{w}^{(N-1)}+\alpha_{N}(e^{(k)}-\hat{w}^{(N-1)}))$ we have
\begin{align*}
    \mc{L}(\hat{w}^{(N)})\leq &~\mc{L}(\hat{w}^{(N-1)})+\alpha_{N}\<\grad_w\mc{L}(\hat{w}^{(N-1)}), \hat{w}-\hat{w}^{(N-1)}\>+\alpha_{N}^2L^\prime_{CV},
\end{align*}
where we use $\ltwo{e^{(k)}-\hat{w}^{(N-1)}}^2\leq 2$. Then for any $N\in\N$
\begin{align*}
    \mc{L}(\hat{w}^{(N)})-\mc{L}(\hat{w})\leq(1-\alpha_{N})\parenBig{\mc{L}(\hat{w}^{(N-1)})-\mc{L}(\hat{w})}+\alpha_{N}^2L^\prime_{CV}.
\end{align*}
Using the same arguments, we can obtain $\mc{L}(\hat{w}^{(0)})-\mc{L}(\hat{w})\leq L^\prime_{CV}$. By induction, if $\mc{L}(\hat{w}^{(N-1)})-\mc{L}(\hat{w})\leq\frac{4L^\prime_{CV}}{N+3}$, then
\begin{align*}
    \mc{L}(\hat{w}^{(N)})-\mc{L}(\hat{w})\leq\parenBig{1-\frac{2}{N+2}}\parenBig{\frac{4L^\prime_{CV}}{N+3}}+\frac{4L^\prime_{CV}}{(N+2)^2}=\frac{4(N^2+3N+3)}{(N+2)^2(N+3)}L^\prime_{CV}\leq\frac{4L^\prime_{CV}}{N+4},
\end{align*}
where it's easy to verify $(N^2+3N+3)(N+4)\leq(N+2)^2(N+3)$ for any $N\in\N$ and the result holds.

\end{proof}
\section{Some useful lemmas}\label{appendix:lemma}
\begin{lemma}\label{lemma:RSC}(Restricted strong convexity) Denote
\begin{align*}
    \delta\mc{L}\paren{\Delta,\beta}\overset{\triangle}{=}&~L_n(Y,\langle \beta+\Delta,X\rangle)-L_n(Y,\langle\beta,X\rangle)-\langle \grad_\beta L_n(Y,\langle \beta,X\rangle),\Delta\rangle,
\end{align*}
where $L_n(Y,\langle \beta,X\rangle)=\dn\sumn L(Y_i,\langle \beta,X_i\rangle)$. Suppose $\beta_0\in\R^p$ satisfies that $\ltwo{\beta_0}\leq C$ for some constant $C>0$. Under the \cref{lasso:covariate,lasso:lossfunc}, for any fixed $R_1,R_2\ge 0$, if $n\ge c_1$, then the following event
    \begin{align*}
        &~\Big\{\delta\mc{L}\paren{\Delta,\beta}
        \ge C_1\ltwo{\Delta}\parenBig{\ltwo{\Delta}-C_2\sqrt{\frac{\log p}{n}}\lone{\Delta}}\\ 
        &~\qquad\text{ for any }\Delta \in \setBig{\Delta\in\R^p\Big|\ltwo{\Delta}\leq R_1} \text{ and } \beta\in\setBig{\beta\in\R^p\Big|\lone{\beta-\beta_0}\leq \frac{R_2}{1\vee\sqrt{\log p}}}\Big\}
    \end{align*}
    holds with probability at least $1-2\exp\paren{-c_2n}$, 
    where $c_1,c_2,C_1,C_2$ are some positive constants not depending on $n$. 
    
    In addition, if $\inf_{\mu\in \R}\inf_{y\in\mc{Y}}\frac{\partial^2L(y,\mu)}{\partial \mu^2}\ge \eta_l$ for some positive constant $\eta_l>0$ and $n\ge c_1$, then for any fixed $R_1>0$, the following event
    \begin{align*}
        \delta\mc{L}\paren{\Delta,\beta}
        \ge C_1\ltwo{\Delta}\parenBig{\ltwo{\Delta}-C_2\sqrt{\frac{\log p}{n}}\lone{\Delta}}\text{ for any }\Delta \in \setBig{\Delta\in\R^p\Big|\ltwo{\Delta}\leq R_1} \text{ and } \beta\in\R^p
    \end{align*}
    holds with probability at least $1-2\exp\paren{-c_2n}$, where $c_1,c_2,C_1,C_2$ are some positive constants not depending on $n$. 
\end{lemma}
%%%%%%%%%%%%%%%%%%%%%%%%%%%%%%%%%%%%%%%%%%%%%%
\begin{lemma}\label{lemma:regularizer-to-lone}
     Suppose $r_\lambda$ satisfies \cref{con:regularizer}, for two vector $a,b\in\R^p$ and a set $S\in\brac{p}$, we have 
     \begin{align*}
         r_\lambda(a+b)-r_\lambda(a)\ge\lambda\lone{b_{S^c}}-\lambda\lone{b_S}-\frac{\kappa_r}{2}\ltwo{b_{S^c}}^2-2\lambda\lone{a_{S^c}}.
     \end{align*}
\end{lemma}

%%%%%%%%%%%%%%%%%%%%%%%%%%%%%%%%%%%%%%%%%%%%%%

\begin{lemma}\label{centering}
    If $Z$ is a random variable with $\norm{Z}_{\psi_j}\leq K_Z$ for some constant $K_Z$ then we have
    \begin{align*}
        \norm{Z-\E Z}_{\psi_j}\leq C\norm{Z}_{\psi_j},
    \end{align*}
    where $j=2$ if $Z$ is sub-gaussian, $j=1$ if $Z$ is sub-exponential, and $C$ is a universal constant.
\end{lemma}
\begin{lemma}\label{indnorm}
    If $Z_1\in\R^p$ is a random vector with $\norm{Z_1}_{\psi_j}\leq K_Z$ for some constant $K_Z$ and $Z_2\in\R^p$ is a random vector independent of $Z_1$ then we have
    \begin{align*}
        \norm{\frac{\<Z_1,Z_2\>}{\ltwo{Z_2}}}_{\psi_j}\leq K_Z,
    \end{align*}
    where $j=2$ if $Z_1$ is sub-gaussian and $j=1$ if $Z_1$ is sub-exponential.
\end{lemma}
%%%%%%%%%%%%%%%%%%%%%%%%%%%%%%%%%%%%%%%%%%%%%%
\begin{lemma}\label{betarate}
Denote \begin{align*}
    T_{n,q}^{(3)}(\delta)=\begin{cases}
        2K_\beta\sqrt{s_u}&~~\text{if }q=0\\
        K_\beta h_{q,u}^{\frac{q}{2}}(\delta+\sqrt{\frac{\log p_u}{n}})^{-\frac{q}{2}}+(K_\beta+1)h_{q,u}^{\frac{2}{2-q}}&~~\text{if }q\in(0,1]\\
    \end{cases}
\end{align*} and 
    \begin{align*}
        \mc{B}_3=\setBig{\max_{m\in\brac{J}}\max_{k\in\brac{K}}\ltwo{\hat{\beta}_{(k)}^{\brac{-m}}-\beta_{(k)}^*}\leq T_{n,q}^{(2)}(\delta)}\bigcap\setBig{\max_{m\in\brac{J}}\max_{k\in\brac{K}}\lone{\hat{\beta}_{(k)}^{\brac{-m}}-\beta_{(k)}^*}\leq T_{n,q}^{(1)}(\delta)}.
    \end{align*} Under \cref{HD-MA:beta-deviation} and $\mc{B}_3$, if $\delta\in(0,c]$, $T_{n,q}^{(2)}(\delta)\leq K_\beta$, then for any $\lone{\Delta}\leq 2$ we have 
    \begin{align*}
\max_{m\in\brac{J}}\lone{\hat{B}^{\brac{-m}}\Delta}\leq T_{n,q}^{(3)}(\delta)\lone{\Delta},~ \max_{m\in\brac{J}}\ltwo{\hat{B}^{\brac{-m}}\Delta}\leq 2K_\beta\lone{\Delta}\leq 4K_\beta.
\end{align*} and $\mc{B}_3$ holds with probability at least $1-4JK\exp\paren{-n_{-1}\delta^2}$.
\end{lemma}
%%%%%%%%%%%%%%%%%%%%%%%%%%%%%%%%%%%%%%%%%%%%%%
\begin{lemma}\label{VGBound} Let $\mc{H}(d,s)=\set{w\in\set{0,1}^d|~\lzero{w}=s}$, where $d,s\in\Z^+$. Define the Hamming distance $\rho(w,w^\prime)\overset{\triangle}{=}\sum_{j=1}^d\indic{w_j\ne w_j^\prime}$. Suppose $d>\frac{5}{2}s$. There exists a subset $\set{w^{(1)},\dots,w^{(M)}}\subseteq\mc{H}(d,s)$ such that $\rho(w^{(i)},w^{(j)})\ge\frac{s}{2}$ for $1\leq i<j\leq M$ and $\log M\ge\frac{s}{2}\log(\frac{2d}{5s})$.
\end{lemma}
%%%%%%%%%%%%%%%%%%%%%%%%%%%%%%%%%%%%%%%%%%%%%%
\begin{lemma}\label{AOPrate}
    Suppose that $\hat{w}$ is the minimizer of $CV(w)$ over $\mc{W}^K$ and the following conditions hold.
    \begin{enumerate}
        \item $\xi_n^{-1}\sup_{w\in\mc{W}^K}\abs{R^*(w)-R(w)}=O_p(a_n)$ and $a_n\to 0$ as $n\to \infty$
        \item $\xi_n^{-1}\sup_{w\in\mc{W}^K}\abs{CV(w)/n-R^*(w)}=O_p(b_n)$ as $n\to \infty$.
    \end{enumerate}Then we have 
    \begin{align*}
        \frac{R(\hat{w})-R_{\beta^*}}{\inf_{w\in\mc{W}^K}R(w)-R_{\beta^*}}=1+O_p(a_n+b_n).
    \end{align*}
    as $n\to\infty$. Furthermore, if 
    \begin{align*}
    \xi_n^{-1}\sup_{w\in\mc{W}^K}\absBig{R(w)-R^*(w)}
\end{align*}
is uniformly integrable and $b_n\to0$ as $n\to\infty$, then we have
\begin{align*}
    \frac{\Bar{R}(\hat{w})-R_{\beta^*}}{\inf_{w\in\mc{W}^K}\Bar{R}(w)-R_{\beta^*}}\to1
\end{align*}
in probability as $n\to\infty$.
\end{lemma}

%%%%%%%%%%%%%%%%%%%%%%%%%%%%%%%%%%%%%%%%%%%%%%
\begin{lemma}\label{supsumsubexp}
    Suppose $\set{Z_{i,k}}_{i=1}^n$ are independent and identically distributed sub-exponential random variables with mean zero for fixed $k\in\brac{K}$ ($K>1$) and there is a constant $K_Z$ such that $\norm{Z_{i,k}}_{\psi_1}\leq K_Z$ for any $k\in\brac{K}$. If $\sqrt{\frac{\log K}{n}}\to0$ as $n\to\infty$, then we have
    \begin{align*}
        \E\sup_{k\in\brac{K}}\abs{\dn\sumn Z_{i,k}}=O(\sqrt{\frac{\log K}{n}})
    \end{align*}
\end{lemma}

\begin{lemma}\label{lemma: CLIME}
    Under \cref{lasso:lossfunc,lasso:covariate,inference: hessian,inference:smoothness}, $J^{-1}$ is feasible for the optimization problem \eqref{CLIME estimator} with probability going to one and $\matrixnorm{\hat{W}-J^{-1}}_1=O_p(s_0\gamma_n^{1-q^\prime})$.
\end{lemma}

%%%%%%%%%%%%%%%%%%%%%%%%%%%%%%%%%%%%%%%%%%%%%%%%%

\section{Proof of Lemma}\label{appendix:proof-of-lemma}

%%%%%%%%%%%%%%%%%%%%%%%%%%%%%%%%%%%%%%%%%%%%%%%%%%%%%
\subsection{Proof of \cref{lemma:RSC}}
\begin{proof} The proof of this lemma generalizes the proof of Proposition 2 in \cite{negahban2012unified}. In this proof, we use the shorthand $\E_n Z=\dn\sumn Z_i$ for a sequence $\set{Z_i}_{i=1}^n$. We prove the first part of lemma and the second part is similar to the first part.

Denote $\mc{H}=\set{\beta\in\R^p\Big|\lone{\beta-\beta_0}\leq \frac{R_2}{1\vee\sqrt{\log p}}}$, $\gamma_i=\sup_{\beta\in\mc{H}}\abs{\langle\beta-\beta_0,X_i\rangle}$  for each $i\in\brac{n}$.
For fixed $T>0$ and any $\beta\in\mc{H}$, by Taylor's expansion we have
\begin{align*}
    \delta\mc{L}\paren{\Delta,\beta}
    \ge&~\frac{1}{2}\E_n \parenBig{\frac{\partial^2L(Y,\langle\beta+v\Delta,X\rangle)}{\partial \mu^2}\langle X_i,\Delta\rangle^2}\\
    \ge&~\frac{1}{2}L(3T)\E_n\langle \Delta,X\rangle^2\indic{\abs{\langle\Delta,X\rangle}\leq T,\abs{\langle \beta_0,X\rangle}\leq T, \gamma\leq T}\\
    \ge&~\frac{1}{2}L(3T)\E_n\varphi_T\parenBig{\<\Delta,X\>\indic{\abs{\<\beta_0,X\>}\leq T, \gamma\leq T}}
\end{align*}
where $v\in[0,1]$, $L(3T)=\inf_{\abs{\mu}\leq3T}\inf_y\frac{\partial^2L(y,\mu)}{\partial \mu^2}>0$ and
\begin{align*}
    \varphi_T\paren{u}=\begin{cases}
        \abs{u}^2&\text{if $\abs{u}\leq\frac{T}{2}$}\\
        \paren{T-\abs{u}}^2&\text{if $\frac{T}{2}\leq\abs{u}\leq T$}\\
        0&\text{otherwise}
    \end{cases}
\end{align*}
is a Lipschitz function with $\norm{\varphi_T}_{\text{Lip}}\leq T$ and $\vphi_T\paren{0}=0$.

For any $\ltwo{\Delta}=R_1$, by \cref{lasso:covariate} we can obtain that
\begin{align*}
    \E\langle \Delta,X\rangle^2\ge\kappa_lR_1^2.
\end{align*}
By Cauchy-Schwarz inequality and Proposition 2.5.2 of \cite{vershynin2018high}, we have
\begin{align*}
    &~\E\parenBig{\langle\Delta,X\rangle^2-\varphi_T\parenBig{\<\Delta,X\>\indic{\abs{\<\beta_0,X\>}\leq T,\gamma\leq T}}}\\
    \leq&~\E\bracBig{\langle \Delta,X\rangle^2\indic{\abs{\langle\beta_0,X\rangle}>T}+\langle \Delta,X\rangle^2\indic{\gamma>T}+\langle \Delta,X\rangle^2\indic{\abs{\langle\Delta,X\rangle}>\frac{T}{2}}}\\
    \leq&~\sqrt{\E\langle\Delta,X\rangle^4\parenBig{\P\paren{\abs{\langle\beta_0,X\rangle}>T}}}+\sqrt{\E\langle\Delta,X\rangle^4\parenBig{\P\paren{\frac{R_2}{1\vee\sqrt{\log p}}\linf{X}>T}}}\\
    &~+\sqrt{\E\langle\Delta,X\rangle^4\parenBig{\P\paren{\abs{\langle\Delta,X\rangle}>\frac{T}{2}}}}\\
    \leq&~C_1K_x^2R_1^2\bracBig{\exp\paren{-\frac{c_{\psi_2}T^2}{2K_x^2\ltwo{\beta_0}^2}}+\exp\paren{-\frac{c_{\psi_2}T^2\paren{1\vee\log p}}{2K_x^2R_2^2}+\frac{1\vee\log p}{2}}+\exp\paren{-\frac{c_{\psi_2}T^2}{8K_x^2\ltwo{\Delta}^2}}}\\
    \leq&~\frac{\kappa_l}{2}R_1^2
\end{align*}
for any $\ltwo{\Delta}=R_1$, where $T^2=\frac{2K_x^2\paren{\ltwo{\beta_0}^2+4R_1^2+R_2^2}}{c_{\psi_2}}\log\paren{\frac{6C_1K_x^2}{\kappa_l}+1}+\frac{K_x^2R_2^2}{c_{\psi_2}}$ and $C_1,c_{\psi_2}>0$ are universal constants. Thus for any $\ltwo{\Delta}=R_1$, we have
\begin{align*}
    \E\varphi_T\parenBig{\<\Delta,X\>\indic{\abs{\<\beta_0,X\>}\leq T,\gamma\leq T}}\ge\frac{\kappa_l}{2}R_1^2
\end{align*}

For $t>0$, denote 
{\small\begin{align*}
    &~Z_t\paren{X_1,\dots X_i,\dots X_n}\\
    =&~\sup_{\substack{\ltwo{\Delta}= R_1, \\ \lone{\Delta}\leq t}}\absBig{\dn\sumn\varphi_T\parenBig{\<\Delta,X_i\>\indic{\abs{\<\beta_0,X_i\>}\leq T,\gamma_i\leq T}}-\E\varphi_T\parenBig{\<\Delta,X\>\indic{\abs{\<\beta_0,X\>}\leq T,\gamma\leq T}}}.
\end{align*}}
Since $\vphi_T/T$ is a contraction, by a standard symmetrization argument and contraction inequality \citep[see Theorem 4.12 of][]{ledoux1991probability}, we have
\begin{align*}
    \E Z_t&~\leq 2\E\sup_{\ltwo{\Delta}= R_1,\lone{\Delta}\leq t}\absBig{\dn\sumn\varepsilon_i\varphi_T\parenBig{\<\Delta,X_i\>\indic{\abs{\<\beta_0,X_i\>}\leq T,\gamma_i\leq T}}}\\
    &~\leq 4T\E\sup_{\ltwo{\Delta}= R_1,\lone{\Delta}\leq t}\absBig{\dn\sumn\varepsilon_i\<\Delta,X_i\>\indic{\abs{\<\beta_0,X_i\>}\leq T,\gamma_i\leq T}}\\
    &~\leq 4Tt\E\linfBig{\dn\sumn\varepsilon_i X_i\indic{\abs{\<\beta_0,X_i\>}\leq T,\gamma_i\leq T}},
\end{align*}
where $\set{\varepsilon_i}_1^n$ is an i.i.d. sequence of Rademacher variables and is independent of $\set{X_i}_1^n$.

Fixed $j\in\brac{p}$, $\set{V_{ij}}_{i=1}^n\overset{\triangle}{=}\set{\eps_i X_{ij}\indic{\abs{\<\beta_0,X_i\>}\leq T,\gamma_i\leq T}}_{i=1}^n$ are independent mean-zero sub-gaussian variables with $\norm{V_{ij}}_{\psi_2}\leq\norm{X_{ij}}_{\psi_2}\leq K_x$. Thus $\dn\sumn V_{ij}$ is sub-gaussian with $\norm{\dn\sumn V_{ij}}_{\psi_2}\lesssim\frac{K_x}{\sqrt{n}}$ by Proposition 2.6.1 of \cite{vershynin2018high} and then 
\begin{align*}
  \E Z_t\leq 4Tt\E\linfBig{\dn\sumn \eps_i X_i\indic{\abs{\<\beta_0,X_i\>}\leq T,\gamma_i\leq T}}\leq C_2K_xT\sqrt{\frac{\log p}{n}}t
\end{align*}
for some universal constants $C_2>0$. 

Since for each $i\in\brac{n}$ we have $\absBig{Z_t\paren{X_1,\dots X_i,\dots X_n}-Z_t\paren{X_1,\dots X_i^\prime,\dots X_n}}\leq \frac{T^2}{2n}$, where $X_i^\prime$ is an independent copy of $X_i$, by McDiarmid's bounded differences inequality, for any $t>0$ we have
\begin{align*}
    \P\parenBig{Z_t>\E Z_t+\frac{\kappa_l}{4}R_1^2+C_2K_xT\sqrt{\frac{\log p}{n}}t}\leq \exp\parenBig{-\frac{\kappa_l^2R_1^4n}{2T^4}-\frac{8C_2^2K_x^2}{T^2}t^2\log p}.
\end{align*}
Then for any $t>0$ we have
\begin{align*}
    \P\parenBig{Z_t>2C_2K_xT\sqrt{\frac{\log p}{n}}t+\frac{\kappa_l}{4}R_1^2}\leq\exp\parenBig{-\frac{\kappa_l^2R_1^4n}{2T^4}-\frac{8C_2^2K_x^2}{T^2}t^2\log p},
\end{align*}
and then
{\small\begin{align*}
    &~\P\parenBig{\E_n\varphi_T\parenBig{\<\Delta,X\>\indic{\abs{\<\beta_0,X\>}\leq T,\gamma\leq T}}\leq\frac{\kappa_l}{4}R_1^2-2C_2K_xT\sqrt{\frac{\log p}{n}}t\text{ for some }\Delta \in \mc{U}(R_1,t)}\\
    \leq&~\exp\parenBig{-\frac{\kappa_l^2R_1^4n}{2T^4}-\frac{8C_2^2K_x^2}{T^2}t^2\log p},
\end{align*}}
where 
\begin{align*}
    \mc{U}(R_1,t)=\setBig{\Delta\in\R^p\Big|\ltwo{\Delta}=R_1,\lone{\Delta}\leq t}.
\end{align*}
Setting $\mu=\frac{\kappa_lR_1^2}{16C_2K_xT}\sqrt{\frac{n}{\log p}}$, by a peeling argument, if $n\ge\frac{16T^4}{\kappa_l^2R_1^4}\log 2$, we have
{
\begin{equation}
\begin{aligned}
    &~\P\Big(\E_n\varphi_T\parenBig{\<\Delta,X\>\indic{\abs{\<\beta_0,X\>}\leq T,\gamma\leq T}}\leq\frac{\kappa_l}{8}R_1^2-4C_2K_xT\sqrt{\frac{\log p}{n}}\lone{\Delta}\\
    &~\qquad\text{ for some }\Delta \in \mc{U}(R_1,+\infty)\Big)\\
    \leq&~\P\Big(\E_n\varphi_T\parenBig{\<\Delta,X\>\indic{\abs{\<\beta_0,X\>}\leq T,\gamma\leq T}}\leq\frac{\kappa_l}{8}R_1^2-4C_2K_xT\sqrt{\frac{\log p}{n}}\lone{\Delta}\\
    &~\qquad\text{ for some }\Delta \in \mc{U}(R_1,\mu)\Big)\\
    &~+\sum_{l=1}^\infty\P\Big(\E_n\varphi_T\parenBig{\<\Delta,X\>\indic{\abs{\<\beta_0,X\>}\leq T,\gamma\leq T}}\leq\frac{\kappa_l}{8}R_1^2-4C_2K_xT\sqrt{\frac{\log p}{n}}\lone{\Delta}\\
    &~\qquad\text{ for some }\Delta \in \mc{U}(R_1,2^{l}\mu)\backslash\mc{U}(R_1,2^{l-1}\mu)\Big)\\
    \leq&~\P\Big(\E_n\varphi_T\parenBig{\<\Delta,X\>\indic{\abs{\<\beta_0,X\>}\leq T,\gamma\leq T}}\leq\frac{\kappa_l}{4}R_1^2-2C_2K_xT\sqrt{\frac{\log p}{n}}\mu\\
    &~\qquad\text{ for some }\Delta \in \mc{U}(R_1,\mu)\Big)\\
    &~+\sum_{l=1}^\infty\P\Big(\E_n\varphi_T\parenBig{\<\Delta,X\>\indic{\abs{\<\beta_0,X\>}\leq T,\gamma\leq T}}\leq\frac{\kappa_l}{4}R_1^2-2C_2K_xT\sqrt{\frac{\log p}{n}}2^l\mu\\
    &~\qquad\text{ for some }\Delta \in \mc{U}(R_1,2^l\mu)\Big)\\
    \leq&~\sum_{l=0}^\infty\exp\parenBig{-\frac{\kappa_l^2R_1^4n}{2T^4}-\frac{8C_2^2K_x^2}{T^2}\paren{2^l\mu}^2\log p}\\
    \leq&~\parenBig{\exp\parenBig{-\frac{\kappa_l^2R_1^4n}{2T^4}}}\parenBig{\sum_{l=0}^\infty\exp\parenBig{-\frac{16lC_2^2K_x^2}{T^2}\mu^2\log p}}\\
    =&~\parenBig{\exp\parenBig{-\frac{\kappa_l^2R_1^4n}{2T^4}}}\parenBig{\frac{1}{1-\exp\parenBig{-\frac{\kappa_l^2R_1^4n}{16T^4}}}}\\
    \leq&~2\exp\parenBig{-\frac{\kappa_l^2R_1^4n}{2T^4}},
\end{aligned}\label{proof:RSC}
\end{equation}
}
where we use $2^{2l}\ge 2l$ for $l\in\N$.

Recall the construction of the function $\vphi_T\paren{\cdot}$. We have $\vphi_{ca}\paren{cz}=\vphi_a\paren{z}c^2$ and $\vphi_{T_1}\paren{z}\ge\vphi_{T_2}\paren{z}$ if $T_1\ge T_2>0$ for any $a,c>0,z\in\R$. Thus for any $\Delta\in\setBig{\Delta\Big|0<\ltwo{\Delta}\leq R_1}$ we have 
\begin{align*}
    \E_n\vphi_T\parenBig{\<\frac{R_1\Delta}{\ltwo{\Delta}},X\>\indic{\abs{\<\beta_0,X\>}\leq T,\gamma\leq T}}\leq\frac{R_1^2}{\ltwo{\Delta}^2}\E_n\vphi_T\parenBig{\<\Delta,X\>\indic{\abs{\<\beta_0,X\>}\leq T,\gamma\leq T}}.
\end{align*}
Then, by \eqref{proof:RSC}, we have
{
\begin{align*}
    &~\P\Big(\E_n\varphi_T\parenBig{\<\Delta,X\>\indic{\abs{\<\beta_0,X\>}\leq T,\gamma\leq T}}\ge\frac{\ltwo{\Delta}}{R_1^2}\parenBig{\frac{\kappa_lR_1^2}{8}\ltwo{\Delta}-4C_2K_xTR_1\sqrt{\frac{\log p}{n}}\lone{\Delta}}\\
    &~\qquad\text{ for any }\Delta \in \setBig{\Delta\in\R^p\Big|\ltwo{\Delta}\leq R_1}\Big)\\
    =&~\P\Big(\E_n\varphi_T\parenBig{\<\Delta,X\>\indic{\abs{\<\beta_0,X\>}\leq T,\gamma\leq T}}\ge\frac{\ltwo{\Delta}^2}{R_1^2}\parenBig{\frac{\kappa_lR_1^2}{8}-4C_2K_xT R_1\sqrt{\frac{\log p}{n}}\frac{\lone{\Delta}}{\ltwo{\Delta}}}\\
    &~\qquad\text{ for any }\Delta \in \setBig{\Delta\in\R^p\Big|0<\ltwo{\Delta}\leq R_1}\Big)\\
    \ge&~\P\Big(\E_n\vphi_T\parenBig{\<\frac{R_1\Delta}{\ltwo{\Delta}},X\>\indic{\abs{\<\beta_0,X\>}\leq T,\gamma\leq T}}\ge\frac{\kappa_lR_1^2}{8}-4C_2K_xT R_1\sqrt{\frac{\log p}{n}}\frac{\lone{\Delta}}{\ltwo{\Delta}}\\
    &~\qquad\text{ for any }\Delta \in \setBig{\Delta\in\R^p\Big|0<\ltwo{\Delta}\leq R_1}\Big)\\
    =&~\P\Big(\E_n\varphi_T\parenBig{\<\Delta,X\>\indic{\abs{\<\beta_0,X\>}\leq T,\gamma\leq T}}\ge\frac{\kappa_lR_1^2}{8}-4C_2K_xT\sqrt{\frac{\log p}{n}}\lone{\Delta}\\
    &~\qquad\text{ for any }\Delta \in \mc{U}(R_1,+\infty)\Big)\\
    >&~1-2\exp\parenBig{-\frac{\kappa_l^2R_1^4n}{2T^4}},
\end{align*}}
which completes the first part of proof.

If $\inf_{\mu\in \R}\inf_{y\in\mc{Y}}\frac{\partial^2L(y,\mu)}{\partial \mu^2}\ge \eta_l>0$, then for any $\beta\in\R^p$ and $T>0$ we have 
\begin{align*}
    \delta\mc{L}\paren{\Delta,\beta}
    \ge&~\frac{1}{2}\E_n\parenBig{\frac{\partial^2L(Y,\langle\beta+v\Delta,X\rangle)}{\partial \mu^2}\langle X,\Delta\rangle^2}\\
    \ge&~\frac{1}{2}\eta_l\E_n\varphi_T\paren{\<\Delta,X\>}.
\end{align*}
The desired result follows by using the same argument as the previous proof.
\end{proof}
%%%%%%%%%%%%%%%%%%%%%%%%%%%%%%%%%%%%%%%%%%%%%%%%%%%%%%%%%%%%%%
\subsection{Proof of \cref{lemma:regularizer-to-lone}}
\begin{proof} By Lemma 4 in \cite{loh2015regularized}, we have
\begin{align*}
    \abs{\sum_{j\in S}\set{r_\lambda(a_j+b_j)-r_\lambda(a_j)}}\leq\lambda\sum_{j\in S}\absBig{\abs{a_j+b_j}-\abs{a_j}}\leq\lambda\lone{b_S}
\end{align*}
and
\begin{align*}
    r_\lambda(b_{S^c})\ge\lambda\lone{b_{S^c}}-\frac{\kappa_r}{2}\ltwo{b_{S^c}}^2.
\end{align*}
Then, by $r_\lambda(0)=0$, we can obtain that
\begin{align*}
    r_\lambda(a+b)-r_\lambda(a)\ge&~-\lambda\lone{b_S}+r_\lambda(a_{S^c}+b_{S^c})-\lambda\lone{a_{S^c}}\\
    \ge&~r_\lambda(b_{S^c})-\lambda\lone{b_S}-2\lambda\lone{a_{S^c}}\\
    \ge&~\lambda\lone{b_{S^c}}-\lambda\lone{b_S}-2\lambda\lone{a_{S^c}}-\frac{\kappa_r}{2}\ltwo{b_{S^c}}^2,
\end{align*}
which completes the proof.
\end{proof}

%%%%%%%%%%%%%%%%%%%%%%%%%%%%%%%%%%%%%%%%%
\subsection{Proof of \cref{centering}}
\begin{proof} If $Z$ is sub-gaussian, then the result follows from Lemma 2.6.8. in \cite{vershynin2018high}. If $Z$ is sub-exponential, then we have
\begin{align*}
    &~\norm{Z-\E Z}_{\psi_1}\\
    \leq&~\norm{Z}_{\psi_1}+\norm{\E Z}_{\psi_1}\\
    \leq&~\norm{Z}_{\psi_1}+\norm{1}_{\psi_1}\E\abs{Z}\\
    \leq&~(1+C\norm{1}_{\psi_1})\norm{Z}_{\psi_1},
\end{align*}
where the last inequality follows from Proposition 2.7.1 of \cite{vershynin2018high} and $C$ is a universal constant. The proof is complete.
\end{proof}
%%%%%%%%%%%%%%%%%%%%%%%%%%%%%%%%%%%%%%%%%%%
\subsection{Proof of \cref{indnorm}}
\begin{proof} We only prove the result for the case $j=1$ and the result for the case $j=2$ is similar. Suppose $\norm{Z_1}_{\psi_1}\leq K_Z$. Conditioning on $Z_2$, we have
\begin{align*}
   \E \bracBig{\exp\paren{\frac{\abs{\<Z_1,Z_2\>}}{K_Z\ltwo{Z_2}}}\Big|Z_2}\leq 2.
\end{align*}
By taking expectation and the definition of sub-exponential random variable, the result holds.
\end{proof}

%%%%%%%%%%%%%%%%%%%%%%%%%%%%%%%%%%%%%%%%

\subsection{Proof of \cref{betarate}}
\begin{proof} For $q\in(0,1]$ and any $k\in\brac{K}$, let $S_{k,q}=\set{j\in\brac{p}|~\abs{\beta^*_{(k),j}}>h_{q,u}^{-q/(2-q)}}$. Since $h_{q,u}^{-q^2/(2-q)}\abs{S_{k,q}}\leq h_{q,u}^q$ and $\lone{\beta^*_{(k),S_{k,q}^c}}h_{q,u}^{\frac{q}{2-q}}\leq \norm{\beta^*_{(k),S_{k,q}^c}}_qh_{q,u}^{\frac{q^2}{2-q}}\leq h_{q,u}^{\frac{2q}{2-q}}$, then $\lone{\beta^*_{(k)}}\leq \sqrt{\abs{S_{k,q}}}\ltwo{\beta^*_{(k)}}+\lone{\beta^*_{(k),S_{k,q}^c}}\leq (K_\beta+1)h_{q,u}^{\frac{q}{2-q}}$.
Note that under the event $\mc{B}_3$ and $T_{n,q}^{(2)}(\delta)\leq K_\beta$, for any $\lone{\Delta}\leq 2$, we have 
\begin{align*}
\max_{m\in\brac{J}}\lone{\hat{B}^{\brac{-m}}\Delta}&\leq\lone{\Delta}\max_{m\in\brac{J}}\max_{k\in\brac{K}}\parenBig{\lone{\hat{\beta}_{(k)}^{\brac{-m}}-\beta_{(k)}^*}+\lone{\beta^*_{(k)}}}\leq T_{n,q}^{(3)}(\delta)\lone{\Delta}\\
\max_{m\in\brac{J}}\ltwo{\hat{B}^{\brac{-m}}\Delta}&\leq\lone{\Delta}\max_{m\in\brac{J}}\max_{k\in\brac{K}}\parenBig{\ltwo{\hat{\beta}_{(k)}^{\brac{-m}}-\beta_{(k)}^*}+\ltwo{\beta^*_{(k)}}}\leq 2K_\beta\lone{\Delta}\leq 4K_\beta.
\end{align*}
By \cref{HD-MA:beta-deviation}, since $J$ is fixed, with probability at least $1-4JK\exp\paren{-n_{-1}\delta^2}$ we have
\begin{align*}
\max_{m\in\brac{J}}\max_{k\in\brac{K}}\ltwo{\hat{\beta}_{(k)}^{\brac{-m}}-\beta_{(k)}^*}\leq T_{n,q}^{(2)}(\delta),~ \max_{m\in\brac{J}}\max_{k\in\brac{K}}\lone{\hat{\beta}_{(k)}^{\brac{-m}}-\beta_{(k)}^*}\leq T_{n,q}^{(1)}(\delta),
\end{align*} which completes the proof.
\end{proof}
%%%%%%%%%%%%%%%%%%%%%%%%%%%%%%%%%%%%%%%%%
\subsection{Proof of \cref{VGBound}}
\begin{proof} This proof is similar to the proof of Lemma 4 in \cite{raskutti2011minimax}. Note that $\abs{\mc{H}(d,s)}=\binom{d}{s}$. For fixed $w$, we have
\begin{align*}
    \absBig{\setBig{w^\prime\in\set{0,1}^d\Big|\lzero{w^\prime}=s, \rho(w,w^\prime)\leq\frac{s}{2}}}\leq\binom{d}{\lfloor\frac{s}{2}\rfloor}2^{\lfloor\frac{s}{2}\rfloor},
\end{align*}
since $\rho(w,w^\prime)\leq\frac{s}{2}$ implies that $w$ and $w^\prime$ have at most $\lfloor\frac{s}{2}\rfloor$ different coordinates. Then for any fixed set $\mc{A}\subseteq\mc{H}(d,s)$ we have
\begin{align*}
    \absBig{\setBig{w^\prime\in\set{0,1}^d\Big|\lzero{w^\prime}=s, \rho(w,w^\prime)\leq\frac{s}{2}~~~\text{ for some } w\in\mc{A}}}\leq\abs{\mc{A}}\binom{d}{\lfloor\frac{s}{2}\rfloor}2^{\lfloor\frac{s}{2}\rfloor}.
\end{align*}
Starting with $\mc{A}=0$, we can inducively add new point $w^\prime\in\mc{H}$ to $\mc{A}$ such that $\rho(w,w^\prime)>\frac{s}{2}$ for any $w\in\mc{A}$. If $\abs{\mc{A}}\binom{d}{\lfloor\frac{s}{2}\rfloor}2^{\lfloor\frac{s}{2}\rfloor}<\abs{\mc{H}}=\binom{d}{s}$, then
\begin{align*}
    \absBig{\setBig{w^\prime\in\set{0,1}^d\Big|\lzero{w^\prime}=s, \rho(w,w^\prime)\leq\frac{s}{2}~~~\text{ for some } w\in\mc{A}}}<\abs{\mc{H}}
\end{align*}
which means we can find new point $w^\prime\in\mc{H}$ such that $\rho(w^\prime,w)>\frac{s}{2}$ for any $w\in\mc{A}$. Hence, there exists $\mc{A}=\set{w^{(1)},\dots,w^{(M)}}\in\mc{H}$ such that $\rho(w^{(i)},w^{(j)})>\frac{s}{2}$ for any $i,j\in\brac{M}$ and 
\begin{align*}
    \abs{\mc{A}}=&~M\ge\frac{\binom{d}{s}}{\binom{d}{\lfloor\frac{s}{2}\rfloor}2^{\lfloor\frac{s}{2}\rfloor}}=\prod_{j=1}^{s-\lfloor\frac{s}{2}\rfloor}\parenBig{\frac{d-\lfloor\frac{s}{2}\rfloor+1-j}{\lfloor\frac{s}{2}\rfloor+j}}2^{-\lfloor\frac{s}{2}\rfloor}\\
    \ge&~\parenBig{\frac{d-\lfloor\frac{s}{2}\rfloor}{s}}^{s-\lfloor\frac{s}{2}\rfloor}2^{-\lfloor\frac{s}{2}\rfloor}\ge\parenBig{\frac{d-\frac{s}{2}}{2s}}^{\frac{s}{2}}\ge\parenBig{\frac{2d}{5s}}^{\frac{s}{2}}
\end{align*}
since $d>\frac{5}{2}s$.
\end{proof}
%%%%%%%%%%%%%%%%%%%%%%%%%%%%%%%%%%%%%%%%%%%%%%%%%%
\subsection{Proof of \cref{AOPrate}}
\begin{proof} Note that
\begin{align*}
\inf_{w\in\mc{W}^K}R(w)-\inf_{w\in\mc{W}^K}R^*(w)=&~\inf_{w\in\mc{W}^K}\parenBig{R(w)-\inf_{w^\prime\in\mc{W}^K}R^*(w^\prime)}\ge
    \inf_{w\in\mc{W}^K}\parenBig{R(w)-R^*(w)}\\
    \ge&~-\sup_{w\in\mc{W}^K}\absBig{R(w)-R^*(w)},\\
    \inf_{w\in\mc{W}^K}R(w)-\inf_{w\in\mc{W}^K}R^*(w)=&~\sup_{w^\prime\in\mc{W}^K}\parenBig{\inf_{w\in\mc{W}^K}R(w)-R^*(w^\prime)}\leq \sup_{w^\prime\in\mc{W}^K}\parenBig{R(w^\prime)-R^*(w^\prime)}\\
    \leq&~\sup_{w\in\mc{W}^K}\absBig{R(w)-R^*(w)}.
\end{align*}
Then 
\begin{align*}
    \absBig{\inf_{w\in\mc{W}^K}R(w)-\inf_{w\in\mc{W}^K}R^*(w)}\leq\sup_{w\in\mc{W}^K}\absBig{R(w)-R^*(w)}
\end{align*} and similarly we have 
\begin{align*}
    \absBig{\inf_{w\in\mc{W}^K}CV(w)/n-\inf_{w\in\mc{W}^K}R^*(w)}\leq\sup_{w\in\mc{W}^K}\absBig{CV(w)/n-R^*(w)}.
\end{align*}
Thus,
\begin{align*}
    \absBig{\frac{\inf_{w\in\mc{W}^K}R(w)-R_{Bayes}}{\inf_{w\in\mc{W}^K}R^*(w)-R_{Bayes}}-1}=&~\xi_n^{-1}\absBig{\inf_{w\in\mc{W}^K}R(w)-\inf_{w\in\mc{W}^K}R^*(w)}\\
    \leq&~\xi_n^{-1}\sup_{w\in\mc{W}^K}\absBig{R(w)-R^*(w)}\\
    =&~o_p(1).
\end{align*}
Since $\hat{w}$ is the minimizer of $CV(w)/n$ over $\mc{W}^K$, we can obtain that
\begin{align*}
    &~\frac{R(\hat{w})-R_{Bayes}}{\inf_{w\in\mc{W}^K}R(w)-R_{Bayes}}-1\\
    =&~\frac{\inf_{w\in\mc{W}^K}R^*(w)-R_{Bayes}}{\inf_{w\in\mc{W}^K}R(w)-R_{Bayes}}\times {\frac{R(\hat{w})-\inf_{w\in\mc{W}^K}R(w)}{\inf_{w\in\mc{W}^K}R^*(w)-R_{Bayes}}}\\
    \leq&~\absBig{\frac{\inf_{w\in\mc{W}^K}R^*(w)-R_{Bayes}}{\inf_{w\in\mc{W}^K}R(w)-R_{Bayes}}}\xi_n^{-1}\bracBig{2\sup_{w\in\mc{
    W}^K}\abs{R(w)-R^*(w)}+2\sup_{w\in\mc{
    W}^K}\abs{CV(w)/n-R^*(w)}}\\
    =&~O_p(a_n+b_n),
\end{align*}
which completes the first part of proof.

Now, we assume that $\xi_n^{-1}\sup_{w\in\mc{W}^K}\absBig{R(w)-R^*(w)}$ is uniformly integrable. Then we have
\begin{align*}
    \xi_n^{-1}\sup_{w\in\mc{W}^K}\abs{\Bar{R}(w)-R^*(w)}\leq\E\parenBig{\xi_n^{-1}\sup_{w\in\mc{W}^K}\abs{R(w)-R^*(w)}}=o(1)
\end{align*}
since $a_n\to 0$ as $n\to\infty$. Similarly to the previous proof, we have
\begin{align*}
    \absBig{\frac{\inf_{w\in\mc{W}^K}\Bar{R}(w)-R_{Bayes}}{\inf_{w\in\mc{W}^K}R^*(w)-R_{Bayes}}-1}=&~\xi_n^{-1}\absBig{\inf_{w\in\mc{W}^K}\Bar{R}(w)-\inf_{w\in\mc{W}^K}R^*(w)}\\
    \leq&~\xi_n^{-1}\sup_{w\in\mc{W}^K}\absBig{\Bar{R}(w)-R^*(w)}\\
    =&~o(1).
\end{align*}Consequently,
\begin{align*}
    &~\frac{\Bar{R}(\hat{w})-R_{Bayes}}{\inf_{w\in\mc{W}^K}\Bar{R}(w)-R_{Bayes}}-1\\
    \leq&~\absBig{\frac{\inf_{w\in\mc{W}^K}R^*(w)-R_{Bayes}}{\inf_{w\in\mc{W}^K}\Bar{R}(w)-R_{Bayes}}}\xi_n^{-1}\bracBig{2\sup_{w\in\mc{
    W}^K}\abs{\Bar{R}(w)-R^*(w)}+2\sup_{w\in\mc{
    W}^K}\abs{CV(w)/n-R^*(w)}}\\
    =&~o_p(1),
\end{align*}
which completes the proof.
\end{proof}
%%%%%%%%%%%%%%%%%%%%%%%%%%%%%%%%%%%%%%%%%%%%%%%%%%%%

\subsection{Proof of \cref{supsumsubexp}}
\begin{proof} By Bernstein's inequality, for any $t>0$ we have
\begin{align*}
    &~\P\parenBig{\sup_{k\in\brac{K}}\abs{\dn\sumn Z_{i,k}}>t}\leq\sum_{k\in\brac{K}}\P\parenBig{\abs{\dn\sumn Z_{i,k}}>t}
    \leq 2K\exp\parenBig{-c\min(\frac{t^2}{K_Z^2},\frac{t}{K_Z})n},
\end{align*}
where $c$ is a universal constant. Then
\begin{align*}
   &~\E\sup_{k\in\brac{K}}\abs{\dn\sumn Z_{i,k}}\\
    =&~\int_0^\infty\P\parenBig{\sup_{k\in\brac{K}}\abs{\dn\sumn Z_{i,k}}>t}\text{d}t\\
    \leq&~\int_0^{K_Z\sqrt{\frac{\log K}{n}}}1\text{d}t+\int_{K_Z\sqrt{\frac{\log K}{n}}}^\infty 2K\exp\parenBig{-c\frac{t^2}{K_Z^2}n}\text{d}t+\int_{K_Z}^\infty2K\exp\parenBig{-c\frac{t}{K_Z}n}\text{d}t\\
    =&~O(\sqrt{\frac{\log K}{n}})+O(\sqrt{\frac{1}{n\log K}})+O(\frac{1}{n})\\
    =&~O(\sqrt{\frac{\log K}{n}}),
\end{align*}
where we use $\int_a^\infty\exp\paren{-cnt^2}\text{d}t\leq\int_0^\infty\exp\paren{-2cnat}\exp\paren{-cna^2}\text{d}t=\frac{\exp\paren{-cna^2}}{2cna}$ for any $a,c>0$ and $\sqrt{\frac{\log K}{n}}=o(1)$.
\end{proof}
%%%%%%%%%%%%%%%%%%%%%%%%%%%%%%%%%%%%%%%%%%%%%%%%%%%%%
\subsection{Proof of \cref{lemma: CLIME}}
\begin{proof} We first show that $J^{-1}$ is feasible for the optimization problem \eqref{CLIME estimator} with probablity going to one. By \cref{inference: hessian}, we have
\begin{equation}\label{proof: lemma-CLIME1}
    \begin{aligned}
    &~\matrixnorm{J\dn\sumn\frac{\partial^2 L(Y_i,\<\hat{\beta}_{MA}(\hat{w}),X_i\>)}{\partial \mu^2}X_iX_i^\top-I_p}_{\max}\\
    \leq&~\matrixnorm{J^{-1}}_1~\matrixnorm{\dn\sumn\frac{\partial^2 L(Y_i,\<\hat{\beta}_{MA}(\hat{w}),X_i\>)}{\partial \mu^2}X_iX_i^\top-J}_{\max}\\
    \lesssim&~\matrixnorm{\dn\sumn\frac{\partial^2 L(Y_i,\<\hat{\beta}_{MA}(\hat{w}),X_i\>)}{\partial \mu^2}X_iX_i^\top-\dn\sumn\frac{\partial^2 L(Y_i,\<\beta^*,X_i\>)}{\partial \mu^2}X_iX_i^\top}_{\max}\\
    &~+\matrixnorm{\dn\sumn\frac{\partial^2 L(Y_i,\<\beta^*,X_i\>)}{\partial \mu^2}X_iX_i^\top-J}_{\max}.
\end{aligned}
\end{equation}
By \cref{inference:smoothness}, we can obtain that
\begin{align*}
    &~\matrixnorm{\dn\sumn\frac{\partial^2 L(Y_i,\<\hat{\beta}_{MA}(\hat{w}),X_i\>)}{\partial \mu^2}X_iX_i^\top-\dn\sumn\frac{\partial^2 L(Y_i,\<\beta^*,X_i\>)}{\partial \mu^2}X_iX_i^\top}_{\max}\\
    \leq&~ C\lone{\hat{\beta}_{MA}(\hat{w})-\beta^*}\max_{1\leq j,k,l\leq p}\abs{\dn\sumn X_{i,j}X_{i,k}X_{i,l}}.
\end{align*}
By a basic Orlicz norm inequality \citep[see, for example, Lemma 2 in the supplementary material of ][]{cai2025statistical}, we have
\begin{align*}
    \parenBig{\E(\max_{1\leq j\leq p}\abs{X_{1j}}^3)}^{1/3}\lesssim\sqrt{\log p}\norm{X_1}_{\psi_2}.
\end{align*}
Combined with \eqref{eq: estimation error}, we have
\begin{align}\label{proof: lemma-CLIME2}
    \matrixnorm{\dn\sumn\frac{\partial^2 L(Y_i,\<\hat{\beta}_{MA}(\hat{w}),X_i\>)}{\partial \mu^2}X_iX_i^\top-\dn\sumn\frac{\partial^2 L(Y_i,\<\beta^*,X_i\>)}{\partial \mu^2}X_iX_i^\top}_{\max}=O_p(\frac{s\log^2(p)}{\sqrt{n}}).
\end{align}
By \cref{lasso:lossfunc,lasso:covariate,supsumsubexp}, we have
\begin{align}\label{proof: lemma-CLIME3}
\E\matrixnorm{\dn\sumn\frac{\partial^2 L(Y_i,\<\beta^*,X_i\>)}{\partial \mu^2}X_iX_i^\top-J}_{\max}=O(\sqrt{\frac{\log p}{n}}).
\end{align}
Combining \eqref{proof: lemma-CLIME1}, \eqref{proof: lemma-CLIME2} and \eqref{proof: lemma-CLIME3}, we have
\begin{align*}
    \matrixnorm{J\dn\sumn\frac{\partial^2 L(Y_i,\<\hat{\beta}_{MA}(\hat{w}),X_i\>)}{\partial \mu^2}X_iX_i^\top-I_p}_{\max}=O_p(\frac{s\log^2(p)}{\sqrt{n}}),
\end{align*}
which implies that $J^{-1}$ is feasible and $\matrixnorm{\hat{W}}_1\leq\matrixnorm{J^{-1}}_1$ with probability going to one.

Next, we are going to show that $\matrixnorm{\hat{W}-J^{-1}}_1\leq O_p(s_0\gamma_n^{1-q^\prime})$. By \eqref{proof: lemma-CLIME2}, \eqref{proof: lemma-CLIME3} and \cref{inference: hessian}, we have
\begin{align*}
   &~\matrixnorm{\hat{W}-J^{-1}}_{\max}\\
    \leq&~\matrixnorm{J^{-1}}_1\matrixnorm{\hat{W}J-I_p}_{\max}\\
    \leq&~\matrixnorm{J^{-1}}_1\Big(\matrixnorm{\hat{W}\dn\sumn\frac{\partial^2 L(Y_i,\<\hat{\beta}_{MA}(\hat{w}),X_i\>)}{\partial \mu^2}X_iX_i^\top-I_p}_{\max}\\
    &~+\matrixnorm{\hat{W}}_1\matrixnorm{\dn\sumn\frac{\partial^2 L(Y_i,\<\hat{\beta}_{MA}(\hat{w}),X_i\>)}{\partial \mu^2}X_iX_i^\top-J}_{\max}\Big)\\
=&~O_p(\gamma_n).
\end{align*}
Denote $\hat{W}=(\hat{\mathbf{b}}_1,\dots,\hat{\mathbf{b}}_p)^\top$. Since $J^{-1}$ is feasible with probability going to one, we can assume that $J^{-1}$ is feasible in the rest of the proof. Then $\lone{\hat{\mathbf{b}}_j}\leq\lone{\tilde{\mathbf{b}}_j}$ for each $j\in\brac{p}$. Consider the following two cases:
\begin{enumerate}
    \item Suppose $q^\prime\in(0,1)$. Let $S_j=\set{i\in\brac{p}|~\abs{\tilde{b}_{j,i}}>\gamma_n}$. By \cref{inference: hessian}, we have $\abs{S_j}\gamma_n^{q^\prime}\leq s_0$ and $\sum_{i\in S_j^c}\abs{\tilde{b}_{j,i}}\leq\gamma_n^{1-q^\prime}\sum_{i\in \brac{p}}\abs{\tilde{b}_{j,i}}^{q^\prime}\leq s_0\gamma_n^{1-q^\prime}$ for each $j\in\brac{p}.$ By $\lone{\hat{\mathbf{b}}_j}\leq\lone{\tilde{\mathbf{b}}_j}$ for each $j\in\brac{p}$, we have
    \begin{equation}\label{proof: lemma-CLIME4}
        \begin{aligned}
        \matrixnorm{\hat{W}-J^{-1}}_1=&~\max_{j\in\brac{p}}\lone{\hat{\mathbf{b}}_j-\tilde{\mathbf{b}}_j}\\
        \leq&~\max_{j\in\brac{p}}\lone{(\hat{\mathbf{b}}_j-\tilde{\mathbf{b}}_j)_{S_j}}+\max_{j\in\brac{p}}\lone{(\hat{\mathbf{b}}_j-\tilde{\mathbf{b}}_j)_{S_j^c}}\\
        \leq&~\matrixnorm{\hat{W}-J^{-1}}_{\max}\max_{j\in\brac{p}}\abs{S_j}+\max_{j\in\brac{p}}\parenBig{\lone{(\tilde{\mathbf{b}}_j)_{S_j^c}}+\lone{\tilde{\mathbf{b}}_j}-\lone{(\hat{\mathbf{b}}_j)_{S_j}}}\\
        \leq&~\matrixnorm{\hat{W}-J^{-1}}_{\max}\max_{j\in\brac{p}}\abs{S_j}+\max_{j\in\brac{p}}\parenBig{2\lone{(\tilde{\mathbf{b}}_j)_{S_j^c}}+\lone{(\hat{\mathbf{b}}_j-\tilde{\mathbf{b}}_j)_{S_j}}}\\
        \leq&~2\matrixnorm{\hat{W}-J^{-1}}_{\max}\max_{j\in\brac{p}}\abs{S_j}+2\max_{j\in\brac{p}}\lone{(\tilde{\mathbf{b}}_j)_{S_j^c}}\\
        =&~O_p(s_0\gamma_n^{1-q^\prime}).
    \end{aligned}
    \end{equation}
    \item Suppose $q^\prime=0$. Let $S_j=\set{i\in\brac{p}|~\abs{\tilde{b}_{j,i}}>0}$. Then we have $\abs{S_j}\leq s_0$ and $\lone{(\tilde{\mathbf{b}}_j)_{S_j^c}}=0$. Using the same arguments as \eqref{proof: lemma-CLIME4}, we have $\matrixnorm{\hat{W}-J^{-1}}_1=O_p(s_0\gamma_n).$
\end{enumerate}
\end{proof}

\end{document}